\documentclass{article}

\usepackage[left=3cm,right=2.5cm,top=3cm,bottom=2.5cm]{geometry}

\usepackage{authblk}

\usepackage[utf8]{inputenc}
\usepackage[english]{babel}
\usepackage{amsthm}
\usepackage{amsmath}
\usepackage{amssymb}
\usepackage{amsfonts}
\usepackage{dsfont}
\usepackage{mathrsfs}

\DeclareMathOperator*{\esssup}{ess\,sup}
\DeclareMathOperator*{\essinf}{ess\,inf}

\usepackage[colorlinks=true,%
  linkcolor=black,anchorcolor=black,%
  citecolor=black,filecolor=black,%
  menucolor=black,%
  urlcolor=black,bookmarks]{hyperref}

\usepackage{enumerate}
\usepackage{indentfirst}
\usepackage{color}
\usepackage{pgfplots}
\pgfplotsset{width=9cm,compat=1.15}

\allowdisplaybreaks

\newcommand{\norm}[2]{\left |  {#1}  \right| _{#2}}
\newcommand{\norma}[2]{\left \|  {#1} \right \|_{#2}}
\newcommand{\prodl}[2]{\left ( {#1} , {#2} \right )}

\newcommand{\D}[1]{\displaystyle{#1}}

\renewcommand{\qedsymbol}{$\blacksquare$}

\usepackage{verbatim}

\newtheorem{teo}{Theorem}
\newtheorem{lema}[teo]{Lemma}
\newtheorem{coro}[teo]{Corollary}
\newtheorem{defi}[teo]{Definition}

\theoremstyle{remark}

\newtheorem{obsv}[teo]{\bf Remark}
\newtheorem{hyp}{\bf Hypothesis}

\newcommand\blfootnote[1]{%
  \begingroup
  \renewcommand\thefootnote{}\footnote{#1}%
  \addtocounter{footnote}{-1}%
  \endgroup
}

\begin{document}
\title{Uniform in time solutions for a chemotaxis with potential consumption model}

\date{}

\author[*]{André Luiz Corrêa Vianna Filho}
\author[*]{Francisco Guillén-González}
\affil[*]{\footnotesize Departament of Differential Equations and Numerical Analysis, Universidad de Sevilla}

\maketitle
\blfootnote{Email addresses: \texttt{acorreaviannafilho@us.es} (A. L. Corrêa Vianna Filho), \texttt{guillen@us.es} (F. Guillén-González).}
\begin{abstract}
{
In this work we investigate the following chemo-attraction with consumption model in bounded domains of \, $\mathbb{R}^N$ ($N=1,2,3$):
 $$ \partial_t u - \Delta u  = - \nabla \cdot (u \nabla v), \quad
        \partial_t v - \Delta v  = - u^s v
        $$
where $s\ge 1$, endowed with isolated boundary conditions and initial conditions for $(u,v)$.   
The main novelty in the model is  the nonlinear potential consumption term $u^sv$.
 Through the convergence of solutions of an adequate truncated model, two main results are established; existence of uniform in time weak solutions in $3D$ domains, and uniqueness and regularity in $2D$ (or $1D$) domains. Both results are proved imposing minimal regularity assumptions on the boundary of the domain.
}

\vspace{12pt}



\end{abstract}


\section{Introduction}
    \label{introduction}
    
    Chemotaxis is the movement of cells in response to the concentration gradient of a chemical signal. Cells can be attracted or repelled by the chemical substance. This phenomenon is present in many physiological events, such as wound healing, immune cells migration, among others. Chemotaxis also plays an important role in undesired events, such as inflammatory diseases and cancer metastasis \cite{murphy2001chemokines} \cite{wang2011signaling}. Besides, there are also studies which call the attention to the effects of chemotaxis in the migration of bacteria towards polluting substances and its usefulness in the degradation of these substances \cite{pandey2002bacterial} \cite{parales2000toluene}. All this makes chemotaxis a matter of very practical interest.
    
    \
    
    In the present work, we are interested in a model where the cells are attracted by a chemical substance that they consume. Let $\Omega$ be a  bounded domain of $\mathbf{R}^N$ ($N=1,2,3$) and let $\Gamma$ be its boundary. Let $u=u(t,x)$ and $v=v(t,x)$ be the density of cell population and the concentration of chemical substance, respectively, on $x \in \Omega$ and $t > 0$. This model is governed by the initial-boundary PDE problem
    \begin{equation}\tag{P}
      \left\{\begin{array}{l}
        \partial_t u - \Delta u  = - \nabla \cdot (u \nabla v), \quad
        \partial_t v - \Delta v  = - u^s v, \\
        \partial_\eta u \Big |_{\Gamma}  =  \partial_\eta v \Big |_{\Gamma} = 0, \quad
        u(0)  = u^0, \quad v(0) = v^0,
      \end{array}\right.
      \label{problema_P} 
    \end{equation}
    where $\nabla \cdot (u \nabla v)$ is the chemotaxis term and $u^s$ is the consumption rate, with $s \geq 1$. $\partial_\eta u$ denotes the normal derivative of $u$ on the boundary. The initial conditions $u^0 \geq 0$ and $v^0 \geq 0$ in $\Omega$.
    
    Theoretical properties of this model have been studied in different contexts when $s = 1$. Existence of global weak solutions which become smooth after a sufficiently large period of time is proved in \cite{tao2012eventual} for smooth and convex $3$D domains. And more recently, a parabolic-elliptic simplification of \eqref{problema_P}, for $s = 1$, is studied in \cite{tao2019global}, yielding results on the existence and long-time behavior of global classical solutions in $n$-dimensional smooth domains.
    
    Still considering $s = 1$, there are some studies on the coupling of \eqref{problema_P} with models for fluids. In \cite{lorz2010coupled}, the author proves local existence of weak solutions for the chemotaxis-Navier-Stokes equations in $3$D smooth domains, while in \cite{duan2010global} the existence of global classical solutions is attained near constant states. In \cite{winkler2012global}, considering smooth and convex domains, existence and uniqueness of a global classical solution for the chemotaxis-Navier-Stokes equations is proved in $2$D and existence of global weak solutions which become smooth after a large enough period of time is proved for the chemotaxis-Stokes equations in $3$D. In \cite{jiang2015global} the results of \cite{winkler2012global} on the existence of solution are extended to non-convex domains. In \cite{winkler2014stabilization} the author studies the assymptotic behavior of the chemotaxis-Navier-Stokes equations in 2D domains with the chemotaxis and consumption terms generalized by using adequate functions depending on the chemical substance, proving the convergence towards constant states in the $L^{\infty}$-norm. Finally, in \cite{winkler2016global} existence of global weak solutions for the chemotaxis-Navier-Stokes equations is established in $3$D smooth and convex domains and in \cite{winkler2017far} the assymptotic behavior of these solutions is studied. An interesting and open question is whether the technique used in the present paper to study \eqref{problema_P} can be extended to the chemotaxis-fluids models. It is probable that the present approach to \eqref{problema_P} could be combined with the usual regularization for the fluid equations used in \cite{winkler2016global}, for example.
    
    \
    
    Another topic of interest in chemotaxis models is the existence or non-existence of blowing-up solutions. When it comes to the model \eqref{problema_P}, with $s = 1$, this question has been answered for $2$D convex domains, because existence and uniqueness of classical and uniformly bounded solutions is proved in \cite{tao2012eventual}. As far as we know, for $3$D domains this question remains open. In fact, there are partial results obtained under the assumption of adequate constraints relating the chemotaxis coefficient with $\norma{v^0}{L^{\infty}(\Omega)}$. On this subject, we refer the interested reader to \cite{tao2011boundedness} and \cite{baghaei2017boundedness}, for the problem \eqref{problema_P} with $s = 1$. In addition, we also have \cite{fuest2019analysis} and \cite{frassu2021boundedness}, where the authors extend these results to other related chemotaxis models with consumption.
    We remark that all the aforementioned works are carried out considering smooth domains. We also observe that, in \cite{jiang2015global}, the authors extend the results of \cite{winkler2012global} to non-convex domains, but some estimates that were time-independent become time-dependent. 
    Now, a first objective of this paper will be to prove time-independent regularity for general non-convex domains, for all models varying the potential consumption term $u^sv$ for any $s\ge 1$.
    
    \
    
    In the context of chemotaxis models, researchers are also interested in the control theory of these models. For the present chemoattraction with consumption model, with $s = 1$, coupled with fluids, an optimal control problem is studied in \cite{lopez2021optimal}.
   
   \
    
   Regarding the numerical approximation, although the numerical simulation of chemotaxis models is a relevant and growing research topic, when we turn to problem \eqref{problema_P}, we still find a relatively small amount of studies on its numerical approximation. This is probably due to the difficulty that is imposed by the complex relation between the chemoattraction and consumption terms. As we can see in Section \ref{section: original problem 3D}, the procedures on which the theoretical analysis relies are difficult to adapt in a numerical method.
   To the best of our knowledge,  for the numerical approximation of \eqref{problema_P},    
   we can cite two studies: \cite{duarte2021numerical} and \cite{FranciscoGiordano}. In \cite{duarte2021numerical} a chemotaxis-Navier-Stokes system is approached via Finite Elements (FE) analyzing optimal error estimates, assuming the existence of a sufficiently regular solution. In \cite{FranciscoGiordano}, inspired by the treatment given to the chemo-repulsion model with linear production in \cite{guillen2019unconditionally}, several finite element schemes are designed to approximate the system, focusing on properties such as conservation of cells, energy stability and approximate positivity rather than convergence.
  
  \
  
    In view of the relatively low number of studies on the numerical approximation of \eqref{problema_P} and the probable reasons for it, a second objective of the present study is to design an adequate background to the development of numerical approximations of the chemotaxis models \eqref{problema_P}. As far as we know, there are not studies on the model \eqref{problema_P} in the case $s > 1$, neither from the theoretical nor from the numerical point of view. 
    
    \
    
    In consequence, this work gives the following main contributions:
    
    (i) generalization of the model with the consumption term $u^sv$ for $s\geq 1$ (instead of only $uv$)
    
    (ii) enlargement of the class of considered domains, maintaining the no blow-up effect in the $2D$ case and some time-independent estimates  in $3D$ domains;
    
    (iii) the introduction and analysis of a regularized model, see \eqref{problema_P_m_intro} below, proving existence, uniqueness, regularity, positivity, a priori estimates and convergence towards the original model \eqref{problema_P}. That will be useful for designing numerical schemes that will be used in order to approximate the original model. This subject will be studied by the authors in a forthcoming paper.
    
    \
    
    We would like to make an observation regarding the rigor of the calculations. We have observed that, in some papers on analysis of chemotaxis models, singular functions are taken as test function (as for instance $log(u)$, without taking care that one only has $u\ge 0$). In our opinion, it should  be considered only as formal computations. Then, in this paper, we have done a great effort in order to guarantee that all of our computations be rigorous. Similarly to \cite{jungel2022analysis} for a cross-diffusion model, in order to make rigorous computations, we rely on a regularization procedure (for instance, taking $log(u+\epsilon)$ as test function). 

    \section{Main results}
    In order to present the main results, we introduce the following regularized problems, which depend on a truncation parameter $m \in \mathbb{N}$,
    \begin{equation}\tag{P$_m$}
    \left\{
      \begin{array}{l}
        \partial_t u_m - \Delta u_m  = - \nabla \cdot (a_m(u_m) \nabla v_m), \quad
        \partial_t v_m - \Delta v_m  = - a_m(u_m)^s v_m, \\
        \partial_\eta u_m \Big |_{\Gamma}  = \partial_\eta v_m \Big |_{\Gamma} = 0, \quad
        u_m(0)  = u^0_m, \quad v_m(0) = v^0_m,
      \end{array}
      \right.
      \label{problema_P_m_intro}
    \end{equation}
    where $u^0_m \geq 0$ and $v^0_m \geq 0$ are suitable regular approximations of $u^0$ and $v^0$, respectively, and $a_m(\cdot)$ is the following truncation of the identity function, from above and from below:
    \begin{equation}
      a_m(u) = \left \{
      \begin{array}{cl}
        - 1, & \mbox{ if } u \leq -1,  \\
        C^2 \mbox{ extension}, & \mbox{ if } u \in (-1,0), \\
        u, & \mbox{ if } u \in [0, m], \\
        C^2 \mbox{ extension}, & \mbox{ if } u \in (m,m + 2), \\
        m + 1, & \mbox{ if } u \geq m + 2.
      \end{array}
      \right .
      \label{truncamento_limitado_da_identidade}
    \end{equation}

    \
        
    With the objective of enlarging the class of considered domains, we state and demonstrate our results in terms of the regularity of the Poisson-Neumann
    \begin{equation}  \label{Neumann_problem}
      \left \{ \begin{array}{l}
        - \Delta w + w  = f  \ \mbox{ in } \Omega, \quad
        \partial_\eta w \Big |_{\Gamma}  = 0
      \end{array} \right.
    \end{equation}
    (see definition \ref{defi_regularidade_H_m} in page \pageref{defi_regularidade_H_m}); and, when necessary, in terms of the following technical hypothesis:
    \begin{hyp}
      For each $z \in H^2(\Omega)$ such that $\partial_{\eta} z \Big |_{\Gamma} = 0$ there is a sequence $\{ \rho_n \} \subset C^2(\overline{\Omega})$ such that $\partial_{\eta} \rho_n \Big |_{\Gamma} = 0$ and $\rho_n \to z$ in $H^2(\Omega)$.
      \label{hypothesis_density}
    \end{hyp}
    
    In order to to show that the Hypothesis \ref{hypothesis_density} is not too restrictive, we prove in Lemma \ref{lema_densidade}, in the appendix, that if the Poisson-Neumann problem has the $W^{3,p}$-regularity (see definition \ref{defi_regularidade_H_m} in page \pageref{defi_regularidade_H_m}), for $p > N$, then Hypothesis \ref{hypothesis_density} is satisfied.
   
   \
    
    Let us consider the average of $u^0$
    \begin{equation*}
      u^\ast = \frac{1}{\norm{\Omega}{}} \int_{\Omega}{u^0(x) \ dx}.
    \end{equation*}
    Now we highlight our main results in this work:   
    \begin{teo}{\bfseries ($\boldsymbol{3}$D. Existence of global weak solutions)}
      Let $\Omega \subset \mathbb{R}^3$ be a bounded domain such that the Neumann problem \eqref{Neumann_problem} has the $H^2$-regularity (see definition \ref{defi_regularidade_H_m} in page \pageref{defi_regularidade_H_m}) and Hypothesis \ref{hypothesis_density} is satisfied. Let $u^0 \in L^{1 + \varepsilon}(\Omega)$, for some $\varepsilon > 0$, if $s = 1$, and $u^0 \in L^s(\Omega)$, if $s > 1$, and $v^0 \in H^1(\Omega) \cap L^{\infty}(\Omega)$ be non-negative functions. Then there is a non-negative weak solution $(u,v)$ of the original problem \eqref{problema_P}, for $s \geq 1$, obtained through a limit of non-negative solutions $(u_m,v_m)$ of the regularized problems \eqref{problema_P_m_intro} as $m \to \infty$ and such that
      \begin{equation}
        \left \{
        \begin{array}{rl}
          \D{\int_{\Omega}}{u(t,x) \ dx} = \D{\int_{\Omega}}{u^0(x) \ dx}, \ a.e. \ t \in (0,\infty) \\
          0 \leq v(t,x) \leq \norma{v^0}{L^{\infty}(\Omega)}, \ a.e. \ (t,x) \in (0,\infty) \times \Omega,
        \end{array}
        \right.
        \label{pointwise_properties_3D}
      \end{equation}
      \begin{equation*}
        \left \{
        \begin{array}{rl}
          u \in L^{\infty}(0,\infty;L^s(\Omega)) \cap L^{5s/3}_{loc}([0,\infty);L^{5s/3}(\Omega)), \mbox{ if } s \geq 1, \\
          u^{s/2} \nabla v \in L^2(0,\infty;L^2(\Omega)), \mbox{ if } s \geq 1,
        \end{array}
        \right.
      \end{equation*}
      \begin{equation*}
        \left \{
        \begin{array}{rl}
          \nabla u \in L^2(0,\infty;L^s(\Omega)) \cap L^{\frac{5s}{3+s}}_{loc}([0,\infty);L^{\frac{5s}{3+s}}(\Omega)), & \mbox{ if } s \in [1,2), \\
          \nabla u \in L^2(0,\infty;L^2(\Omega)), & \mbox{ if } s \geq 2,
        \end{array}
        \right.
      \end{equation*}
      \begin{equation*}
        \left \{
        \begin{array}{rl}
          u \nabla v \in L^2(0,\infty;L^s(\Omega)), & \mbox{ if } s \in [1,2), \\
          u \nabla v \in L^2(0,\infty;L^2(\Omega)), & \mbox{ if } s \geq 2
        \end{array}
        \right.
      \end{equation*}
      and
      \begin{equation*}
        v \in L^{\infty}(0,\infty;H^1(\Omega)) \cap L^2(0,\infty;H^2(\Omega)),
        \quad
        \nabla v \in L^4(0,\infty; L^4(\Omega)).
      \end{equation*}
      \label{teo_3D_existencia}
    \end{teo}
    \begin{proof}[\bf Proof]
      This is proved in Section \ref{section: original problem 3D}.
    \end{proof}
    
    \begin{obsv}
      We remark that, from the regularities of $u$ and $v$ that are listed in Theorem \ref{teo_3D_existencia}, we can conclude that
      \begin{equation*}
        \left \{
        \begin{array}{rl}
          u_t \in L^2 \Big (0,\infty; \big ( W^{1,s/(s-1)}(\Omega) \big )' \Big ), & \mbox{ if } s \in [1,2), \\
          u_t \in L^2 \Big ( 0,\infty;\big ( H^1(\Omega) \big )' \Big ), & \mbox{ if } s \geq 2
        \end{array}
        \right.
      \end{equation*}
      and
      \begin{equation*}
        v_t \in L^2(0,\infty;L^{3/2}(\Omega)).
      \end{equation*}
      Attending to the regularity of $(u,v)$ given so far, one has that $(u,v)$ satisfies the $u$-equation of \eqref{problema_P} in a variational sense, while the $v$-equation is satisfied $a.e.$ in $(0,\infty) \times \Omega$. Moreover, the initial conditions have a sense because, thanks to the regularity of $u$, $v$, $u_t$ and $v_t$, one has that $(u,v)$ is weakly continuous from $[0,\infty)$ to $L^s(\Omega) \times H^1(\Omega)$, if $s \in [1,2]$, and $L^2(\Omega) \times H^1(\Omega)$, if $s \geq 2$ (see Chapter 3 of \cite{Temam}).
      \hfill \qedsymbol
    \end{obsv}

    \begin{obsv}
      Note that, for $s \in [1,2]$, the regularity of the fluxes of the $u$-equation of \eqref{problema_P}, namely, self diffusion $\nabla u$ and chemotaxis $u \nabla v$, increase as $s$ increases. When we consider $s > 2$, the regularity of $\nabla u$ and $u \nabla v$ do not increase as $s$ increases anymore. On the other hand, the regularity of the function $v$ is independent of $s$.
      \hfill \qedsymbol
    \end{obsv}
    
    \begin{teo}{\bfseries ($2$D. Existence and uniqueness of global strong solution)}
      Let $\Omega \subset \mathbb{R}^2$ be a bounded domain such that the Neumann problem \eqref{Neumann_problem} has the $W^{2,3}$-regularity (see definition \ref{defi_regularidade_H_m} in page \pageref{defi_regularidade_H_m}) and Hypothesis \ref{hypothesis_density} is satisfied. Let $u^0 \in H^2(\Omega)$ and $v^0 \in H^2(\Omega)$ be such that $u^0 \geq 0$ and $v^0 \geq 0$. Then there is a unique non-negative solution $(u,v)$ for the original problem \eqref{problema_P}, for $s \geq 1$, satisfying  \eqref{pointwise_properties_3D} and
      \begin{equation*}
        \begin{array}{c}
          u, v \in L^{\infty}(0,\infty;H^2(\Omega)), \quad 
          (u - u^\ast), v \in L^2(0,\infty;W^{2,3}(\Omega))
          \\
          \Delta u, \Delta v, u_t, v_t \in L^2(0,\infty;H^1(\Omega)).
        \end{array}
      \end{equation*}
      In particular, $u$ does not blow-up neither at finite nor infinite time, that is $u \in L^{\infty}(0,\infty;L^{\infty}(\Omega))$. Consequently, there is $m_0 \in \mathbb{N}$ such that, for all $m \in [m_0, \infty)$, the solution of \eqref{problema_P_m_intro} is also the solution of \eqref{problema_P}, that is,
      \begin{equation*}
        (u_m,v_m) = (u,v) \ a.e. \mbox{ in } (0,\infty) \times \Omega.
      \end{equation*}
      \label{teo_2D_existencia_unicidade_estimativas}
    \end{teo}
    \begin{proof}[\bf Proof]
     This is proved in Section \ref{section: original problem 2D}.
    \end{proof}
    In this case, both equations of \eqref{problema_P} are satisfied $a.e.$ in $(t,x) \in (0,\infty) \times \Omega$.
    
    \
    
    In contrast with Theorem \ref{teo_2D_existencia_unicidade_estimativas}, where more regular initial data lead to a more regular solution in general domains, in \cite{wang2021immediate} problem \eqref{problema_P}, with $s = 1$, is studied in $2$D disks but with highly irregular initial data. More precisely, the authors consider $u^0$ to be a measure and $v^0 \in L^2(\Omega)$, both nonnegative and radially symmetric.
    
    \
    
    We observe that in 3D domains it is not possible to state a complete result such as Theorem \ref{teo_2D_existencia_unicidade_estimativas}. This is due to the gap between the regularity of the solutions provided by the existence result and the regularity needed to prove uniqueness. Notice that $v$ does not blow-up neither at finite nor infinite time. On the other hand, to the best of our knowledge, whether $u$ may blow-up or not is an open problem.
    
    \
    
    The rest of the paper is organized as follows. In Section \ref{section: preliminary} we give  the basic theoretical and classical results which will be used in the rest of the paper. Section \ref{section: regularized problem} is devoted to discuss the regularity's properties of the solutions of regularized problem \eqref{problema_P_m_intro}. In Section \ref{section: original problem 3D} we deal with the proof of Theorem \ref{teo_3D_existencia} and, finally, Section \ref{section: original problem 2D} is dedicated to the proof of Theorem \ref{teo_2D_existencia_unicidade_estimativas}.
    

\section{Functional Framework}
  \label{section: preliminary}
 

  In this section we present some technical tools which will be used in the rest of the paper. 
  For $p \in [1,\infty]$, we denote by $L^p(\Omega)$,  the usual Banach spaces of $p$-integrable Lebesgue-mensurable functions, with the norm $\norma{\cdot}{L^p(\Omega)}$. We recall that $L^2(\Omega)$ is a Hilbert space with the inner product
  \begin{equation*}
    \prodl{f}{g} = \int_{\Omega}{f(x) g(x) \ dx}.
  \end{equation*}
  We also denote by $W^{l,p}(\Omega)$, with $l\in \mathbb{N}$,  the usual Sobolev space, equipped with the usual norm $\norma{\cdot}{W^{l,p}(\Omega)}$; for $p=2$, we denote $W^{l,2}(\Omega)$ by $H^{l}(\Omega)$, with norm $\norma{\cdot}{H^{l}(\Omega)}$. 
  
  If $X$ is a Banach space, then $L^p(a,b;X)$ is the Bochner space with the norm
  \begin{equation*}
    \norma{v}{L^p(a,b;X)} = \left ( \int_a^b{\norma{v(t)}{X}^p \ dt} \right )^{1/p}, \quad \norma{v}{L^{\infty}(a,b;X)} = \esssup_{t \in (a,b)} \norma{v(t)}{X}.
  \end{equation*}
  
  If $p = 2$ and $X$ is a Hilbert space then $L^2(a,b;X)$ is a Hilbert space with the inner product
  \begin{equation*}
    \prodl{u}{v}_{L^2(a,b;X)} = \int_a^b{\prodl{u(t)}{v(t)}_X \ dt}, \quad \forall u,v \in L^2(a,b;X),     
  \end{equation*}
  where $\prodl{\cdot}{\cdot}_X$ denotes the inner product of $X$.
  
  Next we present some interpolation inequalities and other results which will be of frequent use in the article. Unless otherwise stated, we consider $\Omega \subset \mathbb{R}^N$ ($N=2,3$) to be an open, bounded and locally Lipschitz domain.
    
  \begin{lema} \label{interpolation} We have the following interpolation inequalities:
    
    \noindent 1.      Let $1 \leq p < q \leq \infty$,    $\theta \in (0,1)$ and $r \in [p,q]$, with $\frac{1}{r} = \frac{\theta}{p} + \frac{1 - \theta}{q}$. If $f \in L^p(\Omega) \cap L^q(\Omega)$ then $f \in L^r(\Omega)$ and
    \begin{equation*}
      \norma{f}{L^r(\Omega)} \leq \norma{f}{L^p(\Omega)}^{\theta} \norma{f}{L^q(\Omega)}^{1 - \theta}.
    \end{equation*}
    
    \noindent 2. There exist (different) constants $\beta > 0$ such that
     
    \textnormal{(i)} if $N = 2$ then
    \begin{equation}
      \norma{v}{L^4(\Omega)} \leq \beta \norma{v}{L^2(\Omega)}^{1/2} \norma{v}{H^1(\Omega)}^{1/2}, \forall v \in H^1(\Omega);
      \label{desig_ladyzhenskaya}
    \end{equation}
          
        \textnormal{(ii)} if $N = 3$ then
          \begin{equation*}
            \norma{v}{L^4(\Omega)} \leq \beta \norma{v}{L^2(\Omega)}^{1/4} \norma{v}{H^1(\Omega)}^{3/4}, \forall v \in H^1(\Omega).
          \end{equation*}
      
      \label{lema_desig_lady2}
       
       \textnormal{(iii)} if $N = 3$ then
          \begin{equation*}
            \norma{v}{L^3(\Omega)} \leq \beta \norma{v}{L^2(\Omega)}^{1/2} \norma{v}{H^1(\Omega)}^{1/2}, \forall v \in H^1(\Omega).
          \end{equation*}
  \end{lema}
  \begin{proof}[\bf Proof]
    See \cite{Brezis_functional} and \cite{Temam}.
  \end{proof}
  
  \begin{lema}[\bfseries Poincare's Inequality] \label{lema_desigualdade_poincare_media_nula}
    Let $v \in W^{1,p}(\Omega)$ and $v^{\ast} = \dfrac{1}{\norm{\Omega}{}} \D{\int_{\Omega} v(x) \ dx}$. Then there is a constant $C > 0$ which is independent of $v$ and such that
    \begin{equation*}
      \norma{v - v^{\ast}}{W^{1,p}(\Omega)} \leq C \norma{\nabla v}{L^p(\Omega)}.
    \end{equation*}
  \end{lema}
  \begin{proof}[\bfseries Proof]
    It can be proved by contradiction as in \cite{Evans2010}.
  \end{proof}
  
  \begin{lema} \label{lema_a_m_elevado_a_s}
    Let $w_1$ and $w_2$ be nonnegative real numbers. For each $s \geq 1$ we have
    \begin{equation*}
      \norm{w_2^s - w_1^s}{} \leq s \norm{w_2 + w_1}{}^{s-1} \norm{w_2 - w_1}{}.
    \end{equation*}
  \end{lema}
  \begin{proof}[\bfseries Proof]
    Without loss of generality, suppose $w_1 \leq w_2$. Then
    \begin{align*}
      \norm{w_2^s - w_1^s}{} & = s \norm{\int_{w_1}^{w_2}{r^{s-1} \ dr}}{}  \leq s \norm{\int_{w_1}^{w_2}{\sup_{r \in (w_1,w_2)}{r}^{s-1} \ dr}}{} \\
      & \leq s \norm{\int_{w_1}^{w_2}{\norm{w_2 + w_1}{}^{s-1} \ dr}}{} = s \norm{w_2 + w_1}{}^{s-1} \norm{w_2 - w_1}{}.
    \end{align*}
  \end{proof}
    
  Using Lemma \ref{lema_a_m_elevado_a_s}, we can prove the following.
  
  \begin{lema} \label{lema_convergencia_w_elevado_a_s}
    Let $p \in (1, \infty)$ and let $\{ w_n \}$ be a sequence of nonnegative functions in $L^p(0,T;L^p(\Omega))$ such that $w_n \to w$ in $L^p(0,T;L^p(\Omega))$ as $n \to \infty$. Then, for every $r \in (1,p)$, $w_n^r \to w^r$ in $L^{p/r}(0,T;L^{p/r}(\Omega))$ as $n \to \infty$.
  \end{lema}
  
  \begin{lema} \label{lema_integral_de_a_m}
    Let $p \in (1,\infty)$. Let $\{ w_m \}$ be a sequence of nonnegative functions which is uniformly bounded in $L^{\infty}(0,\infty;L^p(\Omega))$ with respect to $m$ and defined for every $t \in (0,\infty)$. If there is $\alpha > 0$ such that
    \begin{equation*}
      \int_{\Omega}{w_m(t,x) \ dx} \geq \alpha, \ \forall t \in (0,\infty), \ \forall m \in \mathbb{N},
    \end{equation*}
    then there exist $\beta > 0$ and $m_0 \in \mathbb{N}$ large enough such that
    \begin{equation*}
      \int_{\Omega}{a_m(w_m(t,x)) \ dx} \geq \beta, \ \forall \ t \in (0,\infty), \ \forall m \geq m_0.
    \end{equation*}
  \end{lema}
  \begin{proof}[\bf Proof]
    For every $t \in (0,\infty)$, let
    \begin{equation*}
      S_m(t) = \left \{ x \in \Omega \ \Big | \ w_m(t,x) > m \right \} .
    \end{equation*}
    Then, for all $t \in (0,\infty)$, we have
    \begin{equation*}
      \int_{S_m(t)}{m^p \ dx} + \int_{\Omega \backslash S_m(t)}{a_m(w_m(t,x))^p \ dx} \leq \int_{\Omega}{a_m(w_m(t,x))^p \ dx} \leq C_1(p).
    \end{equation*}
    This implies
    \begin{equation}
      \norm{S_m(t)}{} = \int_{S_m(t)}{dx} \leq C_1(p) \left ( \frac{1}{m} \right )^p.
      \label{eq_estimativa_medida_de_S_m}
    \end{equation}
    We have
    \begin{align*}
      \int_{\Omega}{a_m(w_m(t,x)) \ dx} & = \int_{\Omega \backslash S_m(t)}{w_m(t,x) \ dx} + \int_{S_m(t)}{a_m(w_m(t,x)) \ dx} \\
      & = \int_{\Omega}{w_m(t,x) \ dx} - \int_{S_m(t)}{w_m(t,x) - a_m(w_m(t,x)) \ dx} \\
      & \geq \alpha - \int_{S_m(t)}{w_m(t,x) - a_m(w_m(t,x)) \ dx}.
    \end{align*}
    
    To finish the proof, we show that
    \begin{equation*}
      \lim_{m \to \infty}{\int_{S_m(t)}{w_m(t,x) - a_m(w_m(t,x)) \ dx}} = 0,
    \end{equation*}
    uniformly with respect to $t \in (0,\infty)$. In fact, using Hölder's inequality and \eqref{eq_estimativa_medida_de_S_m}, we obtain
    \begin{align*}
      \int_{S_m(t)}{w_m(t,x) - a_m(w_m(t,x)) \ dx} & \leq \left ( \int_{S_m(t)}{\norm{w_m(t,x) - a_m(w_m(t,x))}{}^p \ dx} \right )^{\frac{1}{p}} \left ( \int_{S_m(t)}{dx} \right )^{\frac{p-1}{p}} \\
      & \leq C_2(p) \left ( \int_{S_m(t)}{dx} \right )^{\frac{p-1}{p}}  \leq C(p) \left ( \frac{1}{m} \right )^{p-1}.
    \end{align*}
    for all $t \in (0,\infty)$. Therefore we have
    \begin{equation*}
      \int_{\Omega}{a_m(w_m(t,x)) \ dx} \geq \alpha - C(p) \left ( \frac{1}{m} \right )^{p-1}
    \end{equation*}
    and then we can choose $m_0$ large enough such that
    \begin{equation*}
      \int_{\Omega}{a_m(w_m(t,x)) \ dx} \geq \beta =  \frac{\alpha}{2},
    \end{equation*}
    for example, completing the proof.
  \end{proof}
  
  The following result (Corollary $4$ of \cite{Simon1986compact}) establishes a criterion of compactness in Bochner spaces.
  
  \begin{lema}[\bfseries Compactness in Bochner spaces]
    Let $X,B$ and $Y$ be Banach spaces, let
    \begin{equation*}
      F = \Big \{ f \in L^1(0,T;Y) \ \Big | \ \partial_t f \in L^1(0,T;Y) \Big \}.
    \end{equation*}
    Suppose that $X \subset B \subset Y$, with compact embedding $X \subset B$ and continuous embedding $B \subset Y$. Let the set $F$ be bounded in $L^q(0,T;B) \cap L^1(0,T;X)$, for $1 < q \leq \infty$, and $\Big \{ \partial_t f, \ \forall f \in F \Big \}$ be bounded in $L^1(0,T;Y)$. Then $F$ is relatively compact in $L^p(0,T;B)$, for $1 \leq p < q$.
    \label{lema_Simon}
  \end{lema}
  
  \begin{lema}[\bfseries Gronwall's inequality]
    Let $f$, $g$ and $h$ be nonnegative functions such that $f \in W^{1,1}(0,T)$ and $g, h \in L^1(0,T)$, for some $T > 0$. Let
    \begin{equation*}
      G(t) = \int_0^t{g(r) \ dr} \qquad \mbox{ and } \qquad H(t) = \int_0^t{h(r) \ dr}.
    \end{equation*}
    If $f$ is such that
    \begin{equation*}
      \dfrac{d}{dt} f(t) \leq g(t) + h(t) f(t), \ a.e. \ t \in (0,T),
    \end{equation*}
    then
    \begin{equation*}
      f(t) \leq \Big ( G(t) + f(0) \Big ) e^{H(t)}, \ a.e. \ t \in (0,T).
    \end{equation*}
    \label{lema_gronwall}
  \end{lema}
    
    \begin{lema}
      Let $\mathcal{X}$ and $\mathcal{Y}$ be Banach spaces. Let $S: \mathcal{X} \longrightarrow \mathcal{Y}$ be a continuous linear transformation. If $f \in L^1((0,\infty);\mathcal{X})$ then $Sf \in L^1((0,\infty);\mathcal{Y})$ and
      \begin{equation*}
        \int_0^{\infty}{Sf \ dt} = S \int_0^{\infty}{f \ dt}.
      \end{equation*}
      \label{lema_integral_da_transformacao}
    \end{lema}
    \begin{proof}[\bf Proof]
      See the section about the Bochner's integral in the book of Yosida \cite{Yosida}.
    \end{proof}
    We will apply this lemma for $S: L^2(\Omega) \longrightarrow \mathbb{R}$ given by $Sf = \D{\int_{\Omega}{f \ dx}}$.
    
    \begin{lema}
      If $w, w_t \in L^1_{loc}((0,\infty);L^2(\Omega))$ then $\D{\dfrac{d}{dt}(\int_{\Omega}{w(\cdot,x) \ dx})}\in L^1_{loc}((0,\infty);L^2(\Omega))$ and
      \begin{equation*}
        \dfrac{d}{dt}(\int_{\Omega}{w(\cdot,x) \ dx}) = \int_{\Omega}{w_t(\cdot,x) \ dx }.
      \end{equation*}
      \label{lema_derivada_da_transformacao}
    \end{lema}
    \begin{proof}[\bf Proof]
      We look at the integral
      \begin{equation*}
        \int_0^{\infty}{\int_{\Omega}{w_t(t,x) \ dx } \ \psi(t) \ dt},
      \end{equation*}
      for every $\psi \in C^{\infty}_c((0,\infty))$. Since the integral over $\Omega$ is a linear transformation and $\psi(t)$ is a real number, for each $t \in (0,\infty)$, we have
      \begin{equation*}
        \int_{\Omega}{w_t(t,x) \ dx} \ \psi(t) = \int_{\Omega}{w_t(t,x) \psi(t) \ dx},
      \end{equation*}
      for $a.e.$ fixed $t$. Then we use lemma \ref{lema_integral_da_transformacao} with $f = w_t \psi \in L^1((0,\infty);L^2(\Omega))$ to write
      \begin{align*}
        & \int_0^{\infty}{\int_{\Omega}{w_t(t,x) \ dx} \ \psi(t) \ dt} = \int_0^{\infty}{\int_{\Omega}{w_t(t,x) \psi(t) \ dx} \ dt} = \int_{\Omega}{ \int_0^{\infty}{w_t(t,x) \psi(t) \ dt} \ dx} \\
        & = - \int_{\Omega}{ \int_0^{\infty}{w(t,x) \psi_t(t) \ dt} \ dx} = - \int_0^{\infty}{\int_{\Omega}{w(t,x) \ dx} \ \psi_t(t) \ dt}, \forall \psi \in C^{\infty}_c((0,{\infty})).
      \end{align*}
    \end{proof}
    
    Let $C^{\infty}_c(0,\infty;L^2(\Omega))$ denote the space of the infinitely differentiable functions defined in $(0,\infty)$ with range in $L^2(\Omega)$ and with compact support in $(0,\infty)$. Consider the space
    \begin{equation*}
      H^1((0,\infty);L^2(\Omega)) = \left \{ w \in L^2(0,\infty;L^2(\Omega)) \ \Big | \ w_t \in L^2(0,\infty;L^2(\Omega)) \right \},
    \end{equation*}
    which is a Hilbert space with the norm
    \begin{equation*}
      \norma{w}{H^1((0,\infty);L^2(\Omega))} = \Big ( \norma{w}{L^2(0,\infty;L^2(\Omega))}^2 + \norma{w_t}{L^2(0,\infty;L^2(\Omega))}^2 \Big )^{1/2}.
    \end{equation*}
    
    \begin{lema}
      $C^{\infty}_c(0,\infty;L^2(\Omega))$ is dense in $H^1((0,\infty);L^2(\Omega))$.
      \label{lema_densidade_Lions_Magenes}
    \end{lema}
    \begin{proof}[\bf Proof]
      See Lions and Magenes \cite{Lions_Magenes}.
    \end{proof}


Due to the nature of our problem, it is convenient to have information on the Neumann problem. At first, we mention that if a function $w \in H^1(\Omega)$ and $\Delta w \in L^2(\Omega)$ then we can define the trace of the normal derivative of $w$, which will be denoted by $\partial_{\eta} w \Big |_{\Gamma}$.

Concerning the regularity of solutions of the Poisson-Neumann problem \eqref{Neumann_problem}, we present the following
    \begin{defi}
      Let $w^{\ast} = \dfrac{1}{\norm{\Omega}{}} \D{\int_{\Omega}}{w(x) \ dx}$ and suppose $w \in H^1(\Omega)$ is a weak solution of \eqref{Neumann_problem} with $f \in L^p(\Omega)$. If this implies that $w \in W^{2,p}(\Omega)$ with
      \begin{equation*}
        \norma{w - w^{\ast}}{W^{2,p}(\Omega)} \leq C \norma{\Delta w}{L^p(\Omega)},
      \end{equation*}
      then we say that the Poisson-Neumann problem \eqref{Neumann_problem} has the $W^{2,p}$-regularity. If $p = 2$ we say that the problem has $H^2$-regularity.
      \label{defi_regularidade_H_m}
    \end{defi}
    
    \begin{obsv}
      We remark that if the Neumann problem has the $W^{2,p}$-regularity then, reminding that $\nabla w = \nabla (w - w^{\ast})$ and $D^2 w = D^2 (w - w^{\ast})$, we have
      \begin{equation}
        \norma{\nabla w}{W^{1,p}(\Omega)} \leq C \norma{\Delta w}{L^p(\Omega)}.
        \label{obsv_norma_de_nabla_v}
      \end{equation}
      \hfill \qedsymbol
    \end{obsv}
    
    \begin{obsv}   
      According to Grisvard \cite{Grisvard}, if $f \in L^p(\Omega)$, $p \in [1,\infty]$, and the boundary $\Gamma$ is at least $C^{1,1}$, then the solution $w$ of the Neumann problem \eqref{Neumann_problem} belongs to $W^{2,p}(\Omega)$ with continuous dependence on the data (or it has the $W^{2,p}$-regularity for all $p \in [1,\infty]$). The aforementioned result is also true if $\Omega$ is a polygon, that is, a polyhedron in $\mathbb{R}^2$, or if $\Omega$ is convex and $p=2$. For more regular domains, if $f \in W^{l,p}(\Omega)$ and the boundary $\Gamma$ is at least $C^{l+1,1}$, then the solution $w$ of the Neumann problem \eqref{Neumann_problem} belongs to $W^{l+2,p}(\Omega)$. \hfill \qedsymbol
    \end{obsv}
    
    
  We end this section recalling the concept of positive and negative parts of a function. For $w \in L^p(\Omega)$, $1 \leq p \leq \infty$,  the positive and negative parts of $w$ are given by
    \begin{equation*}
      w_+(x) = \max{ \{ 0, w(x) \} } \mbox{ and } w_-(x) = \min{ \{ 0, w(x) \} },
    \end{equation*}
    respectively. Then $w = w_+ + w_-$ and $\norm{w}{} = w_+ - w_-$; besides, if $w \in H^1(\Omega)$ then $w_+,w_- \in H^1(\Omega)$ with
\[    \begin{array}{lcl}
      \nabla w_+(x) = \left \{
      \begin{array}{rl}
        \nabla w, & \mbox{ if } w(x) > 0,  \\
           0, & \mbox{ if } w(x) < 0,
      \end{array}
      \right.
 & \hspace{1cm}  \mbox{and} \hspace{1cm} 
&    
      \nabla w_-(x) = \left \{
      \begin{array}{rl}
        \nabla w, & \mbox{ if } w(x) < 0,  \\
           0, & \mbox{ if } w(x) > 0.
      \end{array}
      \right.
    \end{array}
  \]
  
   For more details on truncations applied to $H^1(\Omega)$ functions we suggest Gilbarg and Trudinger \cite{GilbargTrudinger}.
    

\section{The Regularized Problem}
  \label{section: regularized problem}
  
  In this section we define and analyze the regularized problem \eqref{problema_P_m_intro}, based on  the truncation of the identity $a_m(u)$ given in  \eqref{truncamento_limitado_da_identidade}. 
  We remark the following  properties
  \begin{equation*}
    a_m(u) \leq u, \quad \forall u \geq 0, 
  \end{equation*}
  \begin{equation}
    \norm{a_m(u)}{} \leq m, \quad 
    \norm{a_m'(u)}{}, \norm{a_m''(u)}{} \leq C, 
    \quad \forall u\in \mathbb{R},
    \label{limitacao_do_truncamento_a_m}
  \end{equation}
  where  $C > 0$ is a constant independent of $m \in \mathbb{N}$.
  
  For each $m \in \mathbb{N}$, let $(u_m,v_m)$ be the solution of \eqref{problema_P_m_intro} with initial data $u^0_m,v^0_m \in C^\infty(\overline\Omega)$ with $u^0_m$ and $v^0_m$ being mollifier regularizations of $u^0$ and $v^0$ extended to $\mathbb{R}^N$. In fact, $u^0\in L^p(\Omega)$ is extended by zero, while $v^0\in H^1(\Omega)$ is extended in the space $H^1(\mathbb{R}^N)$. In particular, these regularizations have the following properties:
  \begin{equation}
    u^0_m\geq 0,\quad \int_\Omega u^0_m = \int_\Omega u^0,\quad 
    u^0_m \rightarrow u^0 \mbox{ strongly in } L^p(\Omega), \mbox{ as } m \to \infty,
    \label{convergencia_dado_inicial_u_m}
  \end{equation}
  for $p = 1 + \varepsilon$, for some $\varepsilon > 0$, if $s = 1$, and $p = s$, if $s > 1$, and
  \begin{equation} \label{convergencia_dado_inicial_v_m}
    \essinf v^0 \le v^0_m \le \esssup v^0,\quad 
    v^0_m \rightarrow v^0 \mbox{ strongly in } H^1(\Omega), \mbox{ as } m \to \infty.
  \end{equation}
  
  
  \subsection{Existence and uniqueness of problem \texorpdfstring{(\ref{problema_P_m_intro})}{(Pm)}}
    \label{Subsection:regularized_problem_approx_solutions}
          
    We use the Galerkin's method based on the set of eigenfunctions $\{ \varphi_j \}$ of the operator $(- \Delta + I)$ with Neumann homogeneous boundary condition. Unless otherwise stated, we will proceed under the assumption that Poisson-Neumann problem \eqref{Neumann_problem} has the $H^2$-regularity. Then $\{ \varphi_j \}$ is a basis of $H^2(\Omega)$.
     
    Let $X^n$ be the finite $n-$dimensional space generated by the first $n$ elements of the set $\{ \varphi_j \}$. Then, we look for Galerkin solutions $(u_n,v_n)$ of the form
    \begin{equation*}
      u_n(t,x) = \sum_{j=1}^n{g^n_j(t) \varphi_j(x)}
      \hspace{1cm} \mbox{ and } \hspace{1cm}
      v_n(t,x) = \sum_{j=1}^n{h^n_j(t) \varphi_j(x)}
    \end{equation*}
     such that
      \begin{align}
        \prodl{\partial_t u_n}{\varphi_i} + \prodl{\nabla u_n}{\nabla \varphi_i} & = \prodl{a_m(u_n) \nabla v_n}{\nabla \varphi_i}, \label{sistema_aproximado_u_n} \\
        \prodl{\partial_t v_n}{\varphi_i} - \prodl{\Delta v_n}{\varphi_i} & = - \prodl{a_m(u_n)^s v_n}{\varphi_i}, \label{sistema_aproximado_v_n} \\
        u_n(0) & = P_n(u^0_m), \quad v_n(0) = P_n(v^0_m), \label{cond_iniciais_aproximado}
      \end{align}
      for $i = 1, \dots, n$, where $P_n(u^0_m)$ and $P_n(v^0_m)$ are orthogonal projections of $u^0_m$ and $v^0_m$ from $H^1(\Omega)$ into $X^n$.
      Since the application of the Galerkin's method is a very standard procedure, some details (such as the proof of existence of the Galerkin solutions, the obtaining of {\it a priori} estimates and the passage to the limit as $n \to \infty$) will be omitted here.
     
     \
      
      In order to obtain $n$-independent \emph{a priori} estimates for $(u_n, v_n)$, we test \eqref{sistema_aproximado_u_n} by $u_n\in X^n$ and \eqref{sistema_aproximado_v_n} by $v_n\in X^n$ and $- \Delta v_n \in X^n$. Then we can also test \eqref{sistema_aproximado_v_n} by $\Delta^2 v_n\in X^n$ and \eqref{sistema_aproximado_u_n} by $- \Delta u_n\in X^n$. Taking the truncation $a_m(\cdot)$ and its bounds \eqref{limitacao_do_truncamento_a_m} into account, it is not difficult to obtain the following \emph{a priori} bounds (independent of $n$) for each final time $T > 0$:
      \begin{equation}
        \begin{array}{c}
          (u_n,v_n)_n \mbox{ is bounded in } L^{\infty}(0,T;H^1(\Omega) \times H^2(\Omega)), \\[6pt]
          (\Delta u_n, \Delta v_n)_n, (\partial_t u_n, \partial_t v_n)_n \mbox{ are bounded in } L^2(0,T;L^2(\Omega) \times H^1(\Omega)).
        \end{array}
        \label{limitacao_u_n_v_n}
      \end{equation}
      Therefore the Galerkin solution, $(u_n,v_n)$, is defined up to infinity time.
 
 \
     
      Besides, if we assume that the Poisson-Neumann problem \eqref{Neumann_problem} has the $W^{2,3}$-regularity, then we have
      \begin{equation}
        (v_n)_n \mbox{ is bounded in } L^2(0,T;W^{2,3}(\Omega)).
        \label{limitacao_v_n_forte}
      \end{equation}
      We can also test \eqref{sistema_aproximado_u_n} by $\Delta^2 u_n \in X^n$ and, using \eqref{limitacao_v_n_forte}, we obtain the $n$-independent bounds
      \begin{equation}
        \begin{array}{c}
          (u_n)_n \mbox{ is bounded in } L^{\infty}(0,T;H^2(\Omega)) \cap L^2(0,T;W^{2,3}(\Omega)), \\[6pt]
          (\partial_t u_n)_n, (\Delta u_n)_n \mbox{ are bounded in } L^2(0,T;H^1(\Omega)).
        \end{array}
        \label{limitacao_u_n_forte}
      \end{equation}
      
      Now, the {\it a priori} bounds \eqref{limitacao_u_n_v_n}, compactness results in the weak and weak* topologies (see \cite{Brezis}) and compactness results in Bochner spaces (given in Lemma \ref{lema_Simon}), for each $T > 0$, allow us to conclude that there exist limit functions $u_m$ and $v_m$ such that, up to a subsequence,
      \begin{equation*}
        u_n \rightarrow u_m \mbox{ weakly* in } L^{\infty}(0,T;H^1(\Omega)), \mbox{ weakly in } L^2(0,T;H^2(\Omega)) \mbox{ and strongly in } L^2(0,T;H^1(\Omega)),
      \end{equation*}
      \begin{equation*}
        v_n \rightarrow v_m \mbox{ weakly* in } L^{\infty}(0,T;H^2(\Omega)), \mbox{ weakly in } L^2(0,T;H^2(\Omega)) \mbox{ and strongly in } L^2(0,T;H^1(\Omega)).
      \end{equation*}
      
      Then using these convergences and passing to the limit in the approximate system it follows that $(u_m,v_m)$ satisfies \eqref{problema_P_m_intro} $a.e.$ in $(0,\infty) \times \Omega$. One can prove that the solution $(u_m,v_m)$ is unique by straightforward calculations.
      
      Thus, in this subsection, for each fixed $m \in \mathbb{N}$, we have proved the existence and uniqueness of $(u_m, v_m)$, solution of \eqref{problema_P_m_intro}, such that
      \begin{equation*}
        \begin{array}{c}
          u_m \in L^{\infty}_{loc}([0,\infty);H^1(\Omega)) \cap L^2_{loc}([0,\infty);H^2(\Omega)), \quad
          v_m \in L^{\infty}_{loc}([0,\infty);H^2(\Omega)), \\
          \partial_t u_m \in L^2_{loc}([0,\infty);L^2(\Omega)), \quad
          \partial_t v_m, \Delta v_m \in L^2_{loc}([0,\infty);H^1(\Omega)).
        \end{array}
      \end{equation*}
      If we assume that the Poisson-Neumann problem \eqref{Neumann_problem} has the $W^{2,3}$-regularity, then it stems from the stronger $n$-independent bounds \eqref{limitacao_v_n_forte} and \eqref{limitacao_u_n_forte} that
      \begin{equation*}
        \begin{array}{c}
          u_m, v_m \in L^{\infty}_{loc}([0,\infty);H^2(\Omega)) \cap L^2_{loc}([0,\infty);W^{2,3}(\Omega)), \\
          \partial_t u_m, \partial_t v_m, \Delta u_m, \Delta v_m \in L^2_{loc}([0,\infty);H^1(\Omega)).
        \end{array}
      \end{equation*}


    \subsection{Regularity up to infinity time of problem \texorpdfstring{(\ref{problema_P_m_intro})}{(Pm)}}
      \label{subsection_aux_problem_regularity}
      
      Continuing the analysis, in the present subsection we prove the following main result.
      \begin{teo}{\bf (Regularity up to infinity time of \eqref{problema_P_m_intro})}
        Let $u^0_m$ and $v^0_m$ be approximations of $u^0$ and $v^0$ as defined in the beginning of Section \ref{section: regularized problem}. Under the assumption that the Poisson-Neumann problem \eqref{Neumann_problem} has the $H^2$-regularity, there is a unique solution $(u_m,v_m)$ of \eqref{problema_P_m_intro} such that
        \begin{equation*}
          u_m(t,x) \geq 0, \ a.e. \ (t,x) \in (0,\infty) \times \Omega,
        \end{equation*}
        \begin{equation*}
          \D{\essinf_{x \in \Omega}{v^0(x)}} e^{-m^s t} \leq v_m(t,x) \leq \norma{v^0}{L^{\infty}(\Omega)}, \ a.e. \ (t,x) \in (0,\infty) \times \Omega,
        \end{equation*}
        with the following regularity:
        \begin{equation}
          \begin{array}{c}
            (u_m - u^{\ast}) \in L^{\infty}(0,\infty;H^1(\Omega)) \cap L^2(0,\infty;H^2(\Omega)), \quad
            v_m \in L^{\infty}(0,\infty;H^2(\Omega)) \cap L^2(0,\infty;H^2(\Omega)), \\
            \partial_t u_m \in L^2(0,\infty;L^2(\Omega)), \quad
            \partial_t v_m, \Delta v_m \in L^2(0,\infty;H^1(\Omega)),
          \end{array}
          \label{regularidade_solucao_u_m_v_m_reg_H_2}
        \end{equation}
        where $u^{\ast} = \frac{1}{\norm{\Omega}{}} \int_{\Omega}{u^0(x) \ dx}$. Additionally, if we assume that the Poisson-Neumann problem \eqref{Neumann_problem} has the $W^{2,3}$-regularity, then
        \begin{equation}
          \begin{array}{c}
            u_m \in L^{\infty}(0,\infty;H^2(\Omega)), \quad
            \partial_t u_m, \Delta u_m \in L^2(0,\infty;H^1(\Omega)).
          \end{array}
          \label{regularidade_solucao_u_m_v_m_reg_W_2_3}
        \end{equation}
        \label{teo_regularidade_solucao_do_problema_regularizado_ate_infinito}
      \end{teo}
      
      The proof of Theorem 15 will be carried out along the subsection in several steps. We begin with the proof of some pointwise estimates, in Lemma \ref{lemma_positividade_u_v_m}, and some direct $m$-independent estimates for the solution $(u_m,v_m)$ of the regularized problem, in Lemma \ref{lemma_limitacao_u_v_independente_m}. Next we prove the weak regularity up to infinity time, in Lemma \ref{lema_comportamento_assintotico_u_v_m}, and use it to finish the proof of Theorem \ref{teo_regularidade_solucao_do_problema_regularizado_ate_infinito}.
      
      \begin{lema}{\bf (Pointwise $\boldsymbol{m}$-uniform estimates for $\boldsymbol{(u_m,v_m)}$)}
        \begin{enumerate}
          \item If $u^0_m(x) \geq 0 \ a.e. \ x \in \Omega$ then $u_m(t,x) \geq 0 \ a.e. \ (t,x) \in (0,\infty) \times \Omega$;
          \item If $v^0(x) \geq 0 \ a.e. \ x \in \Omega$ and $v^0 \in L^{\infty}(\Omega)$ then
          \begin{equation*}
            \D{\essinf_{x \in \Omega}{\{ v^0(x) \} }} \, \exp(-m^s t) \leq v_m(t,x) \leq \norma{v^0}{L^{\infty}(\Omega)} \ a.e. (t,x) \in (0,\infty) \times \Omega;
          \end{equation*}
        \end{enumerate}
        \label{lemma_positividade_u_v_m}
      \end{lema}
      \begin{proof}[\bf Proof]
        By testing the $u_m$-equation of \eqref{problema_P_m_intro} by $(u_m)_-$ and using that $\norm{a_m((u_m)_-)}{} \leq \norm{(u_m)_-}{}$, we obtain
        \begin{align*}
          & \frac{1}{2} \frac{d}{dt} \norma{(u_m)_-}{L^2(\Omega)}^2 + \norma{\nabla (u_m)_-}{L^2(\Omega)}^2 = \int_{\Omega}{a_m((u_m)_-) \nabla v_m \cdot \nabla (u_m)_- \ dx} \\
          & \leq \norma{(u_m)_-}{L^3(\Omega)} \norma{\nabla v_m}{L^6(\Omega)} \norma{\nabla (u_m)_-}{L^2(\Omega)} \\
          & \leq C \norma{\nabla v_m}{L^6(\Omega)} \Big( \norma{(u_m)_-}{L^2(\Omega)}  \norma{\nabla (u_m)_-}{L^2(\Omega)} + C \norma{(u_m)_-}{L^2(\Omega)}^{1/2}  \norma{\nabla (u_m)_-}{L^2(\Omega)}^{3/2}\Big).
        \end{align*}
        Hence, using Young's inequality and that $v_m \in L^{\infty}(0,T;W^{1,6}(\Omega))$ we can arrive at
        \begin{align*}
          \frac{1}{2} \frac{d}{dt} \norma{(u_m)_-}{L^2(\Omega)}^2 + \frac{1}{2} \norma{\nabla (u_m)_-}{L^2(\Omega)}^2 & \leq C \norma{(u_m)_-}{L^2(\Omega)}^2.
        \end{align*}
        Note that $(u^0_m)_- = 0$, by hypothesis. Therefore, if we apply Gronwall's inequality (Lemma \ref{lema_gronwall}) we conclude that $(u_m)_-(t,x) = 0$, $a.e. (t,x) \in (0,T) \times \Omega$, for all $T > 0$, that is, $u_m(t,x) \geq 0$, $a.e. (t,x) \in (0,\infty) \times \Omega$.
      
      \
      
        In order to establish the positivity of $v_m$ we define the function $V(t) = \D{\min_{x \in \Omega}{\{ v^0_m(x) \} }} \, \exp(- m^s t)$. Clearly, $V$ is a sub-solution of the $v_m$-equation of \eqref{problema_P_m_intro}, because $-a_m(u_m)^s \geq -m^s$. In fact, we have
        \begin{equation*}
          V_t(t) - \Delta V(t) = - m^s V(t) \leq - a_m(u_m(t,x))^s V(t),
        \end{equation*}
        Comparing $V$ and $v_m$ we conclude that $v_m(t,x) \geq V(t) \ a.e.$ in $(0,\infty) \times \Omega$.       
        The upper bound on $v_m$ can be obtained an analogous argument, but now using the super-solution $V(t,x) = \norma{v^0_m}{L^{\infty}(\Omega)}$. This implies that
        \begin{equation*}
          \D{\essinf_{x \in \Omega}{\{ v^0_m(x) \} }} \, \exp(-m^s t) \leq v_m(t,x) \leq \norma{v^0_m}{L^{\infty}(\Omega)} \ a.e. (t,x) \in (0,\infty) \times \Omega,
        \end{equation*}
        and using \eqref{convergencia_dado_inicial_v_m} finally leads us to the desired result.
      \end{proof}
      
    \begin{lema}{\bf ($\boldsymbol{m}$-uniform estimates for $\boldsymbol{(u_m,v_m)}$)}
      \begin{enumerate}
        \item For every $t \geq 0$,
          \begin{equation*}
            \norma{u_m(t)}{L^1(\Omega)} = \int_{\Omega}{u_m(t,x) \ dx} = \norma{u^0_m}{L^1(\Omega)} = \norma{u^0}{L^1(\Omega)} = u^{\ast} \norm{\Omega}{};
          \end{equation*}
        \item For every $t > 0$,
        \begin{equation*}
          \norma{v_m(t)}{L^2(\Omega)}^2 + 2 \int_0^t{\norma{\nabla v_m(s)}{L^2(\Omega)}^2 \ ds} + \int_0^t{\norma{a_m(u_m(s))^{s/2} v_m(s)}{L^2(\Omega)}^2 \ ds} \leq \norma{v^0}{L^2(\Omega)}^2,
        \end{equation*}
        which allows us to conclude in particular that
        \begin{equation}
          \nabla v_m \mbox{ is bounded in } L^2(0,\infty;L^2(\Omega)), \mbox{ independently of } m \in \mathbb{N}.
          \label{limitacao_uniforme_nabla_v_m}
        \end{equation}
      \end{enumerate}
      \label{lemma_limitacao_u_v_independente_m}
    \end{lema}
    \begin{proof}[\bf Proof]
      Taking \eqref{convergencia_dado_inicial_u_m} and \eqref{convergencia_dado_inicial_v_m} into account, to prove the first item, we integrate the the $u_m$-equation of \eqref{problema_P_m_intro} and take into account that $u_m \geq 0$, thanks to Lemma \ref{lemma_positividade_u_v_m}. The second item can be proved by testing the the $v_m$-equation of \eqref{problema_P_m_intro} by $2v_m$.
    \end{proof}
    \begin{lema}{\bf (Weak regularity of $\boldsymbol{(u_m,v_m)}$ up to infinity time)}
      For each fixed $m \in \mathbb{N}$, the following regularity at infinity time to $u$ and $v$ holds:
      \begin{equation*}
        (u_m - u^{\ast}), v_m \in L^{\infty}(0,\infty;L^2(\Omega)) \cap L^2(0, \infty; H^1(\Omega)).
      \end{equation*}
      \label{lema_comportamento_assintotico_u_v_m}
    \end{lema}
    \begin{proof}[\bf Proof]
      From Lemma \ref{lemma_limitacao_u_v_independente_m}.1 we have $ \displaystyle  \int_{\Omega}{(u_m(t) - u^{\ast}) \ dx} = 0, \ \forall t \in [0,\infty)$, that is, $u_m - u^{\ast}$ is a null mean function. Besides, by testing the the $u_m$-equation of \eqref{problema_P_m_intro} by $2 u_m$ and using Lemma \ref{lemma_limitacao_u_v_independente_m}.2, we arrive at
      \begin{align*}
        \norma{u_m(t)}{L^2(\Omega)}^2 + \int_0^t{\norma{\nabla u_m(s)}{L^2(\Omega)}^2 \ ds} & \leq \norma{u^0_m}{L^2(\Omega)}^2 + m^2 \int_0^t{\norma{\nabla v_m(s)}{L^2(\Omega)}^2 \ ds} \\
        & \leq \norma{u^0_m}{L^2(\Omega)}^2 + \frac{m^2}{2} \norma{v^0_m}{L^2(\Omega)}^2.
      \end{align*}
      The latter allows us to conclude that
      \begin{equation*}
        u_m \in L^{\infty}(0,\infty;L^2(\Omega)) \mbox{ and } \nabla u_m \in L^2(0,\infty;L^2(\Omega)).
      \end{equation*}
      Hence, using the Poincaré's type inequality of Lemma \ref{lema_desigualdade_poincare_media_nula} we can prove that
      \begin{equation}
        u_m - u^{\ast} \in L^{\infty}(0,\infty;L^2(\Omega)) \cap L^2(0, \infty; H^1(\Omega)).
        \label{regularidade_tempo_infinito_u_m}
      \end{equation}
      
      Next, we use this fact and take $v$ as a test function in the following reformulation of the $v$-equation of \eqref{problema_P_m_intro}
      \begin{equation*}
        (v_m)_t - \Delta v_m + a_m(u^\ast)^s v_m = - (a_m(u_m)^s - a_m(u^{\ast})^s) v_m,
      \end{equation*}
      we obtain
      \begin{equation*}
        \frac{1}{2} \frac{d}{dt} \norma{v_m(t)}{L^2(\Omega)}^2 + \norma{\nabla v_m(t)}{L^2(\Omega)}^2 + a_m(u^{\ast})^s \norma{v_m(t)}{L^2(\Omega)}^2 = - \int_{\Omega}{(a_m(u_m(t,x))^s - a_m(u^{\ast})^s) v_m(t,x)^2 \ dx}.
      \end{equation*}
      Using Lemma \ref{lema_a_m_elevado_a_s}, the right hand side can be estimated by
      \begin{align*}
        \norm{\int_{\Omega}{(a_m(u_m)^s - a_m(u^{\ast})^s) v^2 \ dx}}{} & \leq \norma{v^0_m}{L^{\infty}(\Omega)} \int_{\Omega}{ \norm{a_m(u_m)^s - a_m(u^{\ast})^s}{} v_m \ dx} \\
        & \leq \norma{v^0}{L^{\infty}(\Omega)} \int_{\Omega}{ \norm{a_m(u_m) + a_m(u^{\ast})}{}^{s-1} \norm{a_m(u_m) - a_m(u^{\ast})}{} v \ dx} \\
        & \leq 2^{s-1} m^{s-1} \norma{v^0}{L^{\infty}(\Omega)} \norma{u_m - u^{\ast}}{L^2(\Omega)} \norma{v_m}{L^2(\Omega)} \\
        & \leq C(u^\ast,v^0) m^{s-1} \norma{u_m - u^{\ast}}{L^2(\Omega)}^2 + \frac{a_m(u^{\ast})^s}{2} \norma{v_m}{L^2(\Omega)}^2,
      \end{align*}
      where $ C(u^\ast,v^0) > 0$ is a constant (independent of $t$ and $x$). Now, considering also the terms of the left hand side, one has
      \begin{equation*}
        \frac{1}{2} \frac{d}{dt} \norma{v_m(t)}{L^2(\Omega)}^2 + \norma{\nabla v_m(t)}{L^2(\Omega)}^2 + a_m(u^{\ast}) \norma{v_m(t)}{L^2(\Omega)}^2 \leq C(u^{\ast},v^0) m^{s-1} \norma{u_m(t) - u^{\ast}}{L^2(\Omega)}^2.
      \end{equation*}
      Note that $a_m(u^{\ast})$ is a fixed positive real number if $u^0 \neq 0$ and $\norma{u_m(t) - u^{\ast}}{L^2(\Omega)}^2 \in L^1(0,\infty)$, because of \eqref{regularidade_tempo_infinito_u_m}. Hence we can conclude that $v_m \in L^2(0,\infty;L^2(\Omega))$ and, together with \eqref{limitacao_uniforme_nabla_v_m}, we finally conclude that $v_m \in L^2(0,\infty;H^1(\Omega))$.
    \end{proof}
    
    The regularity given in Lemma \ref{lema_comportamento_assintotico_u_v_m} allows us to obtain
    the regularity \eqref{regularidade_solucao_u_m_v_m_reg_H_2}. In fact, first we test the $v_m$-equation of \eqref{problema_P_m_intro} by $-\Delta v_m \in L^2_{loc}([0,\infty);H^1(\Omega))$. After some computations, we arrive at
    \begin{align*}
      & \frac{1}{2} \frac{d}{dt} \norma{\nabla v_m}{L^2(\Omega)}^2 + \norma{\Delta v_m}{L^2(\Omega)}^2 + \norma{a_m(u_m)^{s/2} \nabla v_m(s)}{L^2(\Omega)}^2 \\
      & \qquad \leq C(m,\norma{v^0}{L^{\infty}(\Omega)}) (\norma{\nabla u_m}{L^2(\Omega)}^2 + \norma{\nabla v_m}{L^2(\Omega)}^2),
    \end{align*}
    and this allows us to conclude that
    \begin{equation}
      v_m \in L^{\infty}(0,\infty;H^1(\Omega)) \cap L^2(0,\infty;H^2(\Omega))
      \label{regularidade_tempo_infinito_v_m_forte}
    \end{equation}
    because, after Lemmas \ref{lemma_limitacao_u_v_independente_m} and \ref{lema_comportamento_assintotico_u_v_m} we have $\norma{\nabla u_m}{L^2(\Omega)}^2, \norma{\nabla v_m}{L^2(\Omega)}^2 \in L^1(0,\infty)$. Next we take the gradient of the $v_m$-equation of \eqref{problema_P_m_intro} and test the resulting equation by $- \nabla \Delta v_m \in L^2_{loc}([0,\infty);L^2(\Omega))$, obtaining
    \begin{align*}
      & \frac{1}{2} \frac{d}{dt} \norma{\Delta v_m}{L^2(\Omega)}^2 + \norma{\nabla \Delta v_m}{L^2(\Omega)}^2 \\
      & \qquad = s \int_{\Omega}{a_m'(u_m) a_m(u_m)^{s-1} v_m \nabla u_m \cdot \nabla \Delta v_m \ dx} + \int_{\Omega}{a_m(u_m)^s \nabla v_m \cdot \nabla \Delta v_m \ dx} \\
      & \qquad \leq C(s,m,\norma{v^0}{L^{\infty}(\Omega)}) (\norma{\nabla u_m}{L^2(\Omega)}^2 + \norma{\nabla v_m}{L^2(\Omega)}^2) + \frac{1}{2} \norma{\nabla \Delta v_m}{L^2(\Omega)}^2,
    \end{align*}
    which gives us
    \begin{equation}
      v_m \in L^{\infty}(0,\infty;H^2(\Omega)) \mbox{ and } \Delta v_m \in L^2(0,\infty;H^1(\Omega)).
      \label{regularidade_tempo_infinito_v_m_mais_forte}
    \end{equation}
    After regularities \eqref{regularidade_tempo_infinito_v_m_forte} and \eqref{regularidade_tempo_infinito_v_m_mais_forte}, we can go back to the $v$-equation of \eqref{problema_P_m_intro} and conclude that
    \begin{equation*}
      (v_m)_t \in L^2(0,\infty;H^1(\Omega)).
    \end{equation*}
    
    Now, we can test the $u_m$-equation of \eqref{problema_P_m_intro} by $-\Delta u_m \in L^2_{loc}([0,\infty);L^2(\Omega))$. Considering the bounds of the truncation $a_m$ given in \eqref{limitacao_do_truncamento_a_m}, the interpolation inequality of Lemma \ref{lema_desig_lady2}-2 and inequality \eqref{obsv_norma_de_nabla_v}, we have
    \begin{align*}
      \frac{d}{dt} \norma{\nabla u_m}{L^2(\Omega)}^2 + \norma{\Delta u_m}{L^2(\Omega)}^2 & = \int_{\Omega}{a_m'(u_m) \nabla u_m \cdot \nabla v_m \Delta u_m \ dx} + \int_{\Omega}{a_m(u_m) \Delta v_m \Delta u_m \ dx} \\
      & \leq C \norma{\nabla u_m}{L^3(\Omega)} \norma{\nabla v_m}{L^6(\Omega)} \norma{\Delta u_m}{L^2(\Omega)} + m \norma{\Delta v_m}{L^2(\Omega)} \norma{\Delta u_m}{L^2(\Omega)} \\
      & \leq C \norma{\nabla u_m}{L^2(\Omega)}^{1/2} \norma{\nabla v_m}{L^6(\Omega)} \norma{\Delta u_m}{L^2(\Omega)}^{3/2} + m \norma{\Delta v}{L^2(\Omega)} \norma{\Delta u_m}{L^2(\Omega)} \\
      & \leq C \norma{\nabla v_m}{L^6(\Omega)}^4 \norma{\nabla u_m}{L^2(\Omega)}^2 + C(m) \norma{\Delta v_m}{L^2(\Omega)}^2 + \frac{1}{2} \norma{\Delta u_m}{L^2(\Omega)}^2.
    \end{align*}
    After absorbing the term $\dfrac{1}{2} \norma{\Delta u_m}{L^2(\Omega)}^2$, the other terms in the right hand side of the inequality belong to $L^1(0,\infty)$. Hence, integrating the last inequality with respect to $t$,
    \begin{equation*}
      \nabla u_m \in L^{\infty}(0,\infty;L^2(\Omega)) \mbox{ and } \Delta u_m \in L^2(0;\infty;L^2(\Omega)),
    \end{equation*}
    finishing the proof of the regularity \eqref{regularidade_solucao_u_m_v_m_reg_H_2}.
    
    Finally, we consider the case in which we assume that the Poisson-Neumann problem \eqref{Neumann_problem} has the $W^{2,3}$-regularity. In this case, as observed in the end of Subsection \ref{Subsection:regularized_problem_approx_solutions}, the regularity \eqref{regularidade_solucao_u_m_v_m_reg_W_2_3} already holds if we consider finite intervals $(0,T)$, for finite $T > 0$, instead of $(0,\infty)$. Then we can take the gradient of the $u$-equation of \eqref{problema_P_m_intro} and test the resulting equation by $- \nabla \Delta u_m \in L^2_{loc}([0,\infty);L^2(\Omega))$, obtaining
    \begin{align*}
      &  \frac{1}{2} \dfrac{d}{dt} \norma{\Delta u_m}{L^2(\Omega)}^2  + \norma{\nabla \Delta u_m}{L^2(\Omega)}^2 = \int_{\Omega}{a_m(u_m) \nabla \Delta v_m \cdot \nabla \Delta u_m \ dx} \\
      & + \int_{\Omega}{a_m'(u_m) \Delta v_m \nabla u_m \cdot \nabla \Delta u_m \ dx} + \int_{\Omega}{a_m'(u_m) \nabla u_m \cdot D^2 v_m \nabla \Delta u_m \ dx} \\
      & + \int_{\Omega}{a_m'(u_m) \nabla v_m \cdot D^2 u_m \nabla \Delta u_m \ dx} + \int_{\Omega}{a_m''(u_m) \Big ( \nabla u_m \cdot \nabla v_m \Big ) \Big ( \nabla u_m \cdot \nabla \Delta u_m \Big ) \ dx}.
    \end{align*}
    We recall that, from \eqref{limitacao_do_truncamento_a_m}, we have $\norm{a_m(u_m)}{} \leq m$ and $\norm{a_m'(u_m)}{}, \norm{a_m''(u_m)}{} \leq C$, for some $C > 0$. Then, using Hölder's inequality we obtain
    \begin{align*}    
      & \frac{1}{2} \dfrac{d}{dt} \norma{\Delta u_m}{L^2(\Omega)}^2 + \norma{\nabla \Delta u_m}{L^2(\Omega)}^2 \leq \Big ( m \norma{\nabla \Delta v_m}{L^2(\Omega)} \\
      & + C \norma{\Delta v_m}{L^3(\Omega)} \norma{\nabla u_m}{L^6(\Omega)} + C \norma{D^2 v_m}{L^3(\Omega)} \norma{\nabla u_m}{L^6(\Omega)} \\
      & + C \norma{\nabla v_m}{L^6(\Omega)} \norma{D^2 u_m}{L^3(\Omega)} + C \norma{\nabla u_m}{L^6(\Omega)}^2 \norma{\nabla v_m}{L^6(\Omega)} \Big ) \norma{\nabla \Delta u_m}{L^2(\Omega)}
    \end{align*}
    Using the continuous embedding $L^6(\Omega) \subset H^1(\Omega)$, the $W^{2,3}$-regularity, inequality \eqref{obsv_norma_de_nabla_v}, the interpolation inequality for the the $L^3$-norm (Lemma \ref{lema_desig_lady2}) in $3D$  domains and Young's inequality, we arrive at
    \begin{align*}
      \dfrac{d}{dt} \norma{\Delta u_m}{L^2(\Omega)}^2 & + \frac{1}{2} \norma{\nabla \Delta u_m}{L^2(\Omega)}^2 \leq C m \norma{\nabla \Delta v_m}{L^2(\Omega)}^2 \\
      & \quad + C \Big [ \norma{\Delta v_m}{L^3(\Omega)}^2 + \norma{\Delta v_m}{L^2(\Omega)}^4 + \norma{\Delta u_m}{L^2(\Omega)}^2 \norma{\Delta v_m}{L^2(\Omega)}^2 \Big ] \norma{\Delta u_m}{L^2(\Omega)}^2.
    \end{align*}
    Now notice that the first term of the right hand side of the last inequality belongs to $L^1(0,\infty)$ and, because of the regularity obtained so far, we have $C \Big [ \norma{\Delta v_m}{L^3(\Omega)}^2 + \norma{\Delta v_m}{L^2(\Omega)}^2 + \norma{\Delta u_m}{L^2(\Omega)}^4 \norma{\Delta v_m}{L^2(\Omega)}^2 \Big ] \in L^1(0,\infty)$. Therefore, using Grownwall's inequality (Lemma \ref{lema_gronwall}) we can conclude that
    \begin{equation*}
      \Delta u_m \in L^{\infty}(0,\infty;L^2(\Omega)) \cap L^2(0,\infty;H^1(\Omega))
    \end{equation*}
    and hence
    \begin{equation*}
      u_m \in L^{\infty}(0,\infty;H^2(\Omega)).
    \end{equation*}
    Next we can go back to the $u_m$-equation of \eqref{problema_P_m_intro} and conclude that $(u_m)_t \in L^2(0,\infty;H^1(\Omega))$, finishing the proof of the regularity \eqref{regularidade_solucao_u_m_v_m_reg_W_2_3}.
    

\section{Proof of Theorem \ref{teo_3D_existencia}}
    
    \label{section: original problem 3D}
    
    In this section we will obtain $m$-independent estimates to $(u_m,v_m)$, the solution of problem \eqref{problema_P_m_intro}, in order to pass to the limit as $m \to \infty$ and prove the existence of solution to the original problem \eqref{problema_P}.
    
        \subsection{An energy law appears: formal computations}
    \label{subsection: formal computations}
    
    The basic idea to obtain additional \emph{a priori} $m$-independent estimates is that the effects of the consumption  and chemotaxis terms cancel. First of all, we present some formal calculations to illustrate how it works. Suppose $(u,v)$ is a regular enough solution to the original problem \eqref{problema_P} with $u,v > 0$. Consider the change o variable $z = \sqrt{v}$, then \eqref{problema_P} can be rewritten as
    \begin{equation}\tag{P$_{(u,z)}$}
      \begin{array}{rl}
        \partial_t u - \Delta u & = - \nabla \cdot (u \nabla (z)^2 ) \\
        \partial_t z - \Delta z - \dfrac{\norm{\nabla z}{}^2}{z} & = - \dfrac{u^s z}{2} \\
        \partial_\eta u \Big |_{\Gamma} & =  \partial_\eta z \Big |_{\Gamma} = 0 \\
        u(0) & = u^0, \quad z(0) = \sqrt{v^0},
      \end{array}
      \label{problema_P_u_z} \\
    \end{equation}
    We are going to obtain estimates for $u$ and $z$ and then extract estimates for $v$ from the estimates of $z$.
    
    For this, we consider a function $g(u)$ such that $g''(u) = u^{s-2}$. Formally, assuming $u,z > 0$ we can use 
    \begin{equation}
      g'(u) = \left \{
      \begin{array}{rl}
        \dfrac{u^{s-1}}{(s-1)} &, \mbox{ if } s > 1, \\[12pt]
        ln(u) &, \mbox{ if } s = 1.
      \end{array}
      \right.
      \label{funcao_teste_formal}
    \end{equation}
    as a test function in the $u$-equation of \eqref{problema_P_u_z}, obtaining
    \begin{equation*}
      \dfrac{d}{dt} \int_{\Omega}{g(u) \ dx} + \int_{\Omega}{g''(u) \norm{\nabla u}{}^2 \ dx} = \int_{\Omega}{u g''(u) \nabla (z^2) \cdot \nabla u}
    \end{equation*}
    and, since $u g''(u) = u^{s-1}$, we have
    \begin{equation} \label{estimativa_formal_u}
      \dfrac{d}{dt} \int_{\Omega}{g(u) \ dx} + \int_{\Omega}{u^{s-2} \norm{\nabla u}{}^2 \ dx} = \int_{\Omega}{u^{s-1} \nabla (z^2) \cdot \nabla u \ dx} = \frac{1}{s} \int_{\Omega}{\nabla (z^2) \cdot \nabla (u^s) \ dx}.
    \end{equation}
    On the other hand, we can test the $z$-equation of \eqref{problema_P_u_z} by $- \Delta z$. Then we obtain
    \begin{equation} \label{estimativa_formal_z}
      \begin{array}{rl}
        \dfrac{1}{2} \dfrac{d}{dt} \D{\int_{\Omega}}{\norm{\nabla z}{}^2 \ dx} & + \D{\int_{\Omega}}{\norm{\Delta z}{}^2 \ dx} + \D{\int_{\Omega}}{\dfrac{\norm{\nabla z}{}^2}{z} \Delta z \ dx} + \dfrac{1}{2} \D{\int_{\Omega}}{u^s \norm{\nabla z}{}^2 \ dx} \\[12pt] 
        & = - \dfrac{1}{4} \D{\int_{\Omega}}{\nabla (u^s) \cdot \nabla (z^2) \ dx}.
      \end{array}
    \end{equation}
    Hence, if we add \eqref{estimativa_formal_z} to $s/4$ times \eqref{estimativa_formal_u}, then the two terms on the right hand side cancel each other and we obtain the time differential equation
    \begin{equation} \label{eq_formal_u_w}
      \begin{array}{rl}
      &  \dfrac{d}{dt} \left[ \dfrac{s}{4} \D{\int_{\Omega}{g(u) \ dx}}  + \dfrac{1}{2} \D{\int_{\Omega}{\norm{\nabla z}{}^2 \ dx}} \right] + \dfrac{s}{4} \D{\int_{\Omega}{u^{s-2} \norm{\nabla u}{}^2 \ dx}}  \\[12pt]
        & \qquad + \dfrac{1}{2} \D{\int_{\Omega}{u^s \norm{\nabla z}{}^2 \ dx} + \D{\int_{\Omega}{\norm{\Delta z}{}^2 \ dx}} + \int_{\Omega}{\dfrac{\norm{\nabla z}{}^2}{z} \Delta z \ dx}} = 0.
      \end{array}
    \end{equation}
    The main idea now is to estimate from below the term
    \begin{equation} \label{termo_problematico}
      \D{\int_{\Omega}{\norm{\Delta z}{}^2 \ dx}} + \int_{\Omega}{\dfrac{\norm{\nabla z}{}^2}{z} \Delta z \ dx}
    \end{equation}
    (see Lemma \ref{lema_termo_fonte_final} below).

    
    \subsection{Rigorous justification of the energy inequalities}
    \label{subsection: treatment of the problematic term}
    
    In the sequel, we consider the regularized problem \eqref{problema_P_m_intro} and its solution $(u_m,v_m)$ instead of the original problem \eqref{problema_P} and $(u,v)$. In this case, we have to deal with the truncation $a_m(\cdot)$ in the chemotaxis and consumption terms and with the fact that $u_m$ is nonnegative, but not necessarily strictly positive (Lemma \ref{lemma_positividade_u_v_m}). In order to obtain time independent estimates, we will also need to separate $z$ from zero, that's why we are going to consider the change of variables $z = \sqrt{v + \alpha}$, for $\alpha > 0$ to be chosen later, instead of $\sqrt{v}$. With this modification, we will obtain the corresponding version of \eqref{eq_formal_u_w}. We will separate the cases $s = 1$, $s \in (1,2)$ and $s \geq 2$. Note that if $s = 1$ then $g'(u)$, given by \eqref{funcao_teste_formal}, and $g''(u)$ have a singularity at $u = 0$. If $s \in (1,2)$, then only $g''(u)$ is singular at $u = 0$ and, if $s \geq 2$, then neither $g'(u)$ nor $g''(u)$ are singular.
    
    Let us consider the variable $z_m(t,x) = \sqrt{v_m(t,x) + \alpha}$. Taking into account that the pair $(u_m,v_m)$ is a strong solution of \eqref{problema_P_m_intro}, on has by straightforward calculations that $(u_m,z_m)$ satisfies the following equivalent problem:
    \begin{equation}\tag{Q$_m$}
      \begin{array}{rl}
        \partial_t u_m - \Delta u_m & = - \nabla \cdot (a_m(u_m) \nabla (z_m)^2) \\
        \partial_t z_m - \Delta z_m - \dfrac{\norm{\nabla z_m}{}^2}{z_m} & 
        = - \dfrac12 a_m(u_m)^s z_m + \dfrac{\alpha}2 \dfrac{a_m(u_m)^s}{ z_m} \\
        \partial_\eta u_m \Big |_{\Gamma} & =  \partial_\eta z_m \Big |_{\Gamma} = 0 \\
        u_m(0) & = u^0, \quad z_m(0) = \sqrt{v^0 + \alpha}.
      \end{array}
      \label{problema_P_u_m_z_m} \\
    \end{equation}
    
    In the present subsection, we drop the $m$-subscript and write $(u,z)$ for $(u_m,z_m)$ to simplify the notation along the proofs of the forthcoming Lemmas. We remark that all the constants obtained in these Lemmas are independent of the parameter $m \in \mathbb{N}$ and this is why the energy inequalities proved in this section allow us to obtain $m$-independent bounds in the rest of the paper.
    
    \
    
    We use the following lemma in order to estimate \eqref{termo_problematico}, whose proof can be found in the Appendix B.
    
    \begin{lema}.  \label{lema_termo_fonte_final}
      Suppose that the Poisson-Neumann problem \eqref{Neumann_problem} has the $H^2$-regularity and assume that Hypothesis \ref{hypothesis_density} holds. Then there exist positive constants $C_1, C_2 > 0$ such that
      \begin{align*}
        \int_{\Omega}{\norm{\Delta z}{}^2 \ dx} + \int_{\Omega}{\frac{\norm{\nabla z}{}^2}{z} \Delta z \ dx} & \geq C_1 \Big ( \int_{\Omega}{\norm{D^2 z}{}^2 \ dx} + \int_{\Omega}{\frac{\norm{\nabla z}{}^4}{z^2} \ dx} \Big ) - C_2 \int_{\Omega}{\norm{\nabla z}{}^2 \ dx},
      \end{align*}
      for all $z \in H^2(\Omega)$ such that $\partial_\eta z \Big |_{\Gamma} = 0$ and $z \geq \alpha$, for some $\alpha > 0$.
    \end{lema}
    
    Now we prove the following.    
    \begin{lema}
      The solution $(u,z)$ of \eqref{problema_P_u_m_z_m}, satisfies the inequality
      \begin{equation*}
        \begin{array}{c}
          \dfrac{1}{2} \dfrac{d}{dt} \norma{\nabla z}{L^2(\Omega)}^2 + C_1 \Big ( \D{\int_{\Omega}}{\norm{D^2 z}{}^2 \ dx} + \D{\int_{\Omega}}{\frac{\norm{\nabla z}{}^4}{z^2} \ dx} \Big ) + \dfrac{1}{2} \D{\int_{\Omega}{a_m(u)^s \norm{\nabla z}{}^2 \ dx}} \\[12pt]
          \leq \dfrac{s}{4} \D{\int_{\Omega}}{a_m(u)^{s-1} \nabla (z^2) \cdot \nabla a_m(u) \ dx} + \dfrac{s}{2} \sqrt{\alpha} \D{\int_{\Omega}}{a_m(u)^{s-1} \norm{\nabla z}{} \norm{\nabla a_m(u)}{} \ dx} + C_2 \int_{\Omega}{\norm{\nabla z}{}^2 \ dx}.
        \end{array}
      \end{equation*}
      \label{lema_tratamento_equacao_z}
    \end{lema}
    \begin{proof}[\bf Proof]
      We begin by testing the $z$-equation of \eqref{problema_P_u_m_z_m} by $- \Delta z$. This gives us
      \begin{align*}
       & \frac{1}{2} \frac{d}{dt} \norma{\nabla z}{L^2(\Omega)}^2 + \norma{\Delta z}{L^2(\Omega)}^2 + \int_{\Omega}{\frac{\norm{\nabla z}{}^2}{z} \Delta z \ dx} + \frac{1}{2} \int_{\Omega}{a_m(u)^s \norm{\nabla z}{}^2 \ dx} \\
       & \leq \frac{s}{4} \int_{\Omega}{a_m(u)^{s-1} \nabla (z^2) \cdot \nabla a_m(u) \ dx} + \frac{s}{2} \sqrt{\alpha} \int_{\Omega}{a_m(u)^{s-1} \norm{\nabla z}{} \norm{\nabla a_m(u)}{} \ dx}.
      \end{align*}
      Then, applying Lemma \ref{lema_termo_fonte_final}, we obtain the desired inequality.
    \end{proof}
    
    Define, for $\sigma \geq 0$, the functions $g_m$ and $g_{m,j}$, adequate regularizations of the function $g$ that appears in the formal inequality \eqref{eq_formal_u_w}, by
    \begin{equation} \label{definicao_g_m_e_g_m_j}
      g_m(\sigma) = \int_0^{\sigma}{g_m'(r) \ dr} \quad \mbox{ and } \quad g_{m,j}(\sigma) = \int_0^{\sigma}{g_{m,j}'(r) \ dr},
    \end{equation}
    where $g_m'$ and $g_{m,j}'$ are defined for $r \geq 0$ and given by
    \begin{equation} \label{definicao_g_m'}
      g_m'(r) = \left \{
      \begin{array}{rl}
        ln(a_m(r) + 1), & \mbox{if } s = 1, \\
        \dfrac{a_m(r)^{s-1}}{s - 1}, & \mbox{if } s > 1,
      \end{array}
      \right.
    \end{equation}
    and
    \begin{equation} \label{definicao_g_m_j'}
      g_{m,j}'(r) = \dfrac{(a_m(r) + 1/j)^{s-1}}{s - 1}, \mbox{ for } s \in (1,2) \mbox{ and } j \in \mathbb{N}.
    \end{equation}
    In Lemmas \ref{lemma_estimativa_u_v_m_1_s=1}, \ref{lemma_estimativa_u_v_m_1_s_intermediario} and \ref{lemma_estimativa_u_v_m_1_s_geq_2} below we will be interested in the terms
    \begin{equation*}
      \int_{\Omega}{u_t(t,x) \ g_m'(u(t,x)) \ dx} \quad \mbox{ and } \quad \int_{\Omega}{u_t(t,x) \ g_{m,j}'(u(t,x)) \ dx},
    \end{equation*}
    recalling that $u(t,x)$ denotes the function $u_m(t,x)$ of the solution $(u_m,z_m)$ of \eqref{problema_P_u_m_z_m}. We have the following results on these terms.
    \begin{teo} \label{teo_derivada_composicao_com_g}
      The weak time derivatives of $g_{m,j}(u)$ and $g_m(u)$ belong to $L^2(0,\infty;L^2(\Omega))$ and are given by
      \begin{equation}
        \partial_t g_m(u) = g_m'(u) \ u_t
        \label{derivada_de_g_m}
      \end{equation}
      and
      \begin{equation}
        \partial_t g_{m,j}(u) = g_{m,j}'(u) \ u_t.
        \label{derivada_de_g_m_j}
      \end{equation}
    \end{teo}
    \begin{proof}[\bf Proof]
      We are going to prove \eqref{derivada_de_g_m_j}. The proof of \eqref{derivada_de_g_m} is analogous to the proof of \eqref{derivada_de_g_m_j}. Because of Lemma \ref{lema_densidade_Lions_Magenes}, we know that there is a sequence $(u^n)_{n \in \mathbb{N}} \subset C^{\infty}_c(0,\infty;L^2(\Omega))$ that converges to $u$ in the norm of $H^1((0,\infty);L^2(\Omega))$. For these functions $u^n$, which are very regular in the time variable, we can write
      \begin{align*}
        - \int_0^{\infty}{g_{m,j}(u^n) \ \varphi'(t) \ dt} = \int_0^{\infty}{g_{m,j}'(u^n) u^n_t \ \varphi(t) \ dt} = \int_0^{\infty}{\frac{(a_m(u^n) + 1/j)^{s-1}}{s-1} \ u^n_t \ \varphi(t) \ dt},
      \end{align*}
      for all $\varphi \in C^{\infty}_c(0,\infty)$. Then we note that the convergence in the norm of $H^1((0,{\infty});L^2(\Omega))$ is enough to pass to the limit as $n \to \infty$ in the first and in the last terms of the above equality, yielding
      \begin{equation*}
        - \int_0^{\infty}{g_{m,j}(u) \ \varphi'(t) \ dt} = \int_0^{\infty}{\frac{(a_m(u) + 1/j)^{s-1}}{s-1} \ u_t \ \varphi(t) \ dt}, \ \forall \varphi \in C^{\infty}_c(0,\infty).
      \end{equation*}
      By the definition of weak derivative, this is precisely \eqref{derivada_de_g_m_j}, as we wanted to prove.
    \end{proof}
    
    \begin{coro} \label{coro_expressao_derivada_integral_de_g_m_j}
      The weak time derivatives $\D{\frac{d}{dt} \int_{\Omega}{g_m(u(t,x)) \ dx}}$ and $\D{\frac{d}{dt} \int_{\Omega}{g_{m,j}(u(t,x)) \ dx}}$ belong to \newline $L^2(0,\infty)$ and are given by the expressions
      \begin{equation}
        \frac{d}{dt} \int_{\Omega}{g_m(u(t,x)) \ dx} = \int_{\Omega}{u_t(t,x) g_m'(u(t,x)) \ dx}
        \label{expressao_derivada_integral_de_g_m}
      \end{equation}
      and
      \begin{equation}
        \frac{d}{dt} \int_{\Omega}{g_{m,j}(u(t,x)) \ dx} = \int_{\Omega}{u_t(t,x) g_{m,j}'(u(t,x)) \ dx}.
        \label{expressao_derivada_integral_de_g_m_j}
      \end{equation}
    \end{coro}
    \begin{proof}[\bf Proof]
      We are going to prove \eqref{expressao_derivada_integral_de_g_m_j}. The proof of \eqref{expressao_derivada_integral_de_g_m} is analogous to the proof of \eqref{expressao_derivada_integral_de_g_m_j}. Accounting for Lemma \ref{lema_derivada_da_transformacao} and \eqref{derivada_de_g_m_j}, we finally obtain
      \begin{equation*}
        \frac{d}{dt} \int_{\Omega}{g_{m,j}(u(t,x)) \ dx} = \int_{\Omega}{\dfrac{d}{dt} g_{m,j}(u(t,x)) \ dx} = \int_{\Omega}{u_t(t,x) \frac{(a_m(u(t,x)) + 1/j)^{s-1}}{s-1} \ dx}.
      \end{equation*}
      Since $\frac{(a_m(u(t,x)) + 1/j)^{s-1}}{s-1}$ is pointwisely bounded and $u_t(t,x) \in L^2(0,\infty;L^2(\Omega))$ then 
      \begin{equation*}
        \frac{d}{dt} \int_{\Omega}{g_{m,j}(u(t,x)) \ dx} \in L^2(0,\infty).
      \end{equation*}
    \end{proof}
    
    Now, we are in position to prove an energy inequality associated to the formal inequality \eqref{eq_formal_u_w}.
    
    \begin{lema}[\bf Energy inequality for $\boldsymbol{s = 1}$]
        The solution $(u, z)$ of the problem \eqref{problema_P_u_m_z_m} satisfies, for sufficiently small $\alpha > 0$,
        \begin{equation}
          \begin{array}{c}
            \dfrac{d}{dt} \left[ \dfrac{1}{4} \D{\int_{\Omega}}{g_m(u) \ dx} + \dfrac{1}{2} \D{\int_{\Omega}}{\norm{\nabla z}{}^2 \ dx} \right]
            + C \D{\int_{\Omega}}{\norm{\nabla [a_m(u) + 1]^{1/2}}{}^2 \ dx} + \frac{1}{4} \D{\int_{\Omega}}{a_m(u) \norm{\nabla z}{}^2 \ dx} \\[12pt]
            + C_1 \Big ( \D{\int_{\Omega}}{\norm{D^2 z}{}^2 \ dx} + \D{\int_{\Omega}}{\frac{\norm{\nabla z}{}^4}{z^2} \ dx} \Big ) \leq C \int_{\Omega}{\norm{\nabla z}{}^2 \ dx},
          \end{array}
          \label{estimativa_u_v_m_1_s=1}
        \end{equation}
        where $g_m(u)$ is given by \eqref{definicao_g_m_e_g_m_j}.
        \label{lemma_estimativa_u_v_m_1_s=1}
      \end{lema}
      \begin{proof}[\bf Proof]
        In order to prove \eqref{estimativa_u_v_m_1_s=1} we will use the cancellation effect mentioned in Subsection \ref{subsection: formal computations}. Since now we are dealing with the regularized problem (instead of the original problem), we must pay attention to  two technical difficulties that arise:
        the presence of the truncation $a_m(\cdot)$ in the chemotaxis and consumption terms, the fact that $u$ is nonnegative, but not strictly positive, and now we consider $z = \sqrt{v + \alpha}$. For  $s = 1$, this means that instead of using $g'(u) = ln(u)$ as a test function in the $u$-equation of \eqref{problema_P_m_intro}, we must use $g'(a_m(u) + 1) = ln(a_m(u) + 1)$, in order to preserve the cancellation effect, avoid divisions by zero and invalid values for the argument of $ln(\cdot)$.
        
        We begin by using $\varphi = ln(a_m(u) + 1)$ in the $u$-equation of problem \eqref{problema_P_m_intro} to obtain
        \begin{equation*}
          \frac{d}{dt} \int_{\Omega}{g_m(u) \ dx} + \int_{\Omega}{\frac{a_m'(u)}{a_m(u) + 1} \norm{\nabla u}{}^2 \ dx} = \prodl{\frac{a_m(u)}{a_m(u) + 1} \nabla (z^2)}{\nabla a_m(u)},
        \end{equation*}
        where
        \begin{equation*}
          g_m(r) = \int_0^{r}{ln(a_m(\theta) + 1) \ d\theta}
        \end{equation*}
        is a primitive of $ln(a_m(r) + 1)$. 
        Due to the regularity of $u$, $u_t$ and the functions $g_m(r)$ and $g_m'(r)$, we can conclude that the weak derivative $ \frac{d}{dt} \int_{\Omega}{g_m(u) \ dx} =   \int_{\Omega}{g_m'(u) \, u_t\ dx} $ and belongs to $L^2(0,T)$.
        
        Since $0 \leq a_m'(u) \leq C$, we have $(a_m'(u))^2 \leq C a_m'(u)$, and we can write
        \begin{align*}
          \int_{\Omega}{\frac{a_m'(u)}{a_m(u) + 1} \norm{\nabla u}{}^2 dx} \geq C \int_{\Omega}{\frac{(a_m'(u))^2}{a_m(u) + 1} \norm{\nabla u}{}^2 dx} = C \int_{\Omega}{\frac{\norm{\nabla a_m(u)}{}^2}{a_m(u) + 1} \ dx} \geq C \int_{\Omega}{\norm{\nabla [a_m(u) + 1]^{1/2}}{}^2 dx}.
        \end{align*}
        
        Hence, using that $\dfrac{1}{(a_m(u) + 1)} \leq \dfrac{1}{\sqrt{a_m(u) + 1}}$, we have
        \begin{align*}
          \nonumber \frac{d}{dt} \int_{\Omega}{g_m(u) \ dx} & + C \int_{\Omega}{\norm{\nabla [a_m(u) + 1]^{1/2}}{}^2 \ dx} = 2 \prodl{\frac{a_m(u) + 1 - 1}{a_m(u) + 1} z \nabla z}{\nabla a_m(u)} \\
          & = \prodl{\nabla (z^2)}{\nabla a_m(u)} - 2 \prodl{z \nabla z}{\frac{\nabla a_m(u)}{a_m(u) + 1}} \\
          & \leq \prodl{\nabla (z^2)}{\nabla a_m(u)} + 2 \sqrt{\norma{v^0}{L^{\infty}(\Omega)} + \alpha} \norma{\nabla z}{L^2(\Omega)} \norma{\nabla [a_m(u) + 1]^{1/2}}{L^2(\Omega)}, \label{estimativa_u_m_1_s=1}
        \end{align*}
        Then we obtain
        \begin{equation}
          \frac{d}{dt} \int_{\Omega}{g_m(u) \ dx} + C \int_{\Omega}{\norm{\nabla [a_m(u) + 1]^{1/2}}{}^2 \ dx} \leq \prodl{\nabla (z^2)}{\nabla a_m(u)} + C \norma{\nabla z}{L^2(\Omega)}^2
          \label{estimativa_u_m_1_s=1}
        \end{equation}
        
        Now, using Lemma \ref{lema_tratamento_equacao_z} for $s = 1$ we obtain
        \begin{align*}
          \frac{1}{2} \frac{d}{dt} \norma{\nabla z}{L^2(\Omega)}^2 + C_1 \Big ( \int_{\Omega}{\norm{D^2 z}{}^2 \ dx} + \int_{\Omega}{\frac{\norm{\nabla z}{}^4}{z^2} \ dx} \Big ) + \frac{1}{2} \int_{\Omega}{a_m(u) \norm{\nabla z}{}^2 \ dx} \\
          \leq \frac{1}{4} \int_{\Omega}{\nabla (z^2) \cdot \nabla a_m(u) \ dx} + \frac{\sqrt{\alpha}}{2} \int_{\Omega}{\norm{\nabla z}{} \norm{\nabla a_m(u)}{} \ dx} + C_2 \norma{\nabla z}{L^2(\Omega)}^2.
        \end{align*}
        
        If we add the above inequality to $1/4$ times \eqref{estimativa_u_m_1_s=1}, then the terms $\D{\int_{\Omega}{\nabla a_m(u) \cdot \nabla (z^2) \ dx}}$ cancel and we obtain
        \begin{align}
          \nonumber \frac{d}{dt} \Big [ \frac{1}{4} \int_{\Omega}{g_m(u) \ dx} + \frac{1}{2} \norma{\nabla z}{L^2(\Omega)}^2 \Big ] + C \int_{\Omega}{\norm{\nabla [a_m(u) + 1]^{1/2}}{}^2 \ dx} + \frac{1}{2} \int_{\Omega}{a_m(u) \norm{\nabla z}{}^2 \ dx} \\
          + C_1 \Big ( \int_{\Omega}{\norm{D^2 z}{}^2 \ dx} + \int_{\Omega}{\frac{\norm{\nabla z}{}^4}{z^2} \ dx} \Big )  \leq \frac{\sqrt{\alpha}}{2} \int_{\Omega}{\norm{\nabla z}{} \norm{\nabla a_m(u)}{} \ dx} + C_2 \norma{\nabla z}{L^2(\Omega)}^2 \label{aux_inequality_u_z_s=1} \\
          \nonumber \leq \int_{\Omega}{\sqrt{\alpha} \norm{\nabla [a_m(u) + 1]^{1/2}}{} \norm{\sqrt{a_m(u) + 1}}{} \norm{\nabla z}{} \ dx} + C_2 \norma{\nabla z}{L^2(\Omega)}^2.
        \end{align}
        We can deal with the first term in the right hand side of the inequality using Hölder's and Young's inequality,
        \begin{align*}
       &   \int_{\Omega}{\sqrt{\alpha} \norm{\nabla [a_m(u) + 1]^{1/2}}{} \norm{\sqrt{a_m(u) + 1}}{} \norm{\nabla z}{} \ dx} \\
          \leq 
          & \alpha \ C(\delta) \int_{\Omega}{a_m(u) \norm{\nabla z}{}^2 \ dx} + \delta \norma{\nabla [a_m(u) + 1]^{1/2}}{L^2(\Omega)}^2 + \alpha \ C(\delta) \int_{\Omega}{\norm{\nabla z}{}^2 \ dx}.
        \end{align*}
        
        Therefore, we can first choose $\delta > 0$ and then $\alpha > 0$ sufficiently small in order to use the terms on the left hand side of inequality \eqref{aux_inequality_u_z_s=1} to absorb  the first two terms on the right hand side of the above inequality and finally obtain the desired inequality \eqref{estimativa_u_v_m_1_s=1}.
      \end{proof}
      
      \begin{lema}[\bfseries Energy inequality for $\boldsymbol{s \in (1,2)}$] \label{lemma_estimativa_u_v_m_1_s_intermediario}
        The solution $(u, z)$ of the problem \eqref{problema_P_u_m_z_m} satisfies, for sufficiently small $\alpha > 0$,
        \begin{equation} \label{estimativa_u_v_m_1_s_intermediario}
          \begin{array}{c}
            \dfrac{d}{dt} \left[ \dfrac{s}{4} \D{\int_{\Omega}}{g_m(u) \ dx} + \frac{1}{2} \norma{\nabla z}{L^2(\Omega)}^2 \right] + \dfrac{1}{4} \D{\int_{\Omega}}{a_{m}(u)^s \norm{\nabla z}{}^2 \ dx} \\[12pt]
            + C_1 \Big ( \D{\int_{\Omega}}{\norm{D^2 z}{}^2 \ dx} + \D{\int_{\Omega}}{\frac{\norm{\nabla z}{}^4}{z^2} \ dx} \Big ) \leq C \int_{\Omega}{\norm{\nabla z}{}^2 \ dx},
          \end{array}
        \end{equation}
        where $g_m(u)$ is given by \eqref{definicao_g_m_e_g_m_j}.
      \end{lema}
      \begin{proof}[\bfseries Proof]
        Analogously to the case $s = 1$ (Lemma \ref{lemma_estimativa_u_v_m_1_s=1}) in order to preserve the cancellation effect and avoid divisions by zero, for $s \in (1,2)$, instead of using $g'(u)$ as a test function in the $u$-equation of \eqref{problema_P_m_intro}, we should consider the sequence $\{ 1/j \}_{j \in \mathbb{N}}$ and use $g_{m,j}'(u)$ given by \eqref{definicao_g_m_j'}. Due to the complexity of the procedures that are involved, we divide the proof in three main steps: 
        \begin{enumerate}
          \item Obtain an inequality from the $u$-equation of problem \eqref{problema_P_u_m_z_m};
          \item Use this inequality and Lemma \ref{lema_tratamento_equacao_z} in order to obtain the corresponding version of \eqref{eq_formal_u_w};
          \item Pass to the limit as $j \to \infty$.
        \end{enumerate}
          
          \vspace{12pt}
          {\bfseries STEP 1:} In the first step we deal with the $u$-equation of \eqref{problema_P_m_intro}.
          
          We begin by testing the $u$-equation of \eqref{problema_P_u_m_z_m} by $g_{m,j}'(u) = (a_m(u) + 1/j)^{s-1}/(s-1)$ to obtain
        \begin{equation*}
          \dfrac{d}{dt} \int_{\Omega}{g_{m,j}(u) \ dx} + \int_{\Omega}{\dfrac{a_{m}'(u)}{(a_{m}(u) + 1/j)^{2-s}} \norm{\nabla u}{}^2 \ dx} = \prodl{\dfrac{a_{m}(u)}{(a_{m}(u) + 1/j)^{2-s}} \nabla (z^2)}{\nabla a_{m}(u)},
        \end{equation*}
        where
        \begin{equation*}
          g_{m,j}(r) = \int_0^r{\dfrac{(a_{m}(\theta) + 1/j)^{s-1}}{(s-1)} \ d\theta}
        \end{equation*}
        is a primitive of $(a_{m}(r) + 1/j)^{s-1}/(s-1)$. Since $0 \leq a_m'(u) \leq C$, we have $(a_m'(u))^2 \leq C a_m'(u)$, we can write
        \begin{align*}
          \int_{\Omega}{\frac{a_m'(u)}{(a_m(u) + 1/j)^{2-s}} \norm{\nabla u}{}^2 \ dx} \geq C \int_{\Omega}{\frac{(a_m'(u))^2}{(a_m(u) + 1/j)^{2-s}} \norm{\nabla u_m}{}^2 \ dx} \geq C \int_{\Omega}{\norm{\nabla [a_m(u) + 1]^{s/2}}{}^2 \ dx}.
        \end{align*}
        and hence we obtain
        \begin{align*}
        &  \dfrac{d}{dt} \int_{\Omega} g_{m,j}(u) \ dx + C \int_{\Omega}{\norm{\nabla [a_{m}(u) + 1/j]^{s/2}}{}^2 \ dx} = \prodl{\dfrac{a_{m}(u) + 1/j - 1/j}{(a_{m}(u) + 1/j)^{2-s}} \nabla (z^2)}{\nabla a_{m}(u)} \\
       &    = \prodl{(a_{m}(u) + 1/j)^{s-1} \nabla (z^2)}{\nabla a_{m}(u)} - 2 \prodl{z \nabla z}{\left ( \dfrac{1/j}{(a_{m}(u) + 1/j)} \right )^{1-s/2} \dfrac{(1/j)^{s/2} \nabla a_{m}(u)}{(a_{m}(u) + 1/j)^{1-s/2}}} \\
     &      \leq \prodl{(a_{m}(u) + 1/j)^{s-1} \nabla (z^2)}{\nabla a_{m}(u)} + \dfrac{4}{s} \sqrt{\norma{v^0}{L^{\infty}(\Omega)} + \alpha} (1/j)^{s/2} \norma{\nabla z}{L^2(\Omega)} \norma{\nabla [a_{m}(u) + 1/j]^{s/2}}{L^2(\Omega)},
        \end{align*}
        where in the last estimate we use $a_{m}(u) + 1/j \geq 1/j$. Then, using Young's inequality, we can absorb the term $\norma{\nabla [a_{m}(u) + 1/j]^{s/2}}{L^2(\Omega)}$ obtaining
        \begin{equation*}
          \begin{array}{l}
            \dfrac{d}{dt} \D{\int_{\Omega}{g_{m,j}(u) \ dx}} + C \D{\int_{\Omega}{\norm{\nabla [a_{m}(u) + 1/j]^{s/2}}{}^2 \ dx}} \\[12pt]
        \quad     \leq \D{\prodl{(a_{m}(u) + 1/j)^{s-1} \nabla (z^2)}{\nabla a_{m}(u)} + C (1/j)^s \norma{\nabla z}{L^2(\Omega)}^2}.
          \end{array}
        \end{equation*}
        
        \vspace{12pt}
        {\bfseries STEP 2:} We add the inequality of Lemma \ref{lema_tratamento_equacao_z} to $s/4$ times \eqref{estimativa_u_m_1_s=1}, then we obtain
        \begin{align*}
          \dfrac{d}{dt} \Big [ \frac{s}{4} \int_{\Omega}g_{m,j}(u) & dx  + \dfrac{1}{2} \norma{\nabla z}{L^2(\Omega)}^2 \Big ] + C \int_{\Omega}{\norm{\nabla [a_{m}(u) + 1/j]^{s/2}}{}^2  dx} \\
          & \quad + C_1 \Big ( \int_{\Omega}{\norm{D^2 z}{}^2 \ dx} + \int_{\Omega}{\frac{\norm{\nabla z}{}^4}{z^2} \ dx} \Big ) + \frac{1}{2} \int_{\Omega}{a_m(u)^s \norm{\nabla z}{}^2 \ dx} \\
          & \leq \dfrac{s}{2} \sqrt{\alpha} \int_{\Omega}{a_m(u)^{s-1} \norm{\nabla z}{} \norm{\nabla a_m(u)}{} \ dx} + C \norma{\nabla z}{L^2(\Omega)}^2 \\
          & \quad + s \int_{\Omega}{\Big [ (a_{m}(u) + 1/j)^{s-1} - a_{m}(u)^{s-1} \Big ] \nabla a_{m}(u) \cdot z \nabla z \ dx} \\
          & \leq \int_{\Omega}{\sqrt{\alpha} \norm{\nabla [a_{m}(u) + 1/j]^{s/2}}{} \norm{(a_{m}(u) + 1/j)^{s/2}}{} \norm{\nabla z}{} \ dx} + C \norma{\nabla z}{L^2(\Omega)}^2 \\
          & \quad + \frac{s}{4} \int_{\Omega}{\Big [ (a_{m}(u) + 1/j)^{s-1} - a_{m}(u)^{s-1} \Big ] \nabla a_{m}(u) \cdot \nabla (z^2) \ dx}.
        \end{align*}
        
        Next, we deal with the first term in the right hand side of the previous inequality using Hölder's and Young's inequality,
        \begin{align*}
      &    \int_{\Omega}{\sqrt{\alpha} \norm{\nabla [a_{m}(u) +1/j]^{s/2}}{} \norm{(a_{m}(u) + 1/j)^{s/2}}{} \norm{\nabla z}{} \ dx} \\
     &     \leq \sqrt{\alpha} \left[ \int_{\Omega}{(a_{m}(u) + 1/j)^s \norm{\nabla z}{}^2 \ dx} \right ]^{1/2} \norma{\nabla [a_{m}(u) + 1/j]^{s/2}}{L^2(\Omega)} \\
    &      \leq \alpha \ C(\delta) \int_{\Omega}{(a_{m}(u) + 1/j)^s \norm{\nabla z}{}^2 \ dx}  + \delta \norma{\nabla [a_{m}(u) +1/j]^{s/2}}{L^2(\Omega)}^2.
        \end{align*}
        
        Therefore, we first choose $\delta > 0$ and then $\alpha > 0$ sufficiently small in order to obtain
        \begin{align*}
        &    \dfrac{d}{dt} \left[ \dfrac{s}{4}  \D{\int_{\Omega}} g_{m,j}(u) \ dx + \frac{1}{2} \norma{\nabla z}{L^2(\Omega)}^2 \right] + C \int_{\Omega}{\norm{\nabla [a_{m}(u) + 1/j]^{s/2}}{}^2 \ dx} \\
            & \quad + C_1 \Big ( \D{\int_{\Omega}}{\norm{D^2 z}{}^2 \ dx} + \D{\int_{\Omega}}{\frac{\norm{\nabla z}{}^4}{z^2} \ dx} \Big ) + \frac{1}{2}  \D{\int_{\Omega}}{\Big [ a_{m}(u)^s - \dfrac{1}{2} (a_{m}(u) + 1/j)^s \Big ] \norm{\nabla z}{}^2 \ dx} \\
            & \leq C \norma{\nabla z}{L^2(\Omega)}^2 + \frac{s}{4} \D{\int_{\Omega}{\Big [ a_{m}(u)^{s-1} - (a_{m}(u) + 1/j)^{s-1} \Big ] \nabla a_{m}(u) \cdot \nabla (z^2) \ dx}}.
        \end{align*}
        In order to avoid problems with divisions by zero in the term $C \D{\int_{\Omega}}{\norm{\nabla [a_{m}(u) + 1/j]^{s/2}}{}^2 \ dx}$ as we take the limit as $j \to \infty$, we use the fact that this term is nonnegative and write
        \begin{equation}
          \begin{array}{l}
            \dfrac{d}{dt} \left[ \dfrac{s}{4}  \D{\int_{\Omega}} g_{m,j}(u) \ dx + \dfrac{1}{2} \norma{\nabla z}{L^2(\Omega)}^2 \right] + C_1 \Big ( \D{\int_{\Omega}}{\norm{D^2 z}{}^2 \ dx} + \D{\int_{\Omega}}{\frac{\norm{\nabla z}{}^4}{z^2} \ dx} \Big ) \\[12pt]
             \quad + \dfrac{1}{2} \D{\int_{\Omega}}{\Big [ a_{m}(u)^s - \dfrac{1}{2} (a_{m}(u) + 1/j)^s \Big ] \norm{\nabla z}{}^2 \ dx} \\[12pt]
             \leq C \norma{\nabla z}{L^2(\Omega)}^2 + \dfrac{s}{4} \D{\int_{\Omega}{\Big [ a_{m}(u)^{s-1} - (a_{m}(u) + 1/j)^{s-1} \Big ] \nabla a_{m}(u) \cdot \nabla (z^2) \ dx}}.
          \end{array}
          \label{eq_estimativa_j}
        \end{equation}
        
        \vspace{12pt}
        {\bfseries STEP 3:} We pass to the limit as $j \to \infty$.
        
        Now we show that, passing to the limit as $j \to \infty$, we recover the cancellation of the chemotaxis and consumption terms. To deal with the passage to the limit as $j \to \infty$, we remind that $m \in \mathbf{N}$ is fixed and that the solution $(u_m,z_m)$ of \eqref{problema_P_u_m_z_m}, denoted for simplicity as $(u,z)$ in the present subsection, have the regularity
        \begin{equation*}
          (u - u^\ast), z \in L^{\infty}(0,\infty;H^1(\Omega)) \cap L^2(0,\infty;H^2(\Omega)),
        \end{equation*}
        \begin{equation*}
          u_t \in L^2(0,\infty;L^2(\Omega))
          \hspace{5mm}\mbox{and} \hspace{5mm}
          z_t \in L^2(0,\infty;H^1(\Omega)).
        \end{equation*}
        
        This means that there is a zero measure set $\mathcal{N} \subset (0,\infty)$ such that 
        for any $t\in (0,\infty) \setminus \mathcal{N}$ we have
        \begin{equation*}
          u_t(t,\cdot), \nabla z_t(t,\cdot),  u(t,\cdot), z(t,\cdot), \nabla u(t,\cdot), \nabla z(t,\cdot), D^2 z(t,\cdot) \in L^2(\Omega)
        \end{equation*}
        and, by Corollary \ref{coro_expressao_derivada_integral_de_g_m_j}, we have \eqref{expressao_derivada_integral_de_g_m_j}. Therefore each integral of the inequality \eqref{eq_estimativa_j} is well defined and \eqref{eq_estimativa_j} is satisfied for each $t \in (0,\infty) \setminus \mathcal{N}$.
        
        We want to take to the limit as $j \to \infty$ in \eqref{eq_estimativa_j}. We are going to do it term by term. Let $t \in (0,\infty) \setminus \mathcal{N}$ and let us first consider the term \eqref{expressao_derivada_integral_de_g_m_j}. We define the functions $f,F,f_j \in L^1(\Omega)$, for all $j \in \mathbb{N}$, by
        \begin{equation*}
          f_j(x) = u_t(t,x) \frac{(a_m(u(t,x)) + 1/j)^{s-1}}{s-1}, \hspace{1cm}
          f(x) = u_t(t,x) \frac{a_m(u(t,x))^{s-1}}{s-1}
        \end{equation*}
       \begin{equation*}
         \mbox{and} \hspace{1cm} F(x) = \norm{f_1(x)}{} = \norm{u_t(t,x)}{} \frac{(a_m(u(t,x)) + 1)^{s-1}}{s-1}.
       \end{equation*}
       Then, for almost every $x \in \Omega$, $f_j(x) \to f(x)$ as $j \to \infty$ with $\norm{f_j(x)}{} \leq F(x)$ for all $j \in \mathbb{N}$ and, by the Dominated Convergence Theorem, we conclude that $f_j \to f$ in $L^1(\Omega)$ as $j \to \infty$. This implies, in particular, that
       \begin{equation} \label{TCD_consequencia}
         \int_{\Omega}{f_j \ dx } \longrightarrow \int_{\Omega}{f \ dx }, \mbox{ as } j \to \infty.
       \end{equation}
       Therefore, using \eqref{expressao_derivada_integral_de_g_m_j}, \eqref{TCD_consequencia} and then \eqref{expressao_derivada_integral_de_g_m} we conclude that
       \begin{align*}
         \lim_{j \to \infty}{\dfrac{d}{dt} \int_{\Omega}{g_{m,j}(u(t,x)) \ dx}} = \lim_{j \to \infty}{\int_{\Omega}{u_t(t,x) \frac{(a_m(u(t,x)) + 1/j)^{s-1}}{s-1} \ dx}} \\
         = \int_{\Omega}{u_t(t,x) \frac{a_m(u(t,x))^{s-1}}{s-1} \ dx} = \dfrac{d}{dt} \int_{\Omega}{g_m(u(t,x)) \ dx}, \mbox{ for each } t \in (0,\infty) \setminus \mathcal{N}.
       \end{align*}
       
       We can follow this reasoning and take the limit as $j \to \infty$ in the other terms of the \eqref{eq_estimativa_j}. Using the Dominated Convergence Theorem again we conclude that
       \begin{align*}
         \lim_{j \to \infty}{\dfrac{1}{2} \int_{\Omega}{\Big [ a_{m}(u(t,x))^s - \dfrac{1}{2} (a_{m}(u(t,x)) + 1/j)^s \Big ] \norm{\nabla z(t,x)}{}^2 \ dx}} = \dfrac{1}{4} \int_{\Omega}{a_{m}(u(t,x))^s \norm{\nabla z(t,x)}{}^2 \ dx}
        \end{align*}
        and
        \begin{equation*}
          \lim_{j \to \infty}{s \int_{\Omega}{\Big [ a_{m}(u(t,x))^{s-1} - (a_{m}(u(t,x)) + 1/j)^{s-1} \Big ] \nabla a_{m}(u(t,x)) \cdot z(t,x) \nabla z(t,x) \ dx}} = 0,
        \end{equation*}
        for each $t \in (0,\infty) \setminus \mathcal{N}$.
       
       Then, since the limit preserves inequalities, after we take the limit as $j \to \infty$ in \eqref{eq_estimativa_j}, we obtain
       \begin{align*}
         \dfrac{d}{dt} \Big [ \frac{s}{4} \int_{\Omega}{g_m(u(t,x)) \ dx} + \frac{1}{2} \norma{\nabla z(t,x)}{L^2(\Omega)}^2 \Big ] + \frac{1}{4} \int_{\Omega}{a_{m}(u(t,x))^s \norm{\nabla v(t,x)}{}^2 \ dx} \\
         + C_1 \Big ( \D{\int_{\Omega}}{\norm{D^2 z(t,x)}{}^2 \ dx} + \D{\int_{\Omega}}{\frac{\norm{\nabla z(t,x)}{}^4}{z(t,x)^2} \ dx} \Big ) \leq C \norma{\nabla z}{L^2(\Omega)}^2
       \end{align*}
       for all $t \in (0,\infty) \setminus \mathcal{N}$, which means that the inequality is valid for almost every $t \in (0,\infty)$. Therefore \eqref{estimativa_u_v_m_1_s_intermediario} holds.
      \end{proof}
      
      \begin{lema}[\bfseries Energy inequality for $\boldsymbol{s \geq 2}$]
        The solution $(u, z)$ of the problem \eqref{problema_P_u_m_z_m} satisfies, for sufficiently small $\alpha > 0$,
        \begin{equation}
          \begin{array}{c}
            \dfrac{d}{dt} \Big [ \dfrac{s}{4} \D{\int_{\Omega}} g_m(u) \ dx  + \frac{1}{2} \norma{\nabla z}{L^2(\Omega)}^2 \Big ] + \D{\int_{\Omega}}{\norm{\nabla [a_{m}(u)]^{s/2}}{}^2  dx} + \dfrac{1}{4} \D{\int_{\Omega}}{a_m(u)^s \norm{\nabla z}{}^2 \ dx} \\[12pt]
          \quad + C_1 \Big ( \D{\int_{\Omega}}{\norm{D^2 z}{}^2 \ dx} + \D{\int_{\Omega}}{\frac{\norm{\nabla z}{}^4}{z^2} \ dx} \Big ) \leq C \int_{\Omega}{\norm{\nabla z}{}^2 \ dx},
          \end{array}
          \label{estimativa_u_v_m_1_s_geq_2}
        \end{equation}
        where $g_m(u)$ is given by \eqref{definicao_g_m_e_g_m_j}.
        \label{lemma_estimativa_u_v_m_1_s_geq_2}
      \end{lema}
      \begin{proof}[\bfseries Proof]
        We test the $u$-equation of \eqref{problema_P_u_m_z_m} by
        \begin{equation*}
          g_m'(u) = \dfrac{(a_m(u))^{s-1}}{(s-1)}
        \end{equation*}
        and obtain
        \begin{equation*}
          \dfrac{d}{dt} \int_{\Omega}{g_m(u) \ dx} + \int_{\Omega}{(a_{m}(u))^{s-2} a_{m}'(u) \norm{\nabla u}{}^2 \ dx} = \prodl{a_{m}(u) (a_{m}(u))^{s-2} \nabla (z^2)}{\nabla a_{m}(u)},
        \end{equation*}
        where
        \begin{equation*}
          g_m(r) = \int_0^r{\dfrac{(a_{m}(\theta))^{s-1}}{(s-1)} \ d\theta}
        \end{equation*}
        is a primitive of $(a_{m}(r))^{s-1}/(s-1)$. Since $0 \leq a_m'(u) \leq C$, we have $(a_m'(u))^2 \leq C a_m'(u)$, we can write
        \begin{align*}
          \int_{\Omega}{a_m'(u) (a_m(u))^{s-2} \norm{\nabla u}{}^2 \ dx} \geq C \int_{\Omega}{(a_m'(u))^2 (a_m(u))^{s-2} \norm{\nabla u}{}^2 \ dx} \geq C \int_{\Omega}{\norm{\nabla (a_m(u))^{s/2}}{}^2 \ dx}.
        \end{align*}
        Then we obtain
        \begin{equation}
          \dfrac{d}{dt} \D{\int_{\Omega}{g_m(u) \ dx}} + C \D{\int_{\Omega}{\norm{\nabla [a_{m}(u)]^{s/2}}{}^2 \ dx}} \leq \D{\prodl{a_{m}(u)^{s-1} \nabla (z^2)}{\nabla a_{m}(u)}}.
          \label{estimativa_u_m_1_s_geq_2}
        \end{equation}
        
        If we add $s/4$ times \eqref{estimativa_u_m_1_s_geq_2} to the inequality of Lemma \ref{lema_tratamento_equacao_z} then the term
        \begin{equation*}
          \frac{s}{4} \D{\int_{\Omega}}{a_{m}(u)^{s-1} \nabla a_{m}(u) \cdot \nabla (z^2) \ dx},
        \end{equation*}
        which appears in $s/4$ times \eqref{estimativa_u_m_1_s_geq_2} cancels with the term
        \begin{equation*}
          - \frac{s}{4} \D{\int_{\Omega}}{a_{m}(u)^{s-1} \nabla a_{m}(u) \cdot \nabla (z^2) \ dx},
        \end{equation*}
        which comes from the inequality Lemma \ref{lema_tratamento_equacao_z} and we obtain
        \begin{align*}
          \dfrac{d}{dt} \Big [ \frac{s}{4} \int_{\Omega}g_m(u) & dx  + \frac{1}{2} \norma{\nabla z}{L^2(\Omega)}^2 \Big ] + C \int_{\Omega}{\norm{\nabla [a_{m}(u)]^{s/2}}{}^2  dx} + \frac{1}{2} \int_{\Omega}{a_m(u)^s \norm{\nabla z}{}^2 \ dx} \\
          & \quad + C_1 \Big ( \int_{\Omega}{\norm{D^2 z}{}^2 \ dx} + \frac{1}{2} \int_{\Omega}{\frac{\norm{\nabla z}{}^4}{z^2} \ dx} \Big ) \leq \frac{s}{2} \sqrt{\alpha} \int_{\Omega}{a_m(u)^{s-1} \norm{\nabla z}{} \norm{\nabla a_m(u)}{} \ dx} \\
          & \leq \int_{\Omega}{\sqrt{\alpha} \norm{\nabla [a_{m}(u)]^{s/2}}{} \norm{a_{m}(u)^{s/2}}{} \norm{\nabla z}{} \ dx}.
        \end{align*}
        
        Next, we deal with the second term in the right hand side of the previous inequality using Hölder's and Young's inequality,
        \begin{align*}
          & \int_{\Omega}{\sqrt{\alpha} \norm{\nabla [a_{m}(u)]^{s/2}}{} \norm{a_{m}(u)^{s/2}}{} \norm{\nabla z}{} \ dx} \\
          & \leq \alpha \ C(\delta) \int_{\Omega}{a_{m}(u)^s \norm{\nabla z}{}^2 \ dx}  + \delta \norma{\nabla [a_{m}(u)]^{s/2}}{L^2(\Omega)}^2.
        \end{align*}
        
        Therefore, choosing $\alpha, \delta > 0$ sufficiently small we finally obtain the desired inequality \eqref{estimativa_u_v_m_1_s_geq_2}.
      \end{proof}
      
      The energy inequalities \eqref{estimativa_u_v_m_1_s=1}, \eqref{estimativa_u_v_m_1_s_intermediario} and \eqref{estimativa_u_v_m_1_s_geq_2} allow us to obtain $m$-independent estimates for the function $v$ that are valid up to infinity time in the next Subsection.
    
      \begin{obsv}
        In the next subsection, the aforementioned $m$-independent estimates will obtained upon integration of the energy inequalities \eqref{estimativa_u_v_m_1_s=1}, \eqref{estimativa_u_v_m_1_s_intermediario} and \eqref{estimativa_u_v_m_1_s_geq_2} with respect to the time variable. Therefore we find it appropriate to remark that, for each $T > 0$, we have $\D{\int_{\Omega}g_m(u) \ dx} \in L^2(0,T)$ and $\D{\frac{d}{dt} \int_{\Omega}g_m(u) \ dx} \in L^2(0,T)$ and it implies, in particular, that
        \begin{equation*}
          \int_0^T{\frac{d}{dt} \int_{\Omega}g_m(u(t,x)) \ dx \ dt} = \int_{\Omega}g_m(u(T,x)) \ dx - \int_{\Omega}g_m(u(0,x)) \ dx.
        \end{equation*}
        See \cite{Brezis_functional}.
        \hfill \qedsymbol
      \end{obsv}


    \subsection{\texorpdfstring{$\boldsymbol{m}$}{m}-independent estimates and passage to the limit as \texorpdfstring{$\boldsymbol{m \to \infty}$}{m goes to infinity}}
    \label{subsec:passage_limit_m_to_infty}
    
    Now we use again the notation $(u_m,v_m)$ for the solution of the regularized problem \eqref{problema_P_m_intro}, $z_m = \sqrt{v_m + \alpha}$ and $(u,v)$ for the solution of the original problem \eqref{problema_P}. In this subsection we are going to obtain $m$-independent estimates for $(u_m,v_m)$ that will allow us to pass to the limit in the problem \eqref{problema_P_m_intro} as $m \to \infty$ and prove the existence of solution to the original problem \eqref{problema_P}.
    
    First, in Subsection \ref{subsec:estimativas_v_m}, we will obtain some $m$-independent estimates for $\nabla v_m$ that can be extracted from the energy inequalities \eqref{estimativa_u_v_m_1_s=1}, \eqref{estimativa_u_v_m_1_s_intermediario} and \eqref{estimativa_u_v_m_1_s_geq_2}. Next, in Subsections \ref{subsec:estimativas_u_m_s<2} and \ref{subsec:estimativas_u_m_s_geq_2} we obtain $m$-independent bounds for $(u_m,v_m)$ and pass to the limit in \eqref{problema_P_m_intro} as $m \to \infty$, considering the case $s \in [1,2)$ and $s \geq 2$, respectively.
    
    
    \subsubsection{\texorpdfstring{$\boldsymbol{m}$}{m}-independent estimates for \texorpdfstring{$\boldsymbol{\nabla v_m}$}{grad(v\_m)}}
    \label{subsec:estimativas_v_m}
      
      Let us remind that, for $s \geq 1$, we have defined $g_m'$ as
      \begin{equation*}
        g_m'(u) = \left \{
        \begin{array}{rl}
          ln(a_m(u) + 1) & \mbox{ if } s = 1,  \\
          a_m(u)^{s-1}/(s-1) & \mbox{ if } s > 1.
        \end{array}
        \right .
      \end{equation*}
      And let us define the energy
      \begin{equation}
        E_m(t) = \frac{s}{4} \int_{\Omega}{g_m(u_m(t,x)) \ dx} + \frac{1}{2} \int_{\Omega}{\norm{\nabla z_m(t,x)}{}^2} \ dx.
        \label{energia_E_m}
      \end{equation}
      
      We remark that, since $0 \leq v_m(t,x) \leq \norma{v^0}{L^{\infty}(\Omega)} \ a.e. \ (t,x) \in (0,\infty) \times \Omega$, we have
      \begin{equation*}
        0 < \sqrt{\alpha} \leq z_m(t,x) \leq \sqrt{\norma{v^0}{L^{\infty}(\Omega)} + \alpha} \ a.e. \ (t,x) \in (0,\infty) \times \Omega
      \end{equation*}
      and, by straightforward calculations, we can prove the following lemma.
      
      \begin{lema}
        There are $\beta_1, \beta_2 > 0$, depending on $\alpha$, such that
        \begin{equation}
          \beta_1 \norm{\nabla z_m(t,x)}{} \leq \norm{\nabla v_m(t,x)}{} \leq \beta_2 \norm{\nabla z_m(t,x)}{}
          \label{estimativa_independente_m_nabla_z}
        \end{equation}
        and
        \begin{equation}
          \beta_1 \Big ( \norm{\Delta z_m(t,x)}{} + \norm{\nabla z_m(t,x)}{}^2 \Big ) \leq \norm{\Delta v_m(t,x)}{} \leq \beta_2 \Big ( \norm{\Delta z_m(t,x)}{} + \norm{\nabla z_m(t,x)}{}^2 \Big ),
          \label{estimativa_independente_m_Delta_z}
        \end{equation}
        $a.e. \ (t,x) \in (0,\infty) \times \Omega.$
      \end{lema}
      
      We will integrate the energy inequalities \eqref{estimativa_u_v_m_1_s=1}, \eqref{estimativa_u_v_m_1_s_intermediario} and \eqref{estimativa_u_v_m_1_s_geq_2} with respect to $t$, from $0$ to some $T > 0$. We take into account that $\D{\int_0^T{\norma{\nabla z_m(t)}{L^2(\Omega)}^2 \ dt}}$ is bounded, independently of $T$ and $m$, because of \eqref{estimativa_independente_m_nabla_z} and \eqref{limitacao_uniforme_nabla_v_m}, and we use the $m$-uniform bounds which stem from \eqref{convergencia_dado_inicial_u_m} and \eqref{convergencia_dado_inicial_v_m} on the initial data $u^0_m$ and $v^0_m$ in order to conclude that the energy given in \eqref{energia_E_m} in time $t = 0$, $E_m(0)$, is also bounded, independently of $m$. Thus we can conclude that
      \begin{equation*}
        \nabla z_m \mbox{ is bounded in } L^{\infty}(0,\infty;L^2(\Omega)) \cap L^4(0,\infty;L^4(\Omega)),
      \end{equation*}
      \begin{equation*}
        a_m(u_m)^{s/2} \nabla z_m \mbox{ and } \Delta z_m \mbox{ are bounded in } L^2(0,\infty;L^2(\Omega)).
      \end{equation*}
      But using the fact that $z_m = \sqrt{v_m + \alpha}$, \eqref{estimativa_independente_m_nabla_z} and \eqref{estimativa_independente_m_Delta_z} we can conclude that
      \begin{equation}
        \nabla v_m \mbox{ is bounded in } L^{\infty}(0,\infty;L^2(\Omega)) \cap L^4(0,\infty;L^4(\Omega)),
      \label{limitacao_aux_Dv_m_s_intermediario}
      \end{equation}
      \begin{equation}
        a_m(u_m)^{s/2} \nabla v_m \mbox{ and } \Delta v_m \mbox{ are bounded in } L^2(0,\infty;L^2(\Omega)).
        \label{limitacao_aux_Delta_v_m_s_intermediario}
      \end{equation}
      
      In particular, since $v_m(t) \in H^2(\Omega)$, for each $t \in (0,\infty)$, and $\dfrac{\partial}{\partial \eta} v_m \Big |_{\Gamma} = 0$, it stems from \eqref{limitacao_aux_Delta_v_m_s_intermediario}, the $H^2$-regularity of the Poisson-Neumann problem \eqref{Neumann_problem} and \eqref{obsv_norma_de_nabla_v} that
      \begin{equation}
        \nabla v_m \mbox{ is bounded in } L^2(0,\infty;H^1(\Omega)).
        \label{limitacao_aux_Dv_m_s_intermediario_2}
      \end{equation}
      
      Using the results obtained until this point we analyze the existence of solutions of \eqref{problema_P}, first for $s \in [1,2)$ and then for $s \geq 2$.
      
      
      \subsubsection{\texorpdfstring{$\boldsymbol{m}$}{m}-independent estimates for \texorpdfstring{$\boldsymbol{(u_m,v_m)}$}{(u\_m,v\_m)} and passage to the limit for \texorpdfstring{$\boldsymbol{s \in [1,2)}$}{m}}
      
      \label{subsec:estimativas_u_m_s<2}
      
      Let
      \begin{equation*}
        \forall r > 0, \qquad g'(r) =  \left \{
        \begin{array}{cc}
          ln(r) & \mbox{ if } s = 1, \\
          r^{s-1}/(s-1) & \mbox{ if } s \in (1,2) 
        \end{array}
        \right .
      \end{equation*}
      and let
      \begin{equation*}
        g(r) =  \int_0^r{g'(\theta) \ d\theta} = \left \{
        \begin{array}{cc}
          r ln(r) - r & \mbox{ if } s = 1, \\
          r^s/s(s-1) & \mbox{ if } s \in (1,2).
        \end{array}
        \right .
      \end{equation*}
      Notice that $g''(r) = r^{s-2}, \ \forall r > 0$, in all cases.
      
      We test the $u_m$-equation of \eqref{problema_P_m_intro} by $g'(u_m + 1)$ and obtain
      \begin{align*}
        \dfrac{d}{dt} \int_{\Omega} g(u_m + 1) \ dx & + \dfrac{4}{s^2} \int_{\Omega}{\norm{\nabla [u_m + 1]^{s/2}}{}^2 \ dx} \\
        & = \int_{\Omega}{a_m(u_m) (u_m + 1)^{s/2-1} \nabla v_m  \cdot \nabla u_m\,  (u_m + 1)^{s/2-1} \ dx} \\
        & = \dfrac{2}{s} \int_{\Omega}{\dfrac{a_m(u_m)^{1-s/2}}{(u_m + 1)^{1-s/2}} a_m(u_m)^{s/2} \nabla v_m \cdot \nabla [u_m + 1]^{s/2} \ dx} \\
        & \leq \dfrac{2}{s} \Big ( \int_{\Omega}{a_m(u_m)^s \norm{\nabla v_m}{}^2 \ dx} \Big )^{1/2} \Big ( \int_{\Omega}{\norm{\nabla [u_m + 1]^{s/2}}{}^2 \ dx} \Big )^{1/2}
      \end{align*}
      and thus we have
 $$
        \dfrac{d}{dt} \int_{\Omega}{g(u_m + 1) \ dx}  + \dfrac{2}{s^2} \int_{\Omega}{\norm{\nabla [u_m + 1]^{s/2}}{}^2 \ dx}
        \leq \dfrac{1}{4} \int_{\Omega}{a_m(u_m)^s \norm{\nabla v_m}{}^2 \ dx}.
 $$
      Integrating with respect to $t$ from $0$ to $T$, for any fixed $T \in (0,\infty)$, we obtain
      \begin{align*}
        \int_{\Omega}{g(u_m(T) + 1) \ dx} + \dfrac{2}{s^2} \int_{0}^T{\int_{\Omega}{\norm{\nabla [u_m + 1]^{s/2}}{}^2 \ dx} \ dt} \\
        \leq \dfrac{1}{4} \int_0^T{\int_{\Omega}{a_m(u_m)^s \norm{\nabla v_m}{}^2 \ dx} \ dt} + \int_{\Omega}{g(u^0_m + 1) \ dx}.
      \end{align*}
      
      Then, because of Lemma \ref{lemma_limitacao_u_v_independente_m}.1, \eqref{convergencia_dado_inicial_u_m} and the definition of $g$ and \eqref{limitacao_aux_Delta_v_m_s_intermediario} we conclude that
      \begin{equation}
          (u_m + 1)^{s/2} \mbox{ is bounded in } L^{\infty}(0,\infty;L^2(\Omega)), 
          \label{limitacao_u_m_L^s}
      \end{equation}
      in particular,
      \begin{equation}
        u_m \mbox{ is bounded in } L^{\infty}(0,\infty;L^s(\Omega))
        \label{limitacao_u_m_L^s_2},
      \end{equation}
      and
      \begin{equation}
        \nabla [u_m + 1]^{s/2} \mbox{ is bounded in } L^2(0,\infty;L^2(\Omega)).
        \label{limitacao_nabla_u_m_elevado_a_s_sobre_2}
      \end{equation}
      
      Consider the relation
      \begin{equation}
        \nabla u_m = \nabla (u_m + 1) = \nabla \big ( (u_m + 1)^{s/2} \big )^{2/s} = \dfrac{2}{s} (u_m + 1)^{1 - s/2} \ \nabla (u_m + 1)^{s/2}.
        \label{gradiente_u_m_em_funcao_de_gradiente_u_m^s/2}
      \end{equation}
      Taking into account that we are considering $s \in [1,2)$, we can use \eqref{limitacao_u_m_L^s} to obtain 
      \begin{equation*}
        (u_m + 1)^{1-s/2} \mbox{ is bounded in } L^{\infty}(0,\infty;L^{2s/(2-s)}(\Omega))
      \end{equation*}
      and then \eqref{limitacao_nabla_u_m_elevado_a_s_sobre_2} and \eqref{gradiente_u_m_em_funcao_de_gradiente_u_m^s/2} to conclude that 
      \begin{equation}
        \nabla u_m \mbox{ is bounded in } L^2(0,\infty;L^s(\Omega)).
        \label{limitacao_nabla_u_m_L^s}
      \end{equation}
      In conclusion, using \eqref{limitacao_u_m_L^s_2}, \eqref{limitacao_nabla_u_m_L^s} and the Poincaré inequality for zero mean functions (Lemma \ref{lema_desigualdade_poincare_media_nula}),
      \begin{equation}
        u_m - u^{\ast} \mbox{ is bounded in } L^{\infty}(0,\infty;L^s(\Omega)) \cap L^2(0,\infty;W^{1,s}(\Omega)).
        \label{limitacao_W_1_s}
      \end{equation}
      
      Considering the chemotaxis term of the $u_m$-equation of \eqref{problema_P_m_intro}, we can write $a_m(u_m) \nabla v_m$ as
  \begin{equation*}
    a_m(u_m) \nabla v_m = a_m(u_m)^{1-s/2} a_m(u_m)^{s/2} \nabla v_m.
  \end{equation*}
  Then, we have $a_m(u_m)^{1-s/2}$ bounded in $L^{\infty}(0,\infty;L^{2s/(2-s)}(\Omega))$, because of \eqref{limitacao_u_m_L^s}, and $a_m(u_m)^{s/2} \nabla v_m$ bounded in $L^2(0,\infty;L^2(\Omega))$, because of \eqref{limitacao_aux_Delta_v_m_s_intermediario}, and hence we can conclude that
  \begin{equation}
    a_m(u_m) \nabla v_m \mbox{ is bounded in } L^2(0,\infty;L^s(\Omega)).
    \label{limitacao_termo_chemotaxis_s}
  \end{equation}
  Then, if we consider the $u_m$-equation of \eqref{problema_P_m_intro}, from \eqref{limitacao_W_1_s} and \eqref{limitacao_termo_chemotaxis_s} we conclude that 
  \begin{equation*}
    \partial_t u_m \mbox{ is bounded in } L^2 \Big (0,\infty;\big ( W^{1,s/(s-1)}(\Omega) \big )' \Big ).
  \end{equation*}
      
      Now we turn to the $v_m$-equation, rewritten as
      \begin{equation}
        \partial_t v_m - \Delta v_m + a_m(u^{\ast})^s v_m = - (a_m(u_m)^s - a_m(u^{\ast})^s) v_m.
        \label{eq_v_m_equivalente_s<2}
      \end{equation}
      Analyzing the term on the right hand side of \eqref{eq_v_m_equivalente_s<2}, we have
      \begin{equation}
        a_m(u_m)^s - a_m(u^{\ast})^s \mbox{ is bounded in } L^2(0,\infty;L^{3/2}(\Omega)).
        \label{limitacao_u_s_menos_u_estrela_s}
      \end{equation}
      In fact, using Lemma \ref{lema_a_m_elevado_a_s}, we obtain
      \begin{align*}
        \norm{a_m(u_m)^s - a_m(u^{\ast})^s}{} \leq s \norm{a_m(u_m) + a_m(u^{\ast})}{}^{s-1} \norm{u_m - u^{\ast}}{}.
      \end{align*}
      Then, considering the $m$-uniform bound \eqref{limitacao_W_1_s} and the Sobolev embedding $L^{3s/(3-s)}(\Omega)) \subset W^{1,s}(\Omega)$ we obtain \eqref{limitacao_u_s_menos_u_estrela_s}.
      
      With this information, now we can test \eqref{eq_v_m_equivalente_s<2} by $v_m$, obtaining
      \begin{align*}
        \frac{1}{2} \frac{d}{dt} \norma{v_m}{L^2(\Omega)}^2 & + \norma{\nabla v_m}{L^2(\Omega)}^2 + \frac{a_m(u^{\ast})^s}{2} \norma{v_m}{L^2(\Omega)}^2 \leq C \int_{\Omega}{\norm{a_m(u_m)^s - a_m(u^{\ast})^s}{} v_m^2 \ dx} \\
        & \leq C \norma{v^0_m}{L^{\infty}(\Omega)} \norma{a_m(u_m)^s - a_m(u^{\ast})^s}{L^{3/2}(\Omega)} \norma{v_m}{L^3(\Omega)} \\
        & \leq C(\delta) \norma{a_m(u_m)^s - a_m(u^{\ast})^s}{L^{3/2}(\Omega)}^2 + \delta \norma{v_m}{L^2(\Omega)}^2 + \delta \norma{\nabla v_m}{L^2(\Omega)}^2.
      \end{align*}
      Note that if $u_0 \not\equiv 0$ then have $a_m(u^{\ast}) = u^{\ast} > 0$ all $m \geq u^{\ast}$. Hence, choosing $\delta > 0$ small enough, we can conclude that, for $m \geq u^{\ast}$, there is $\beta > 0$ such that
      \begin{align*}
        \frac{1}{2} \frac{d}{dt} \norma{v_m}{L^2(\Omega)}^2 + \beta \norma{\nabla v_m}{L^2(\Omega)}^2 + \beta \norma{v_m}{L^2(\Omega)}^2 \leq C \norma{a_m(u_m)^s - a_m(u^{\ast})^s}{L^{3/2}(\Omega)}^2.
      \end{align*}
      Therefore, integrating with respect to $t$ and using \eqref{limitacao_u_s_menos_u_estrela_s} we obtain
      \begin{equation}
        v_m \mbox{ is bounded in } L^2(0,\infty;L^2(\Omega)).
        \label{estimativa_v_m_L^2}
      \end{equation}
      Hence, in view of \eqref{limitacao_aux_Dv_m_s_intermediario}, \eqref{limitacao_aux_Delta_v_m_s_intermediario} and \eqref{estimativa_v_m_L^2} we have
      \begin{equation*}
        v_m \mbox{ is bounded in } L^2(0,\infty;H^2(\Omega)).
      \end{equation*}
      
      With the $m$-uniform bounds obtained so far we can obtain a $m$-uniform bound for the function $\partial_t v_m$ in $L^2(0,\infty;L^{3/2}(\Omega))$. In fact, going back to \eqref{eq_v_m_equivalente_s<2}, reminding that $v_m$ is uniformly bounded in $L^{\infty}(0,\infty;L^{\infty}(\Omega))$ with respect to $m$ and considering \eqref{estimativa_v_m_L^2} and \eqref{limitacao_u_s_menos_u_estrela_s}, we conclude that
      \begin{equation*}
        \partial_t v_m \mbox{ is bounded in } L^2(0,\infty;L^{3/2}(\Omega)).
      \end{equation*}
      
      Now we are going to obtain compactness for $\{ u_m \}$ which are necessary in order to pass to the limit as $m \to \infty$ in the nonlinear terms of the equations of \eqref{problema_P_m_intro}. Because of \eqref{limitacao_u_m_L^s} and \eqref{limitacao_nabla_u_m_elevado_a_s_sobre_2}, we have that
      \begin{equation*}
        (u_m + 1)^{s/2} \mbox{ is bounded in } L^{\infty}(0,\infty;L^2(\Omega)) \cap L^2(0,T;H^1(\Omega)), \mbox{ for every finite } T > 0.
      \end{equation*}
      Using the Sobolev inequality $H^1(\Omega) \subset L^6(\Omega)$ and interpolation inequalities we obtain
      \begin{equation*}
        (u_m)^{s/2} \mbox{ is bounded in } L^{10/3}(0,T;L^{10/3}(\Omega)),
      \end{equation*}
      which is equivalent to
      \begin{equation}
        u_m \mbox{ is bounded in } L^{5s/3}(0,T;L^{5s/3}(\Omega)).
        \label{limitacao_u_m_L_5s/3}
      \end{equation}
      By using \eqref{limitacao_u_m_L^s} and \eqref{limitacao_u_m_L_5s/3} in  \eqref{gradiente_u_m_em_funcao_de_gradiente_u_m^s/2} (remind that $s \in [1,2)$), we also have
      \begin{equation*}
        u_m \mbox{ is bounded in } L^{5s/(3+s)}(0,T;W^{1,5s/(3+s)}(\Omega)).
      \end{equation*}
  
  We observe that $W^{1,5s/(3+s)}(\Omega) \subset L^q(\Omega)$, with continuous embedding for $q = 15s/(9-2s)$ and compact embedding for $q \in [1,15s/(9-2s))$. Then, since $s \in [1,2)$, we have $5s/3 < 15s/(9-2s)$ and therefore the embedding $W^{1,5s/(3+s)}(\Omega) \subset L^{5s/3}(\Omega)$ is compact. Note also that $q = 5s/3 \geq 5/3 > 1$.
  
  Now we can use Lemma \ref{lema_Simon} with 
 \[ X = W^{1,5s/(3+s)}(\Omega),\hspace{1cm}
      B = L^{5s/3}(\Omega),  \hspace{1cm}  Y = \big ( H^3(\Omega) \big )'\]
  and $q = 5s/3$, to conclude that there is a subsequence of $\{ u_m \}$ (still denoted by $\{ u_m \}$) and a limit function $u$ such that
  \begin{equation*}
    u_m \longrightarrow u \mbox{ weakly in } L^{5s/(3+s)}(0,T;W^{1,5s/(3+s)}(\Omega)), \ \forall T > 0,
  \end{equation*}
  and
  \begin{equation}
    u_m \longrightarrow u \mbox{ strongly in } L^p(0,T;L^p(\Omega)), \ \forall p \in [1,5s/3), \ \forall T > 0.   \label{convergencia_u_m_s_intermediario}
  \end{equation}
  
  Using the Dominated Convergence Theorem we can conclude from \eqref{convergencia_u_m_s_intermediario} that
  \begin{equation}
    a_m(u_m) \rightarrow u \mbox{ strongly in } L^p(0,T;L^p(\Omega)), \ \forall p \in (1,5s/3), \ \forall T > 0.
    \label{convergencia_a_m_u_s_intermediario}
  \end{equation}
  It stems from the convergence \eqref{convergencia_a_m_u_s_intermediario} and Lemma \ref{lema_convergencia_w_elevado_a_s} that
  \begin{equation}
    (a_m(u_m))^s \rightarrow u^s \mbox{ strongly in } L^q(0,T;L^q(\Omega)), \ \forall q \in (1,5/3), \ \forall T > 0.
    \label{convergencia_a_m_u_elevado_a_s_s_intermediario}
  \end{equation}
  
  The convergence of $v_m$ is better. There is a subsequence of $\{ v_m \}$ (still denoted by $\{ v_m \}$) and a limit function $v$ such that
  \begin{equation}
    \begin{array}{c}
      v_m \rightarrow v \mbox{ weakly* in } L^{\infty}((0,\infty) \times\Omega) \cap L^{\infty}(0,\infty;H^1(\Omega)), \\
      v_m \rightarrow v \mbox{ weakly in } L^2(0,\infty;H^2(\Omega)), \\
      \nabla v_m \rightarrow \nabla v \mbox{ weakly in } L^4(0,\infty;L^4(\Omega)), \\
      \mbox{and } \partial_t v_m \rightarrow \partial_t v \mbox{ weakly in } L^2(0,\infty;L^{3/2}(\Omega)).
    \end{array}
    \label{convergencia_v_m_s_intermediario}
  \end{equation}
  
  Now we are going to use the weak and strong convergences obtained so far to pass to the limit as $m \to \infty$ in the equations of problem \eqref{problema_P_m_intro}. We are going to identify the limits of the nonlinear terms related to chemotaxis and consumption,
  \begin{equation*}
    a_m(u_m) \nabla v_m \mbox{ and } a_m(u_m)^s v_m,
  \end{equation*}
  respectively, with 
  \begin{equation*}
    u \nabla v \mbox{ and } u^s v.
  \end{equation*}
  In fact, considering the chemotaxis term, because of \eqref{convergencia_a_m_u_s_intermediario}, \eqref{limitacao_aux_Dv_m_s_intermediario} and \eqref{convergencia_v_m_s_intermediario}, we can conclude that
  \begin{equation*}
    a_m(u_m) \nabla v_m \longrightarrow u \nabla v \mbox{ weakly in } L^{20s/(5s+12)}(0,T;L^{20s/(5s+12)}(\Omega)), \ \forall T > 0.
  \end{equation*}
  
  Considering now the consumption term, considering \eqref{convergencia_a_m_u_elevado_a_s_s_intermediario} and \eqref{convergencia_v_m_s_intermediario} we conclude that
  \begin{equation*}
    a_m(u_m)^s v_m \longrightarrow u^s v \mbox{ weakly in } L^{5/3}(0,T;L^{5/3}(\Omega)), \ \forall T > 0.
  \end{equation*}
  
  With these identifications and all previous convergences, it is possible to pass to the limit as $m \to \infty$ in each term of the equations of \eqref{problema_P_m_intro}. In order to finish the proof of Theorem \ref{teo_3D_existencia} we must obtain the regularity (up to infinite time) which is claimed for $u$.
  
  From \eqref{limitacao_u_m_L^s} and \eqref{limitacao_nabla_u_m_elevado_a_s_sobre_2} there exists a subsequence of $\{ (u_m + 1)^{s/2} \}$, still denoted by $\{ (u_m + 1)^{s/2} \}$, and a limit function $w$ such that
  \begin{equation*}
    \begin{array}{c}
      (u_m + 1)^{s/2} \longrightarrow \varphi \mbox{ \ weakly* in } L^{\infty}(0,\infty;L^2(\Omega)) \\
      \nabla (u_m + 1)^{s/2} \longrightarrow \nabla \varphi \mbox{ \ weakly in } L^2(0,\infty;L^2(\Omega)).
    \end{array}
  \end{equation*}
  Then, using the strong convergence \eqref{convergencia_u_m_s_intermediario}, the continuity of the function $u_m \mapsto f(u_m) = (u_m + 1)^{s/2}$ and the Dominated Convergence Theorem, we prove that $\varphi = (u + 1)^{s/2}$.
  
  Analogously, because of \eqref{limitacao_aux_Delta_v_m_s_intermediario} we can conclude that, up to a subsequence, there is a limit function $\phi$ such that
  \begin{equation*}
    a_m(u_m)^{s/2} \nabla v_m \longrightarrow \phi \mbox{ weakly in } L^2(0,\infty;L^2(\Omega)).
  \end{equation*}
  And using the convergences \eqref{convergencia_a_m_u_elevado_a_s_s_intermediario} and \eqref{convergencia_v_m_s_intermediario} we can conclude that $\phi = u^{s/2} \nabla v$.
  
  Therefore we have proved the global in time regularity
  \begin{equation}
    \begin{array}{c}
      (u + 1)^{s/2} \in L^{\infty}(0,\infty;L^2(\Omega)), \quad
      \nabla (u + 1)^{s/2} \in L^2(0,\infty;L^2(\Omega)), \\
      u^{s/2} \nabla v \in L^2(0,\infty;L^2(\Omega)).
    \end{array}
    \label{regularidade_termos_limite}
  \end{equation}
  Considering \eqref{regularidade_termos_limite} and proceeding as in the obtaining of \eqref{limitacao_W_1_s} and \eqref{limitacao_termo_chemotaxis_s} we conclude the global in time regularity
  \begin{equation*}
    \begin{array}{c}
      u \in L^{\infty}(0,\infty;L^s(\Omega)), \quad
      \nabla u, \ u \nabla v \in L^2(0,\infty;L^s(\Omega)), \\
    \end{array}
   \end{equation*}
   finishing the proof of Theorem \ref{teo_3D_existencia} in the case $s \in [1,2)$.
  
      
      \subsubsection{\texorpdfstring{$\boldsymbol{m}$}{m}-independent estimates for \texorpdfstring{$\boldsymbol{(u_m,v_m)}$}{um} and passage to the limit for \texorpdfstring{$\boldsymbol{s \geq 2}$}{m}}
      
      \label{subsec:estimativas_u_m_s_geq_2}
      
        The procedure for the case $s \geq 2$ is slightly different. First we note that, integrating the energy inequality \eqref{estimativa_u_v_m_1_s_geq_2} from Lemma \ref{lemma_estimativa_u_v_m_1_s_geq_2} with respect to $t$, we have
        \begin{equation}
          \nabla a_m(u_m)^{s/2} \mbox{ is bounded in } L^2(0,\infty;L^2(\Omega)).
          \label{limitacao_nabla_a_m_u_m_s_geq_2}
        \end{equation}
        We also remind that we defined $g_m'(r) = a_m(r)^{s-1}/(s-1)$, for $s \geq 2$. Then we have
        \begin{align*}
          a_m(r)^s = s \int_0^r{a_m'(\theta) a_m(\theta)^{s-1} \ d\theta} \leq C s \int_0^r{a_m(\theta)^{s-1} \ d\theta} = C s(s - 1) g_m(r).
        \end{align*}
        Therefore it also stems from integrating the energy inequality \eqref{estimativa_u_v_m_1_s_geq_2} with respect to $t$ that
        \begin{equation}
          a_m(u_m)^{s/2} \mbox{ is bounded in } L^{\infty}(0,\infty;L^2(\Omega)).
          \label{limitacao_a_m_u_m_s_geq_2}
        \end{equation}
        From \eqref{limitacao_a_m_u_m_s_geq_2} and \eqref{limitacao_nabla_a_m_u_m_s_geq_2} we can conclude that
        \begin{equation*}
          a_m(u_m)^{s/2} \mbox{ is bounded in } L^{10/3}(0,T;L^{10/3}(\Omega)),
        \end{equation*}
        that is,
        \begin{equation}
          a_m(u_m) \mbox{ is bounded in } L^{5s/3}(0,T;L^{5s/3}(\Omega)).
          \label{limitacao_a_m_u_m_5s/3_s_geq_2}
        \end{equation}
        
        For each fixed $m \in \mathbb{N}$, consider the zero measure set $\mathcal{N} \subset (0, \infty)$ such that
        \begin{equation*}
          u_m(t^{\ast},\cdot), v_m(t^{\ast},\cdot) \in H^1(\Omega), \ \forall t^{\ast} \in (0,\infty) \setminus \mathcal{N}.
        \end{equation*}
        Then, for each fixed $t^{\ast} \in (0,\infty) \setminus \mathcal{N}$, let us consider the sets
        \begin{equation*}
          \{ 0 \leq u_m \leq 1 \} = \Big \{ x \in \Omega \ \Big | \ 0 \leq u_m(t^{\ast},x) \leq 1 \Big \}
        \end{equation*}
        and
        \begin{equation*}
          \{ u_m \geq 1 \} = \Big \{ x \in \Omega \ \Big | \ u_m(t^{\ast},x) \geq 1 \Big \}.
        \end{equation*}
        Now note that, since $s \geq 2$, we have
        \begin{align*}
        &  \int_{\Omega} a_m(u_m(t^{\ast},x)^2 \norm{\nabla v_m(t^{\ast},x)}{}^2 \ dx \\
      &    \leq \int_{\{ 0 \leq u_m \leq 1 \}}{\norm{\nabla v_m(t^{\ast},x)}{}^2 \ dx} + \int_{\{ u_m \geq 1 \}}{a_m(u_m(t^{\ast},x))^s \norm{\nabla v_m(t^{\ast},x)}{}^2 \ dx} \\
       &   \leq \int_{\Omega}{\norm{\nabla v_m(t^{\ast},x)}{}^2 \ dx} + \int_{\Omega}{a_m(u_m(t^{\ast},x))^s \norm{\nabla v_m(t^{\ast},x)}{}^2 \ dx}.
        \end{align*}
        The last inequality is valid for all $t^{\ast} \in (0,\infty) \setminus \mathcal{N}$, then if we integrate in the variable $t$ we obtain
        \begin{align*}
          \int_0^\infty{\int_{\Omega}{a_m(u_m(t,x)^2 \norm{\nabla v_m(t,x)}{}^2 \ dx} \ dt} & \leq \int_0^\infty{\int_{\Omega}{\norm{\nabla v_m(t,x)}{}^2 \ dx} \ dt} \\
          &  + \int_0^\infty{\int_{\Omega}{a_m(u_m(t,x))^s \norm{\nabla v_m(t,x)}{}^2 \ dx} \ dt}.
        \end{align*}
        Therefore by \eqref{limitacao_uniforme_nabla_v_m} and \eqref{limitacao_aux_Delta_v_m_s_intermediario} we can conclude that
        \begin{equation}
        a_m(u_m) \nabla v_m \mbox{ is bounded in } L^2(0,\infty;L^2(\Omega)).
        \label{limitacao_termo_misto_u_v_m_s_geq_2}
      \end{equation}
      
      Now we test the $u_m$-equation of problem \eqref{problema_P_m_intro} by $u_m$. This gives us
      \begin{align*}
        \dfrac{1}{2} \dfrac{d}{dt} \norma{u_m}{L^2(\Omega)}^2 & + \norma{\nabla u_m}{L^2(\Omega)}^2 = \int_{\Omega}{a_m(u_m) \nabla v_m \cdot \nabla u_m \ dx} \\
        & \leq \dfrac{1}{2} \int_{\Omega}{a_m(u_m)^2 \norm{\nabla v_m}{}^2 \ dx} + \dfrac{1}{2} \norma{\nabla u_m}{L^2(\Omega)}^2,
      \end{align*}
      hence we have
      \begin{equation*}
        \dfrac{d}{dt} \norma{u_m}{L^2(\Omega)}^2 + \norma{\nabla u_m}{L^2(\Omega)}^2 \leq \int_{\Omega}{a_m(u_m)^2 \norm{\nabla v_m}{}^2 \ dx}.
      \end{equation*}
      Integrating with respect to $t$, we conclude from \eqref{limitacao_termo_misto_u_v_m_s_geq_2} that
      \begin{equation}
        u_m \mbox{ is bounded in } L^{\infty}(0,\infty;L^2(\Omega))
        \label{limitacao_u_m_s_geq_2}
      \end{equation}
      and
      \begin{equation}
        \nabla u_m \mbox{ is bounded in } L^2(0,\infty;L^2(\Omega)).
        \label{limitacao_nabla_u_m_s_geq_2}
      \end{equation}
      
      Then, if we consider the $u_m$-equation of \eqref{problema_P_m_intro}, by applying \eqref{limitacao_nabla_u_m_s_geq_2} and \eqref{limitacao_termo_misto_u_v_m_s_geq_2} we conclude that
      \begin{equation}
        \partial_t u_m \mbox{ is bounded in } L^2(0,\infty;(H^1(\Omega))').
        \label{limitacao_u_m_derivada_tempo_s_geq_2}
      \end{equation}
      
      Let $(a_m(u_m)^s)^{\ast} = \dfrac{1}{\norm{\Omega}{}} \D{\int_{\Omega}}{a_m(u_m)^s \ dx}$, from \eqref{limitacao_a_m_u_m_s_geq_2} and \eqref{limitacao_nabla_a_m_u_m_s_geq_2}, we can also conclude that
        \begin{equation*}
          \nabla a_m(u_m)^s \mbox{ is bounded in } L^2(0,\infty;L^1(\Omega)).
        \end{equation*}
        In view of Lemma \ref{lema_desigualdade_poincare_media_nula}, the latter implies
        \begin{equation*}
          a_m(u_m)^s - (a_m(u_m)^s)^{\ast} \mbox{ is bounded in } L^2(0,\infty;W^{1,1}(\Omega))
        \end{equation*}
        and, in particular, by the Sobolev embedding, we have
        \begin{equation}
          a_m(u_m)^s - (a_m(u_m)^s)^{\ast} \mbox{ is bounded in } L^2(0,\infty;L^{3/2}(\Omega)).
          \label{limitacao_a_m_u_m_elevado_a_s_L_3/2}
        \end{equation}
        
        Now we consider the $v_m$-equation of \eqref{problema_P_m_intro} written as
        \begin{equation}
          \partial_t v_m - \Delta v_m + (a_m(u_m)^s)^{\ast} v_m = - (a_m(u_m)^s - (a_m(u_m)^s)^{\ast}) v_m.
          \label{eq_v_m_equivalente_s_geq_2}
        \end{equation}
        Testing \eqref{eq_v_m_equivalente_s_geq_2} by $v_m$ and using Hölder's inequality we can obtain
        \begin{align*}
          \frac{1}{2} \frac{d}{dt} \norma{v_m}{L^2(\Omega)}^2 & + \norma{\nabla v_m}{L^2(\Omega)}^2 + (a_m(u_m)^s)^{\ast} \norma{v_m}{L^2(\Omega)}^2 \\
          & \leq C \norma{v^0_m}{L^{\infty}(\Omega)} \norma{a_m(u_m)^s - (a_m(u_m)^s)^{\ast}}{L^{3/2}(\Omega)} \norma{v_m}{L^3(\Omega)} \\
          & \leq C \norma{a_m(u_m)^s - (a_m(u_m)^s)^{\ast}}{L^{3/2}(\Omega)}^2 + \delta \norma{v_m}{L^2(\Omega)}^2 + \delta \norma{\nabla v_m}{L^2(\Omega)}^2.
        \end{align*}
        In order to bound $(a_m(u_m)^s)^{\ast}$ from below, we will apply Lemma \ref{lema_integral_de_a_m}. Indeed,
        \begin{equation*}
          (a_m(u_m)^s)^{\ast} \geq C \Big( \int_{\Omega}{a_m(u_m) \ dx} \Big )^s
        \end{equation*}
        and applying Lemma \ref{lema_integral_de_a_m} (with $w_m = u_m$, $p = 2$, and using \eqref{limitacao_u_m_s_geq_2}) we conclude that there exist $\beta > 0$ and $m_0$ large enough such that $(a_m(u_m)^s)^{\ast} \geq \beta > 0$, $a.e. \ t \in (,\infty)$, for all $m \geq m_0$. Therefore
        \begin{align*}
          \frac{1}{2} \frac{d}{dt} \norma{v_m}{L^2(\Omega)}^2 + (1 - \delta) \norma{\nabla v_m}{L^2(\Omega)}^2 + (\beta - \delta) \norma{v_m}{L^2(\Omega)}^2 \leq C \norma{a_m(u_m)^s - (a_m(u_m)^s)^{\ast}}{L^{3/2}(\Omega)}^2.
        \end{align*}
        Now, choosing $\delta$ small enough, integrating the last inequality with respect to $t$ and using \eqref{limitacao_a_m_u_m_elevado_a_s_L_3/2} we obtain
        \begin{equation}
          v_m \mbox{ is bounded in } L^2(0,\infty;L^2(\Omega)).
          \label{estimativa_v_m_L^2_s_geq_2}
        \end{equation}
        
        With the $m$-independent \emph{a priori} bounds obtained so far we can also give an $m$-independent \emph{a priori} bound for $\partial_t v_m$. In fact, if we consider again the equation \eqref{eq_v_m_equivalente_s_geq_2}, then the $m$-independent estimate in the $L^{\infty}$-norm for $v_m$ given by Lemma \ref{lemma_positividade_u_v_m}-2 and the $m$-independent \emph{a priori} bounds
        \eqref{estimativa_v_m_L^2_s_geq_2}, \eqref{limitacao_a_m_u_m_elevado_a_s_L_3/2} and \eqref{limitacao_aux_Delta_v_m_s_intermediario} allow us to conclude that
        \begin{equation*}
          \partial_t v_m \mbox{ is bounded in } L^2(0,\infty;L^{3/2}(\Omega)).
        \end{equation*}
      
      Now, using \eqref{limitacao_u_m_s_geq_2}, \eqref{limitacao_nabla_u_m_s_geq_2} and \eqref{limitacao_u_m_derivada_tempo_s_geq_2} we can conclude that there is a subsequence of $\{ u_m \}$, still denoted by $\{ u_m \}$, and a limit function $u$ such that
      \begin{equation*}
        \begin{array}{c}
          u_m \longrightarrow u \mbox{ weakly* in } L^{\infty}(0,\infty;L^2(\Omega)), \\
          \nabla u_m \longrightarrow \nabla u \mbox{ weakly in } L^2(0,\infty;L^2(\Omega)), \\
          \partial_t u_m \longrightarrow u \mbox{ weakly in } L^2 \Big ( 0,\infty;\big (H^1(\Omega) \big )' \Big ).
        \end{array}
      \end{equation*}
      By applying the compactness result Lemma \ref{lema_Simon}, one has
      \begin{equation*}
        u_m \longrightarrow u \mbox{ strongly in } L^2(0,T;L^2(\Omega)), \ \forall T > 0.
      \end{equation*}
      Using the Dominated Convergence Theorem and \eqref{limitacao_a_m_u_m_5s/3_s_geq_2} we can also prove that
      \begin{equation*}
        a_m(u_m) \longrightarrow u \mbox{ strongly in } L^p(0,T;L^p(\Omega)), \ \forall p \in (1,5s/3),
      \end{equation*}
      and using Lemma \ref{lema_convergencia_w_elevado_a_s},
      \begin{equation*}
        a_m(u_m)^s \longrightarrow u^s \mbox{ strongly in } L^p(0,T;L^p(\Omega)), \ \forall p \in (1,5/3).
      \end{equation*}
      
      From the global in time estimate \eqref{limitacao_a_m_u_m_s_geq_2} we can conclude that, up to a subsequence,
        \begin{equation*}
          a_m(u_m) \to u \mbox{ weakly* in } L^{\infty}(0,\infty;L^s(\Omega)),
       \end{equation*}
       hence, in particular,
       \begin{equation*}
         u \in L^{\infty}(0,\infty;L^s(\Omega)).
       \end{equation*}
       
       For $s \geq 2$, if we consider the functions $v_m$, we have the same $m$-independent estimates that we had for $s \in [1,2)$. Then we have the same convergences given in \eqref{convergencia_v_m_s_intermediario}.
       
       Following the ideas of Subsection \ref{subsec:estimativas_u_m_s<2}, we can identify the limits of $a_m(u_m) \nabla v_m$ and $a_m(u_m)^s v_m$ with $u \nabla v$ and $u^s v$, respectively.
       
       This finishes the proof of existence of solution to the original problem \eqref{problema_P} as a limit of solutions of the regularized problems \eqref{problema_P_m_intro} for $s \geq 2$.
      

\section{Regularity and Uniqueness in \texorpdfstring{$\boldsymbol{2}$}{2}D}
    
    \label{section: original problem 2D}
      
      In this section, we show that, for two dimensional domains, we can improve the results on the uniqueness and regularity of the solution of \eqref{problema_P}. The key point is the inequality \eqref{desig_ladyzhenskaya}, which allows us to improve the \emph{a priori} estimates of $u_m$ and then of $v_m$, where $(u_m,v_m)$ is the solution of \eqref{problema_P_m_intro}.
      
      \subsection{Uniqueness in \texorpdfstring{$\boldsymbol{2}$}{2}D}
      \begin{teo} \label{teo_unicidade_2_D}
        In the two dimensional case, 
        we have uniqueness of solution in the class of functions $(u,v)$ such that
        \begin{equation}
          u \in L^{\infty}(0,T;L^2(\Omega)) \cap L^2(0,T;H^1(\Omega))
          \cap L^{4s - 4}(0,T;L^{4s - 4}(\Omega)),
          \label{regularidade_u_para_unicidade_2_dim_c}
        \end{equation}
         \begin{equation} \label{regularidade_u_para_unicidade_s=2_extra}
          u \in L^{4}(0,T;L^{4+\epsilon}(\Omega))  \quad \hbox{if $s=2$}
        \end{equation}
        and
        \begin{equation}
          v \in L^{\infty}(0,T;L^{\infty}(\Omega)) \cap
          L^{\infty}(0,T;H^1(\Omega)) \cap L^2(0,T;H^2(\Omega)).
          \label{regularidade_v1_para_unicidade_2_dim}
        \end{equation}
      \end{teo}
      \begin{obsv}
        The regularities \eqref{regularidade_u_para_unicidade_2_dim_c} and \eqref{regularidade_v1_para_unicidade_2_dim} imply in particular that
        \begin{equation*}
          u_t \in L^2(0,T;(H^1(\Omega))')
        \end{equation*}
        and therefore, the solution $u$ can be taken as test function in the $u$-equation of \eqref{problema_P}.
        
        The regularity $u\in L^{4s - 4}(0,T;L^{4s - 4}(\Omega))$ is an additional hypothesis only if $s>2$, because in $2D$ domains, the regularity $u \in L^{\infty}(0,T;L^2(\Omega)) \cap L^2(0,T;H^1(\Omega))$ implies $u\in L^{4}(0,T;L^{4}(\Omega))$.
        \hfill \qedsymbol
      \end{obsv}
      
      \begin{proof}[\bf Proof of Theorem \ref{teo_unicidade_2_D}]
        Suppose $(u_1,v_1)$ and $(u_2,v_2)$ are two solutions of the original problem \eqref{problema_P} with the regularity given in \eqref{regularidade_u_para_unicidade_2_dim_c}-\eqref{regularidade_v1_para_unicidade_2_dim}. Define $(u,v) = (u_2 - u_1, v_2 - v_1)$. Then $(u,v)$ satisfies
        \begin{equation}
          \begin{array}{rl}
            \prodl{u_t(t)}{\varphi} + \prodl{\nabla u(t)}{\nabla \varphi} = \prodl{u(t) \nabla v_2(t)}{\nabla \varphi} + \prodl{u_1(t) \nabla v}{\nabla \varphi}, \ \forall \varphi \in H^1(\Omega), \ a.e \ t \in (0,T), 
          \end{array}
          \label{equacao_u_fraca_unicidade}
        \end{equation}
        and 
      \begin{equation}
        \begin{array}{rl}
          v_t(t) - \Delta v(t) = - [(u_2(t))^s - (u_1(t))^s] v_2(t) - (u_1(t))^s v(t), \ a.e \ t \in (0,T), 
        \end{array}
        \label{equacao_v_forte}
      \end{equation}
      with $(u(0),v(0)) = (0,0)$. Note that we can conclude from \eqref{equacao_u_fraca_unicidade} that $u$ is a zero mean function.
      
      Now we test \eqref{equacao_u_fraca_unicidade} by $u$. We obtain, first using the interpolation inequality from Lemma \ref{lema_desig_lady2}-2, for $2$D domains, and the Poincaré inequality for zero mean functions from Lemma \ref{lema_desigualdade_poincare_media_nula}, and after Young's inequality
      \begin{equation*}
        \begin{array}{l} 
          \displaystyle
          \frac12 \frac{d}{dt} \norma{u}{L^2(\Omega)}^2 +  \norma{\nabla u}{L^2(\Omega)}^2 \leq \norma{u}{L^4(\Omega)} \norma{\nabla v_2}{L^4(\Omega)} \norma{\nabla u}{L^2(\Omega)} + \norma{u_1}{L^4(\Omega)} \norma{\nabla v}{L^4(\Omega)} \norma{\nabla u}{L^2(\Omega)} \\ 
          \displaystyle
          \qquad \leq C \norma{u}{L^2(\Omega)}^{1/2} \norma{\nabla v_2}{L^4(\Omega)} \norma{\nabla u}{L^2(\Omega)}^{3/2} + C \norma{u_1}{L^4(\Omega)} \norma{\nabla v}{L^2(\Omega)}^{1/2} \norma{\nabla v}{H^1(\Omega)}^{1/2} \norma{\nabla u}{L^2(\Omega)} \\
          \displaystyle
          \qquad \leq C(\delta) \norma{\nabla v_2}{L^4(\Omega)}^4 \norma{u}{L^2(\Omega)}^2 
          + C(\delta) \norma{u_1}{L^4(\Omega)}^4 \norma{\nabla v}{L^2(\Omega)}^2 
          + \delta \norma{\nabla u}{L^2(\Omega)}^2 
          + \delta \norma{\nabla v}{H^1(\Omega)}^2,
        \end{array}
      \end{equation*}
      for each $\delta > 0$. Then, accounting for \eqref{obsv_norma_de_nabla_v}, we get
      \begin{equation} \label{unicidade_u_2D_s<2}
        \begin{array}{l} 
        \displaystyle
        \dfrac{1}{2} \frac{d}{dt} \norma{u}{L^2(\Omega)}^2 + \norma{\nabla u}{L^2(\Omega)}^2 \leq C(\delta) \norma{\nabla v_2}{L^4(\Omega)}^4 \norma{u}{L^2(\Omega)}^2 \\
        \displaystyle
        \qquad + C(\delta) \norma{u_1}{L^4(\Omega)}^4 \norma{\nabla v}{L^2(\Omega)}^2 
        + \delta \norma{\nabla u}{L^2(\Omega)}^2
        + \delta C \norma{\Delta v}{L^2(\Omega)}^2.
      \end{array}
    \end{equation}
    
    Next we test \eqref{equacao_v_forte} by $v - \Delta v$. Taking into account that $v_2 \in L^{\infty}(\Omega)$, we obtain
    \begin{align*}
      & \frac{1}{2} \frac{d}{dt} \norma{v}{H^1(\Omega)}^2 + \norma{\Delta v}{L^2(\Omega)}^2 
      + \norma{\nabla v}{L^2(\Omega)}^2 + \int_{\Omega} u_1^{s} v^2 \ dx
      \\ &  \leq \int_{\Omega}{\norm{(u_2)^s - (u_1)^s}{} \norm{v_2}{} \norm{v - \Delta v}{}} \ dx + \int_{\Omega}{\norm{u_1^{s/2} v}{} \norm{\Delta v}{}} \ dx \\
      &  \leq C(\delta) \norma{((u_2)^s - (u_1)^s)}{L^2(\Omega)}^2 + C(\delta) \norma{(u_1)^s v}{L^2(\Omega)}^2 + \delta \norma{v}{L^2(\Omega)}^2 + 2 \delta \norma{\Delta v}{L^2(\Omega)}^2,
    \end{align*}
    for each $\delta > 0$. We must estimate the terms $\norma{((u_2)^s - (u_1)^s)}{L^2(\Omega)}^2$ and $\norma{(u_1)^s v}{L^2(\Omega)}^2$. For the first of these two terms we use
    \begin{equation*}
      \norm{u_2^s - u_1^s}{} \leq s \norm{u_2 + u_1}{}^{s-1} \norm{u_2 - u_1}{},
    \end{equation*}
    from Lemma \ref{lema_a_m_elevado_a_s} and, considering Lemma \ref{lema_desigualdade_poincare_media_nula} applied to zero mean function $u$, we find
    \begin{align*}
      & \norma{[(u_2)^s - (u_1)^s]}{L^2(\Omega)}^2 = \int_{\Omega}{\norm{(u_2)^s - (u_1)^s}{}^2 \ dx} \\
      & \qquad \leq s^2 \int_{\Omega}{\norm{u_2 + u_1}{}^{2s - 2} \norm{u_2 - u_1}{}^2 \ dx} \leq s^2 \norma{u_2 + u_1}{L^{4s-4}(\Omega)}^{2s - 2} \norma{u}{L^4(\Omega)}^2 \\
      & \qquad \leq \norma{u_2 + u_1}{L^{4s-4}(\Omega)}^{2s - 2} \norma{u}{L^2(\Omega)} \norma{\nabla u}{L^2(\Omega)} \leq \norma{u_2 + u_1}{L^{4s-4}(\Omega)}^{4s - 4} \norma{u}{L^2(\Omega)}^2 + \delta \norma{\nabla u}{L^2(\Omega)}^2.
    \end{align*}
    For the second term we have, for any $\varepsilon > 0$,
    \begin{align*}
      & \norma{(u_1)^s v}{L^2(\Omega)}^2 = \int_{\Omega}{(u_1)^{2s} v^2 \ dx} \leq \norma{u_1}{L^{2s+\varepsilon}(\Omega)}^{2s} \norma{v}{L^{(2s+\varepsilon)/\varepsilon}(\Omega)}^2 \leq C(\varepsilon) \norma{u_1}{L^{2s+\varepsilon}(\Omega)}^{2s} \norma{v}{H^1(\Omega)}^2.
    \end{align*}
    Using the estimates of these two terms we obtain
    \begin{equation} \label{unicidade_v_2_D_forte}
      \begin{array}{rl}
        & \dfrac{1}{2} \dfrac{d}{dt} \norma{v}{H^1(\Omega)}^2 + \norma{\Delta v}{L^2(\Omega)}^2 \leq C(\delta) \norma{u_2 + u_1}{L^{4s-4}(\Omega)}^{4s - 4} \norma{u}{L^2(\Omega)}^2 \\
        & \qquad + \delta \norma{\nabla u}{L^2(\Omega)}^2 + C(\delta,\varepsilon) \norma{u_1}{L^{2s+\varepsilon}(\Omega)}^{2s} \norma{v}{H^1(\Omega)}^2 + \delta \norma{v}{L^2(\Omega)}^2 + \delta \norma{\Delta v}{L^2(\Omega)}^2.
      \end{array}
    \end{equation}
    If we sum up \eqref{unicidade_u_2D_s<2} and \eqref{unicidade_v_2_D_forte} and choose $\delta > 0$ small enough so that the terms that are multiplied by $\delta$ on the right hand side can be absorbed by the corresponding nonnegative terms on the left hand side, we obtain
    \begin{equation*}
      \begin{array}{rl}
        & \dfrac{1}{2} \dfrac{d}{dt} \Big ( \norma{u}{L^2(\Omega)}^2 + \norma{v}{H^1(\Omega)}^2 \Big ) \leq C \norma{\nabla v_2}{L^4(\Omega)}^4 \norma{u}{L^2(\Omega)}^2 \\
        \displaystyle
        & \qquad + C \norma{u_1}{L^4(\Omega)}^4 \norma{\nabla v}{L^2(\Omega)}^2 + C \norma{u_2 + u_1}{L^{4s-4}(\Omega)}^{4s - 4} \norma{u}{L^2(\Omega)}^2 \\
        & \qquad + C(\varepsilon) \norma{u_1}{L^{2s+\varepsilon}(\Omega)}^{2s} \norma{v}{H^1(\Omega)}^2 + C \norma{v}{L^2(\Omega)}^2
      \end{array}
    \end{equation*}
    Now, taking into account that
    \begin{equation*}
      \norma{v}{L^2(\Omega)}^2, \norma{\nabla v}{L^2(\Omega)}^2 \leq \norma{v}{H^1(\Omega)}^2
    \end{equation*}
    and grouping the common factors, we have
    \begin{equation}\label{unicidade_u_v_2_D}
      \begin{array}{rl}
        & \dfrac{1}{2} \dfrac{d}{dt} \Big ( \norma{u}{L^2(\Omega)}^2 + \norma{v}{H^1(\Omega)}^2 \Big ) \leq C ( \norma{\nabla v_2}{L^4(\Omega)}^4 + \norma{u_2 + u_1}{L^{4s-4}(\Omega)}^{4s - 4}) \norma{u}{L^2(\Omega)}^2 \\
        \displaystyle
        & \qquad + (C \norma{u_1}{L^4(\Omega)}^4 + C(\varepsilon) \norma{u_1}{L^{2s+\varepsilon}(\Omega)}^{2s} + C) \norma{v}{H^1(\Omega)}^2.
      \end{array}
    \end{equation}
    Finally, we recall from the regularity hypotheses that we have, in particular,
    \begin{equation*}
      u_1, u_2 , \nabla v_1, \nabla v_2 \in L^4(0,T;L^4(\Omega)) \mbox{ and } u_1, u_2 \in L^{4s-4}(0,T;L^{4s-4}(\Omega)).
    \end{equation*}
    Therefore, it suffices to verify that there exists $\varepsilon > 0$ small enough such that
    \begin{equation} \label{regularidade_necesaria_para_u_1_u_2}
      u_1, u_2 \in L^{2s}(0,T;L^{2s+\varepsilon}(\Omega)).
    \end{equation}
    For $s \in [1,2)$, since $u_1,u_2 \in L^4(0,T;L^4(\Omega))$, then in particular one has \eqref{regularidade_necesaria_para_u_1_u_2}. For $s = 2$, \eqref{regularidade_necesaria_para_u_1_u_2} is in fact the  hypothesis \eqref{regularidade_u_para_unicidade_s=2_extra}. And for $s > 2$, hypothesis $u_1,u_2\in L^{4s-4}(0,T;L^{4s-4}(\Omega))$ implies in particular that there is a $\varepsilon > 0$ such that \eqref{regularidade_necesaria_para_u_1_u_2} holds. Therefore, recalling that $u(0) = v(0) = 0$, we are able to apply Gronwall's inequality (Lemma \ref{lema_gronwall}) to \eqref{unicidade_u_v_2_D} and conclude that $u = v = 0$, that is, $u_1 = u_2$ and $v_1 = v_2$.
  \end{proof}

      
      \subsection{Proof of Theorem \ref{teo_2D_existencia_unicidade_estimativas}}
      
        In the two dimensional case, we can study the regularity of the solution $(u,v)$ of \eqref{problema_P} for all $s \geq 1$ at the same time.These solutions can be obtained as a limit of the regularized solutions $(u_m,v_m)$ of \eqref{problema_P_m_intro} as $m \to \infty$, considering initial data $(u^0,v^0) \in H^2(\Omega) \times H^2(\Omega)$. In this case, it is not necessary to regularize the initial data, taking directly $(u_m^0,v_m^0)=(u^0,v^0)$.  
        
        In order to prove that the solution $(u,v)$ of \eqref{problema_P} provided by Theorem \ref{teo_3D_existencia} is in fact more regular, it suffices to prove the corresponding extra $m$-independent estimates for $(u_m,v_m)$ in the spaces given in Theorem \ref{teo_2D_existencia_unicidade_estimativas}.
        
        We take $u_m^p$, for any $1 \leq p < \infty$, as a test function in the $u_m$-equation to obtain
        \begin{align*}
          \dfrac{1}{p + 1} \dfrac{d}{dt} & \int_{\Omega}{u_m^{p+1}(x) \ dx} + p \int_{\Omega}{u_m^{p-1}(x) \norm{\nabla u_m(x)}{}^2 \ dx} \\
          & = p \int_{\Omega}{a_m(u_m(x)) \nabla v(x) \cdot \nabla u_m(x) u_m^{p-1}(x) \ dx} \\
          & \leq p \int_{\Omega}{u_m^p(x) \norm{\nabla v_m(x)}{} \norm{\nabla u_m(x)}{} \ dx} \\
          & \leq p \int_{\Omega}{u_m^{p/2 + 1/2}(x) \norm{\nabla v_m(x)}{} u_m^{p/2 - 1/2} \norm{\nabla u_m(x)}{} \ dx} \\
          & \leq p \norma{u_m^{p/2 + 1/2}}{L^4(\Omega)} \norma{\nabla v_m}{L^4(\Omega)} \Big ( \int_{\Omega}{u_m^{p-1}(x) \norm{\nabla u_m(x)}{}^2 \ dx} \Big )^{1/2} \\
          & \leq C p \norma{u_m^{p/2 + 1/2}}{L^2(\Omega)} \norma{\nabla v_m}{L^4(\Omega)} \Big ( \int_{\Omega}{u_m^{p-1}(x) \norm{\nabla u_m(x)}{}^2 \ dx} \Big )^{1/2} \\
          & \ + C p \sqrt{p + 1} \norma{u_m^{p/2 + 1/2}}{L^2(\Omega)}^{1/2} \norma{\nabla v_m}{L^4(\Omega)} \Big ( \int_{\Omega}{u_m^{p-1}(x) \norm{\nabla u_m(x)}{}^2 \ dx} \Big )^{3/4}.
        \end{align*}
        By using Young's inequality, we obtain
        \begin{align*}
          & \dfrac{1}{p + 1} \dfrac{d}{dt} \int_{\Omega}{u_m^{p+1}(x) \ dx} + \dfrac{p}{2} \int_{\Omega}{u_m^{p-1}(x) \norm{\nabla u_m(x)}{}^2 \ dx}  \\
          & \leq C p \norma{\nabla v_m}{L^4(\Omega)}^2 \int_{\Omega}{u_m^{p+1}(x) \ dx}  + C p (p+1)^2 \norma{\nabla v_m}{L^4(\Omega)}^4 \int_{\Omega}{u_m^{p+1}(x) \ dx}.
        \end{align*}
        By estimates \eqref{limitacao_aux_Dv_m_s_intermediario} and \eqref{limitacao_aux_Dv_m_s_intermediario_2} we have 
        \begin{equation*}
          \norma{\nabla v_m}{L^4(\Omega)}^2, \norma{\nabla v_m}{L^4(\Omega)}^4 \mbox{ are bounded in } L^1(0,\infty),
        \end{equation*}
        then, Gronwall's inequality (Lemma \ref{lema_gronwall}) leads us to
        \begin{equation} \label{limitacao_u_m_L_p_a}
          u_m \mbox{ is bounded in } L^{\infty}(0,\infty;L^{p+1}(\Omega)), \mbox{ for } 1 \leq p < \infty,
        \end{equation}
        \begin{equation} \label{limitacao_u_m_L_p_b}
          u_m^{(p-1)/2} \nabla u_m \mbox{ is bounded in } L^2(0,\infty;L^2(\Omega)), \mbox{ for } 1 \leq p < \infty.
        \end{equation}
        \begin{obsv}
          The $m$-independent bounds in \eqref{limitacao_u_m_L_p_a} and \eqref{limitacao_u_m_L_p_b} depend exponentially on $p(p+1)^2$.
          \hfill \qedsymbol
        \end{obsv}
        
        Then we recall from \eqref{regularidade_solucao_u_m_v_m_reg_H_2} that
      \begin{equation*}
          \partial_t v_m, \Delta v_m \in L^2(0,\infty;H^1(\Omega)).
        \end{equation*}
        Hence,  the following system is satisfied a.e. $(t,x) \in (0, \infty) \times \Omega$:
        \begin{equation} \label{equacao_grad_v_m}
          \nabla (\partial_t v_m) - \nabla \Delta v_m  = - s \ a_m(u_m)^{s-1} \nabla a_m(u_m) v_m - a_m(u_m)^s \nabla v_m.
        \end{equation}
        Since all the terms in this system are in $L^2(0,\infty;L^2(\Omega))$, we can take the inner product with $- \nabla \Delta v_m \in L^2(0,\infty;L^2(\Omega))$ and integrate over $\Omega$. Using integration by parts, Hölder's and Young's inequalities and the estimate of $v_m$ in $L^{\infty}(0,\infty;L^{\infty}(\Omega))$, we obtain
        \begin{equation} \label{estimativa_Delta_v_m_2}
          \dfrac{d}{dt} \norma{\Delta v_m}{L^2(\Omega)}^2 + \norma{\nabla \Delta v_m}{L^2(\Omega)}^2 \leq C \norma{a_m(u_m)^{s-1} \nabla a_m(u_m)}{L^2(\Omega)}^2 + C \norma{a_m(u_m)}{L^{4s}(\Omega)}^{2s} \norma{\nabla v_m}{L^4(\Omega)}^2.
        \end{equation}
        Then, integrating \eqref{estimativa_Delta_v_m_2} with respect to $t$ and using \eqref{limitacao_u_m_L_p_b} with $p = 2s - 1$, \eqref{limitacao_u_m_L_p_a} with $p = 4s - 1$ and \eqref{limitacao_aux_Dv_m_s_intermediario} we conclude that
        \begin{equation}
          \begin{array}{c}
            v_m \mbox{ is bounded in } L^{\infty}(0,\infty;H^2(\Omega)), \\
            \Delta v_m \mbox{ is bounded in } L^2(0,\infty;H^1(\Omega)).
          \end{array}
          \label{limitacao_v_m_final_2D}
        \end{equation}
        Now we can use \eqref{equacao_grad_v_m} to write
        \begin{equation*}
          \nabla \partial_t v_m = \nabla \Delta v_m - s \ a_m(u_m)^{s-1} \nabla a_m(u_m) v_m - a_m(u_m)^s \nabla v_m.
        \end{equation*}
        Then, using the estimate in the $L^{\infty}$-norm for $v_m$ given by Lemma \ref{lemma_positividade_u_v_m}-2, \eqref{limitacao_u_m_L_p_a} with $p = 4s - 1$, \eqref{limitacao_u_m_L_p_b} $p = 2s - 1$ and \eqref{limitacao_v_m_final_2D}, we also conclude that
        \begin{equation}
          \partial_t v_m \mbox{ is bounded in } L^2(0,\infty;H^1(\Omega)).
          \label{limitacao_derivada_v_m_2D}
        \end{equation}
        
        Then, because of the regularity of the solutions $(u_m,v_m)$, we can use $- \Delta u_m \in L^2(0,\infty;L^2(\Omega))$ as a test function in the $u_m$-equation of \eqref{problema_P_m_intro},
        \begin{align*}
          \dfrac{1}{2} \dfrac{d}{dt} \norma{\nabla u_m}{L^2(\Omega)}^2 + \norma{\Delta u_m}{L^2(\Omega)}^2 = \int_{\Omega}{a_m(u_m) \Delta v_m \Delta u_m \ dx} + \int_{\Omega}{\nabla a_m(u_m) \cdot \nabla v_m \Delta u_m \ dx} \\
          \leq \norma{u_m}{L^4(\Omega)} \norma{\Delta v_m}{L^4(\Omega)} \norma{\Delta u_m}{L^2(\Omega)} + \norma{\nabla u_m}{L^4(\Omega)} \norma{\nabla v_m}{L^4(\Omega)} \norma{\Delta u_m}{L^2(\Omega)} \\
          \leq \norma{u_m}{L^4(\Omega)} \norma{\Delta v_m}{L^4(\Omega)} \norma{\Delta u_m}{L^2(\Omega)} + C \norma{\nabla u_m}{L^2(\Omega)}^{1/2} \norma{\nabla v_m}{L^4(\Omega)} \norma{\Delta u_m}{L^2(\Omega)}^{3/2}.
        \end{align*}
        Using Young's inequality we get
        \begin{equation}
          \begin{array}{rl}
            \dfrac{d}{dt} \norma{\nabla u_m}{L^2(\Omega)}^2 + \norma{\Delta u_m}{L^2(\Omega)}^2 \leq C \norma{u_m}{L^4(\Omega)}^2 \norma{\Delta v_m}{L^4(\Omega)}^2 + C \norma{\nabla v_m}{L^4(\Omega)}^4 \norma{\nabla u_m}{L^2(\Omega)}^2.
          \end{array}
          \label{estimativa_Delta_u_m_1}
        \end{equation}
        
        Therefore, using \eqref{limitacao_aux_Dv_m_s_intermediario}, \eqref{limitacao_u_m_L_p_a} for $p = 3$, \eqref{limitacao_v_m_final_2D} and Gronwall's inequality (Lemma \ref{lema_gronwall}) in \eqref{estimativa_Delta_u_m_1} we conclude that
        \begin{equation}
          \nabla u_m \mbox{ is bounded in } L^{\infty}(0,\infty;L^2(\Omega)).
          \label{limitacao_nabla_u_m_2D}
        \end{equation}
        \begin{equation}
          \Delta u_m \mbox{ is bounded in } L^2(0,\infty;L^2(\Omega)).
          \label{limitacao_Delta_u_m_2D}
        \end{equation}
        
        Considering the $u_m$equation of \eqref{problema_P_m_intro},
        \begin{equation*}
          \partial_t u_m - \Delta u_m = - a_m(u_m) \Delta v_m - \nabla a_m(u_m) \nabla v_m 
        \end{equation*}
        and estimates \eqref{limitacao_u_m_L_p_a}, \eqref{limitacao_v_m_final_2D}, \eqref{limitacao_nabla_u_m_2D} and \eqref{limitacao_Delta_u_m_2D} we obtain
        \begin{equation}
          \partial_t u_m \mbox{ is bounded in } L^2(0,\infty;L^2(\Omega)),
          \label{limitacao_derivada_u_m_2D}
        \end{equation}
        
        We can use \eqref{limitacao_u_m_L_p_a},
        \eqref{limitacao_nabla_u_m_2D}, \eqref{limitacao_Delta_u_m_2D}, \eqref{limitacao_v_m_final_2D}, \eqref{limitacao_derivada_u_m_2D} and \eqref{limitacao_derivada_v_m_2D} and compactness results in the weak*, weak and strong topologies and the uniqueness of the limit problem \eqref{problema_P} for functions satisfying \eqref{regularidade_u_para_unicidade_2_dim_c}-\eqref{regularidade_v1_para_unicidade_2_dim} to conclude that there is a unique limit $(u,v)$ satisfying \eqref{problema_P} $a.e.$ in $(0,\infty) \times \Omega$.
        
        \
        
        Finally, if now we suppose that $\Omega$ has the $W^{2,3}$-regularity then we can prove better $m$-independent estimates for $u_m$. We would like to test the $u_m$-equation of \eqref{problema_P_m_intro} by $\Delta^2 u_m$, but we do not have enough regularity about $\Delta^2 u_m$. Instead of it, we argue as in \eqref{equacao_grad_v_m}, first we take the gradient of the $u_m$-equation of \eqref{problema_P_m_intro} and after we test the resulting equation by $\nabla \Delta u_m$.
      
        Before doing this, we recall that, in case the Poisson-Neumann problem \eqref{Neumann_problem} has the $W^{2,3}$-regularity, the solution $(u_m,v_m)$ of \eqref{problema_P_m_intro} have the regularity \eqref{regularidade_solucao_u_m_v_m_reg_W_2_3}. Hence, if we take the gradient in the $u_m$-equation of \eqref{problema_P_m_intro} we obtain
        \begin{equation}
          \begin{array}{rl}
            \nabla (\partial_ t u_m) - \nabla \Delta u_m & = - a_m(u_m) \nabla \Delta v_m - \nabla a_m(u_m) \Delta v_m \\
            & \quad - D^2 v_m \nabla a_m(u_m) - D^2 a_m(u_m) \nabla v_m.
          \end{array}
          \label{equacao_grad_u_m_2D}
        \end{equation}
        Now we test \eqref{equacao_grad_u_m_2D} by $\nabla\Delta u_m$. Using the $m$-uniform bounds obtained so far we can conclude that
        \begin{equation*}
          \Delta u_m \mbox{ is bounded in } L^{\infty}(0,\infty;L^2(\Omega)) \cap L^2(0,\infty;H^1(\Omega)).
        \end{equation*}
        Finally, if we look at \eqref{equacao_grad_u_m_2D} again we can also conclude that
        \begin{equation*}
          \partial_t u_m \mbox{ is bounded in } L^2(0,\infty;H^1(\Omega)).
        \end{equation*}
        This finishes the proof of Theorem \ref{teo_2D_existencia_unicidade_estimativas}.


\section*{Appendix}
\label{appendix}
\addcontentsline{toc}{section}{Appendix}
  
  \subsection*{A. Hypothesis \ref{hypothesis_density}}
  \addcontentsline{toc}{subsection}{A. Hypothesis \ref{hypothesis_density}}
  
    In the proof of Lemma \ref{lema_termo_fonte_final} we will need Hypothesis \ref{hypothesis_density} (see page \pageref{hypothesis_density}). Therefore, in order to show that this hypothesis is not too restrictive, we show that Hypothesis \ref{hypothesis_density} holds if the Poisson-Neumann problem \eqref{Neumann_problem}  has the $W^{3,p}$-regularity (see definition \ref{defi_regularidade_H_m} in page \pageref{defi_regularidade_H_m}), for $p > N$. According to \cite{Grisvard}, this is true if $\Gamma$ is at least $C^{2,1}$, for example. Up to our knowledge, the validity of Hypothesis \ref{hypothesis_density} in other domains of practical interest, such as polyhedra and polygons, is an open question.
    
    \begin{lema}
      Suppose that the Poisson-Neumann problem \eqref{Neumann_problem} has the $W^{3,p}$-regularity, for some $p > N$, and let $z \in H^2(\Omega)$ such that $\partial_{\eta} z \Big |_{\Gamma} = 0$. Then there is a sequence $\{ \rho_n \} \subset C^2(\overline{\Omega})$, with $\partial_{\eta}\rho_n \Big |_{\Gamma} = 0$, which converges to $z$ in $H^2(\Omega)$.
      \label{lema_densidade}
    \end{lema}
    \begin{proof}[\bf Proof]
      For any fixed $z \in H^2(\Omega)$ such that $\partial_{\eta} z \Big |_{\Gamma} = 0$, define $f = - \Delta z + z$. Note that $f \in L^2(\Omega)$ and $z \in H^2(\Omega)$ is the solution of
      \begin{equation}
        \left \{ \begin{array}{rl}
          - \Delta z + z & = f \\
          \partial_{\eta} z \Big |_{\Gamma} & = 0.
        \end{array} \right.
        \label{lema_densidade_equacao_z}
      \end{equation}
      
      Let $ \{ f_n \}_{n \in \mathbb{N}}$ be a sequence of mollifiers of $f$, that is, $f_n \in C^{\infty}_c(\Omega)$ and $f_n \to f$ in $L^2(\Omega)$ as $n \to \infty$. Then, for each fixed $n \in \mathbb{N}$, consider  the following problem: {\it Find }$\rho_n:\Omega \rightarrow \mathbb{R} $ {\it such that}
      \begin{equation}
        \left \{ \begin{array}{rl}
          - \Delta \rho_n + \rho_n & = f_n \\
          \partial_{\eta} \rho_n \Big |_{\Gamma} & = 0.
        \end{array} \right.
        \label{lema_densidade_equacao_rho_n}
      \end{equation}
      Considering the hypothesis that the Poisson-Neumann problem \eqref{Neumann_problem} has the $W^{3,p}$-regularity, for some $p > N$, we can conclude that, for each $n \in \mathbb{N}$, there is one, and only one, function $\rho_n \in W^{3,p}(\Omega) \subset C^2(\overline{\Omega})$ such that $\partial_{\eta} \rho_n \Big |_{\Gamma} = 0$ which solves problem \eqref{lema_densidade_equacao_rho_n}.
      
      Since the functions $z$ and $\rho_n$ solve the problems \eqref{lema_densidade_equacao_z} and \eqref{lema_densidade_equacao_rho_n}, respectively, the functions $z - \rho_n$ solve the problem
      \begin{equation}
        \left \{ \begin{array}{rl}
          - \Delta (z - \rho_n) + (z - \rho_n) & = (f - f_n) \\
          \partial_{\eta} (z - \rho_n) \Big |_{\Gamma} & = 0.
        \end{array} \right.
        \label{lema_densidade_equacao_z_menos_rho_n}
      \end{equation}
      Then, using $(z - \rho_n)$ and $- \Delta (z - \rho_n)$ as test functions in \eqref{lema_densidade_equacao_z_menos_rho_n} we can conclude that
      \begin{equation*}
        \norma{z - \rho_n}{L^2(\Omega)}, \norma{\nabla(z - \rho_n)}{L^2(\Omega)}, \norma{\Delta (z - \rho_n)}{L^2(\Omega)} \leq \norma{f - f_n}{L^2(\Omega)}.
      \end{equation*}
      Since $f_n \to f$ in $L^2(\Omega)$ as $n \to \infty$, the latter implies that $\rho_n \to z$ in $H^2(\Omega)$ as $n \to \infty$, finishing the proof.
    \end{proof}
    
  
  \subsection*{B. Proof of Lemma \ref{lema_termo_fonte_final}}
  \addcontentsline{toc}{subsection}{B. Proof of Lemma \ref{lema_termo_fonte_final}}
  
    Before proving Lemma \ref{lema_termo_fonte_final}, we must present some technical results.
    \begin{lema}
      Let $z: \Omega \rightarrow \mathbb{R}$ be a $C^2(\overline{\Omega})$ function such that $\partial_{\eta} z \Big |_{\Gamma} = 0$. Then
      \begin{align*}
        \int_{\Omega}{\norm{\Delta z}{}^2 \ dx} = \int_{\Omega}{\norm{D^2 z}{}^2 \ dx} - \dfrac{1}{2} \int_{\Gamma}{\nabla ( \norm{\nabla z}{}^2 ) \cdot  \eta \ d \Gamma}.
      \end{align*}
      \label{lema_dal_passo}
    \end{lema}
    \begin{proof}[\bf Proof]
      It suffices to prove the inequality for sufficiently regular functions and then pass to the limit. Integrating by parts we have
      \begin{align*}
        \int_{\Omega}{\norm{\Delta z}{}^2 \ dx} = - \int_{\Omega}{\nabla z \cdot \nabla \Delta z \ dx} & = \int_{\Omega}{\norm{D^2 z}{}^2 \ dx} - \int_{\Gamma}{[ (\nabla z)^T D^2 z ] \cdot \eta \ d\Gamma} \\
        & = \int_{\Omega}{\norm{D^2 z}{}^2 \ dx} - \frac{1}{2} \int_{\Gamma}{\nabla \norm{\nabla z}{}^2 \cdot \eta \ d\Gamma}.
      \end{align*}
    \end{proof}
    
    \begin{lema} \label{lema_winkler_0}
      There is a constant $C > 0$ such that
      \begin{enumerate}
        \item \begin{equation*}
          \int_{\Omega}{\norm{D^2 z}{}^2 \ dx} + \int_{\Omega}{\frac{\norm{\nabla z}{}^2}{z} \Delta z \ dx} = 4 \int_{\Omega}{z \norm{D^2 \sqrt{z}}{}^2 \ dx} + \frac{3}{4} \int_{\Omega}{\frac{\norm{\nabla z}{}^4}{z^2} \ dx},
        \end{equation*}
      \item \begin{equation*}
          \int_{\Omega}{\norm{D^2 z}{}^2 \ dx} \leq C \Big ( \int_{\Omega}{z \norm{D^2 \sqrt{z}}{}^2 \ dx} + \int_{\Omega}{\frac{\norm{\nabla z}{}^4}{z^2} \ dx} \Big ),
        \end{equation*}
      \end{enumerate}
      for all $z \in H^2(\Omega)$ such that $\partial_\eta z \Big |_{\Gamma} = 0$ and $z \geq \alpha$, for some $\alpha > 0$.
    \end{lema}
    \begin{proof}[\bf Proof]
      See lemma 3.3 of \cite{winkler2012global} for $1$. The inequality in $2$ is a direct consequence of the identity $1$.
    \end{proof}
    
    The next two results will allow us to estimate the boundary integral.
    
    \begin{lema} \label{lema_MARB_FGG}
      Let $\Gamma = \displaystyle\bigcup_{i=1}^m{\Gamma_i},$ each $\Gamma_i$ defined through a parametrization of one variable of $\mathbb{R}^3$ by the other two. Then there is $C > 0$ such that, for all $i$, one has
      \begin{equation*}
        \norm{\int_{ \Gamma_i}{\nabla \norm{\nabla z}{}^2 \cdot \eta \ d\Gamma_i}}{} \leq C \int_{\Gamma_i}{\norm{\nabla z}{}^2 d \Gamma_i},
      \end{equation*}
      for all $z \in C^2(\overline{\Omega})$ such that $\partial_{\eta} z \Big |_{\Gamma} = 0$.
    \end{lema}
    \begin{proof}[\bf Proof]
      See \cite{FGG_MARB_uniq_reg_Q_tensor}.
    \end{proof}
      
      \begin{lema}
        Let $\Omega$ be a Lipschitz domain. Then, for each $\delta > 0$, there is a $C(\delta) > 0$ such that
        \begin{equation*}
          \norma{\nabla z}{L^2(\Gamma)} \leq C(\delta) \norma{z}{L^2(\Omega)} + \delta \norma{D^2 z}{L^2(\Omega)}, \ \forall \ z \in H^2(\Omega).
        \end{equation*}
        \label{lema_dal_passo_2}
      \end{lema}
      \begin{proof}[\bfseries Proof]
        This is Lemma 2.4 in \cite{passo1998fourth}. It can be proved by contradiction as Lemma \ref{lema_desigualdade_poincare_media_nula}.
      \end{proof}
      
      Now we are in position of proving Lemma \ref{lema_termo_fonte_final}.
      
      \begin{proof}[\bf Proof of Lemma \ref{lema_termo_fonte_final}]
        We recall that, by hypothesis, $ z(x)\geq \alpha >0 , \ a.e. \ x \in \Omega$. Now we divide the proof in two main steps.
        
\
        
        \noindent {\bf STEP 1:}
         First of all, we are going to obtain the inequality
        \begin{equation}
          \begin{array}{rl}
            2 \D{\int_{\Omega}{\norm{\Delta z}{}^2 \ dx}} + 2 \D{\int_{\Omega}{\frac{\norm{\nabla z}{}^2}{z} \Delta z \ dx}} & \geq 8 \D{\int_{\Omega}}{z \norm{D^2 \sqrt{z}}{}^2 \ dx} + \dfrac{3}{2} \D{\int_{\Omega}}{\frac{\norm{\nabla z}{}^4}{z^2}} \\[12pt]
            & \ - C \delta \norma{D^2 z}{L^2(\Omega)}^2 - C(\delta) \norma{\nabla z}{L^2(\Omega)}^2.
          \end{array}
          \label{desigualdade_aux_lema_termo_fonte_final}
        \end{equation}
        
        In fact, from Hypothesis \ref{hypothesis_density}, we have the existence of a sequence $\{ z_j \}$ such that $z_j \in C^2(\overline{\Omega})$, the trace of the normal derivative of $z_j$ is zero and 
        $\displaystyle;
          \norma{z - z_j}{H^2(\Omega)} \rightarrow 0.
        $
        We can choose the sequence $\{ z_j \}$ to be such that
        \begin{equation*}
          \frac{\alpha}{2} \le  z_j(x) \ a.e. \ x \in \Omega, \forall j \geq j_0.
        \end{equation*}
        
        Now, applying Lemmas \ref{lema_dal_passo} and \ref{lema_winkler_0}-1, we have
        \begin{align*}
          2 \int_{\Omega}{\norm{\Delta z_j}{}^2 \ dx} + 2 \int_{\Omega}{\frac{\norm{\nabla z_j}{}^2}{z_j} \Delta z_j \ dx} & = 2 \int_{\Omega}{\norm{D^2 z_j}{}^2 \ dx} + 2 \int_{\Omega}{\frac{\norm{\nabla z_j}{}^2}{z_j} \Delta z_j \ dx} - \int_{\Gamma}{\nabla \norm{\nabla z_j}{}^2 \cdot \eta \ d\Gamma} \\
          & = 8 \int_{\Omega}{z_j \norm{D^2 \sqrt{z_j}}{}^2 \ dx} + \frac{3}{2} \int_{\Omega}{\frac{\norm{\nabla z_j}{}^4}{z_j^2}} - \int_{\Gamma}{\nabla \norm{\nabla z_j}{}^2 \cdot \eta \ d\Gamma }.
        \end{align*}
        Now we apply Lemma \ref{lema_MARB_FGG} and Lemma \ref{lema_dal_passo_2} (in this order) to obtain
        \begin{align*}
          2 \int_{\Omega}{\norm{\Delta z_j}{}^2 \ dx} & + 2 \int_{\Omega}{\frac{\norm{\nabla z_j}{}^2}{z_j} \Delta z_j \ dx} \geq 8 \int_{\Omega}{z_j \norm{D^2 \sqrt{z_j}}{}^2 \ dx} + \frac{3}{2} \int_{\Omega}{\frac{\norm{\nabla z_j}{}^4}{z_j^2}} - C \norma{\nabla z_j}{L^2(\Gamma)}^2 \\
          & \qquad \geq 8 \int_{\Omega}{z_j \norm{D^2 \sqrt{z_j}}{}^2 \ dx} + \frac{3}{2} \int_{\Omega}{\frac{\norm{\nabla z_j}{}^4}{z_j^2}} - C \delta \norma{D^2 z_j}{L^2(\Omega)}^2 - C(\delta) \norma{\nabla z_j}{L^2(\Omega)}^2,
        \end{align*}
        where $\delta > 0$ is supposed to be a sufficiently small number to be chosen later. Now we use the fact that $z_j \to z$ in $H^2(\Omega)$ as $j \to \infty$ and thus we obtain inequality \eqref{desigualdade_aux_lema_termo_fonte_final}, finishing the first step of the proof.
        
\
        
        \noindent {\bf STEP 2:}    
        Next, we apply Lemma \ref{lema_winkler_0}-2 to the right hand side of \eqref{desigualdade_aux_lema_termo_fonte_final} and choose $\delta > 0$ small enough, then there exist two constants $C_1,C_2>0$ such that
        \begin{align*}
          2 \int_{\Omega}{\norm{\Delta z}{}^2 \ dx} + 2 \int_{\Omega}{\frac{\norm{\nabla z}{}^2}{z} \Delta z \ dx} & \geq C_1 \Big ( \int_{\Omega}{\norm{D^2 z}{}^2 \ dx} + \int_{\Omega}{\frac{\norm{\nabla z}{}^4}{z^2} \ dx} \Big ) - C_2 \norma{\nabla z}{L^2(\Omega)}^2
        \end{align*}
        and the proof of Lemma \ref{lema_termo_fonte_final} is finished.
      \end{proof}

  \subsection*{Acknowledgments} 
   This work has been partially supported by Grant 
   PGC2018-098308-B-I00 (MCI/AEI/FEDER, UE).
   FGG has also been financed in part by the Grant
US-1381261 (US/JUNTA/FEDER, UE) and Grant P20-01120 (PAIDI/JUNTA/FEDER,
UE).
  ALCVF has also been financed in part by the Coordena\c c\~ao de Aperfei\c coamento de Pessoal de N\'{\i}vel Superior - Brasil (CAPES) - Finance Code 001.
  


  \addcontentsline{toc}{section}{References}
  \bibliographystyle{abnt_ord} 

\end{document}